\documentclass[10pt]{article}
\usepackage{amsmath,amssymb,amsthm,mathrsfs}
\usepackage[all]{xy}
\title{Independence-friendly cylindric set algebras}
\author{ALLEN L.~MANN}
\let	\implies		=	\rightarrow

\let	\iso			=	\cong

\let	\meet		=	\wedge
\let	\join			=	\vee

\let	\given		=	|

\let\abs=\envert

\let\norm=\enVert

\newcommand{\set}[1]{\{#1\}}

\newcommand{\setof}[2]{\{\,#1 \mid #2\,\}}

\newcommand{\tuple}[1]{(#1)}
\newcommand{\seq}[1]{\langle#1\rangle}

\newcommand{\seqof}[2]{\langle\,#1 \mid #2\,\rangle}
\newcommand{\concat}{^\frown}
\let	\tuple	=	\seq

\newcommand{\restrict}{\!\upharpoonright\!}

\DeclareMathOperator{\dom}{dom}

\DeclareMathOperator{\pr}{pr}

\DeclareMathOperator{\Sub}{Sub}

\DeclareMathOperator{\Cyls}{\mathfrak{Cs}}

\DeclareMathOperator{\Root}{Root}
\DeclareMathOperator{\Suit}{Suit}
\DeclareMathOperator{\DSuit}{DSuit}

\newcommand{\powerset}{\mathscr P}

\newcommand{\A}{\mathfrak{A}}
\newcommand{\B}{\mathfrak{B}}
\newcommand{\C}{\mathfrak{C}}

\newcommand{\NN}{\mathbb{N}}

\newcommand{\RR}{\mathbb{R}}

\newcommand{\book}[1]{\textit{#1}} 
\newcommand{\url}[1]{\texttt{#1}} 

\newtheorem{theorem}{Theorem}[section]
\newtheorem{corollary}[theorem]{Corollary}
\newtheorem{lemma}[theorem]{Lemma}
\newtheorem{proposition}[theorem]{Proposition}

\theoremstyle{definition}
\newtheorem*{definition}{Definition}%[section]
\newtheorem*{conjecture}{Conjecture}%[section]

\theoremstyle{remark}
\newtheorem*{example}{Example}

\newcommand{\thmref}[1]{Theorem~\ref{#1}}

\newcommand{\lemref}[1]{Lemma~\ref{#1}}
\newcommand{\propref}[1]{Proposition~\ref{#1}}
\newcommand{\corref}[1]{Corollary~\ref{#1}}

\newcommand{\hneg}{\sim\!}
\newcommand{\hand}[1]{\land\!_{/#1}\,}
\newcommand{\hor}[1]{\lor\!\!_{/#1}\,}

\newcommand{\hforall}[2]{\forall {#1}_{/#2}}
\newcommand{\hexists}[2]{\exists {#1}_{/#2}}
\newcommand{\modelt}{\models^+}
\newcommand{\modelf}{\models^-}
\newcommand{\modeltf}{\models^\pm}
\newcommand{\modelft}{\models^\mp}

\newcommand{\trump}[1]{\norm{#1}^+}
\newcommand{\cotrump}[1]{\norm{#1}^-}
\newcommand{\game}{G}
\newcommand{\abelard}{Ab\'elard}
\newcommand{\eloise}{Elo\"ise}
\newcommand{\toind}[1]{\underset{#1}{\to}}
\newcommand{\longtoind}[1]{\underset{#1}{\longrightarrow}}
\newcommand{\Cd}[1]{C^\partial_{#1}}

\newcommand{\n}[1]{{#1^{\hbox{\sm\char91}}}}
\renewcommand{\^}[2]{\,\!^{#1}\!{#2}}
\renewcommand{\cong}{\equiv}
\renewcommand{\cong}{\equiv}

\font\sm = cmsy5

\begin{document}
\maketitle
%% Abstract
\begin{abstract}
Independence-friendly logic is a conservative extension of first-order logic that has the same expressive power as existential second-order logic. In her Ph.D. thesis, Dechesne introduces a variant of independence-friendly logic called IFG logic. We attempt to algebraize IFG logic in the same way that Boolean algebra is the algebra of propositional logic and cylindric algebra is the algebra of first-order logic. 

We define independence-friendly cylindric set algebras and prove two main results. First, every independence-friendly cylindric set algebra over a structure has an underlying Kleene algebra. Moreover, the class of such underlying Kleene algebras generates the variety of all Kleene algebras. Hence the equational theory of the class of Kleene algebras that underly an independence-friendly cylindric set algebra is finitely axiomatizable. Second, every one-dimensional independence-friendly cylindric set algebra over a structure has an underlying monadic Kleene algebra. However, the class of such underlying monadic Kleene algebras does not generate the variety of all monadic Kleene algebras. Finally, we offer a conjecture about which subvariety of monadic Kleene algebras the class of such monadic Kleene algebras does generate.
\end{abstract}

%\Keywords{independence-friendly logic, cylindric algebra, De Morgan algebra, monadic De Morgan algebra}

%% Body

%\input{1-IFG_Logic}
%LaTeX by Allen L. Mann

% IFG Logic
\section{IFG Logic}
Imagine a sentence of first-order logic. Most likely, it has an initial block of quantifiers. If there are multiple quantifiers in the block, then some of the quantifiers are dependent on other quantifiers. In fact, there is only one possible dependence relation among the quantifiers: later quantifiers depend on prior quantifiers. Hence the dependence relation is a linear order. The first attempt to allow dependence relations other than the usual one on the quantifiers of a first-order sentence was made by Henkin\index{Henkin, L.} \cite{Henkin:1961}. He allowed the dependence relation to be a partial order, rather than a linear order. For example, in the Henkin sentence
\[
\left(
\begin{array}{cc}
\forall x & \exists y \\
\forall z & \exists w 
\end{array}
\right)
\phi(x,y,z,w)
\]
the variable $y$ depends only on $x$, while $w$ depends only on $z$. It is a result due to Ehrenfeucht\index{Ehrenfeucht, A.} \cite{Henkin:1961} that the above quantifier is not definable in ordinary first-order logic. Later it was shown independently by Enderton\index{Enderton, H.~B.} \cite{Enderton:1970} and Walkoe\index{Walkoe Jr., W.~J.} \cite {Walkoe:1970} that first-order logic with these ``branching\index{branching quantifiers} quantifiers'' has the same expressive power as existential second-order logic.

Independence-friendly logic\index{independence-friendly logic|(} (IF logic\index{IF logic|see{independence-friendly logic}}) was introduced by Hintikka\index{Hintikka, J.} and Sandu\index{Sandu, G.} \cite{Hintikka:1989} as a way to allow arbitrary dependence relations between the quantifiers of a first-order sentence. Independence-friendly logic is a conservative extension of ordinary first-order logic in that every ordinary first-order sentence has an independence-friendly counterpart that is true in exactly the same models. However, the ability to specify arbitrary dependence relations between the quantifiers (and even the connectives) means that there are far more sentences than before. It turns out that a given independence-friendly sentence has the same expressive power as a pair of contrary existential second-order sentences. In \book{The Principles of Mathematics Revisited} \cite{Hintikka:1996}, Hintikka\index{Hintikka, J.} argues that IF logic is the correct first-order logic and advocates for its adoption as the foundation of mathematics.

In IF logic the truth of a sentence (or formula) is defined in terms of games. Let $\phi$ be a sentence, and let $\A$ be a suitable structure. Informally, the semantic game\index{semantic game} $\game(\A, \phi)$\index{$\game(\A, \phi)$} is played between two players, player 0 (\abelard) and player 1 (\eloise). \eloise's goal is to verify the sentence $\phi$ in the structure $\A$, while \abelard's goal is to falsify it. A familiar example from calculus is the definition of continuity:
\[
\forall x\forall\varepsilon (\varepsilon > 0 \implies \exists\delta (\delta >0 \land \forall y (\abs{x - y} < \delta \implies \abs{f(x) - f(y)} < \varepsilon))).
\]
Given a function $f$, \abelard\ picks an $x$ and an $\varepsilon$. If \(\varepsilon \leq 0\), \eloise\ wins. If \(\varepsilon > 0\), \eloise\ chooses a $\delta$. If \(\delta \leq 0\), \abelard\ wins. If \(\delta > 0\), \abelard\ chooses a $y$. If \(\abs{x - y} \geq \delta\), \eloise\ wins. Otherwise, \eloise\ wins if \(\abs{f(x) - f(y)} < \varepsilon\), and \abelard\ wins if \(\abs{f(x) - f(y)} \geq \varepsilon\).

Notice that \eloise\ winning a single play of the game is not sufficient to guarantee that $f$ is continuous. Likewise, \abelard\ winning a single play is not sufficient to show that $f$ is not continuous. What is required for $f$ to be continuous is for \eloise\ to have a way to win every play of the game, given correct play by both players---that is, \eloise\ must have a winning strategy. Dually, $f$ is not continuous if and only if \abelard\ has a winning strategy. By extension, a first-order sentence is true in a model $\A$ if and only if \eloise\ has a winning strategy for the game \(\game(\A, \phi)\), and $\phi$ is false in $\A$ if and only if \abelard\ has a winning strategy for \(\game(\A, \phi)\).

In IF logic a sentence is defined to be true in a model if \eloise\ has a winning strategy for the corresponding semantic game. Dually, a sentence is false in a model if \abelard\ has a winning strategy. For any ordinary first-order sentence $\phi$ the semantic game $\game(\A, \phi)\) is a two-player, win-loss game of perfect information. Hence Zermelo's\index{Zermelo, E.} theorem tells us that one of the two players must have a winning strategy. Thus every ordinary first-order sentence is either true of false. That is, one can prove the principle of bivalence\index{principle of bivalence} for ordinary first-order sentences from the game-theoretical definition of truth.

The extension beyond ordinary first-order logic comes from the fact that in IF logic one can write down a sentence whose corresponding semantic game is not a game of perfect information by restricting the information available to the existential player. For such a sentence it is no longer the case that one of the two players must have a winning strategy. Thus it is possible for a sentence in IF logic to be neither true nor false. For example, consider the sentence 
\[
\forall{x}\hexists{y}{x}(x = y).
\]
First \abelard\ chooses an element of the universe and assigns it as the value of the variable $x$. Then \eloise\ chooses an element of the universe and assigns it to $y$, but \eloise\ must make her choice without knowing the value of $x$. If \eloise\ correctly guesses which element \abelard\ chose, she wins. Otherwise \abelard\ wins. It should be clear that in any structure with at least two elements, \eloise\ does not have a winning strategy. Therefore the sentence is not true. But neither does \abelard\ have a winning strategy because there is always the possibility that \eloise\ will guess correctly. Therefore the sentence is not false. It is worth noting that whether or not a sentence is undetermined depends on the structure in which the semantic game is played. For example, the sentence
\[
\forall{x}\hexists{y}{x}(y \leq x)
\]
is neither true nor false in $\RR$, but it is true in $\NN$ because \eloise\ can always choose 0. \index{independence-friendly logic|)}

In her Ph.D. thesis \cite{Dechesne:2005} Dechesne\index{Dechesne, F.} provides a rigorous mathematical foundation for an extension of IF logic in which one allows the information available to \eloise\ and \abelard\ to be restricted. She calls her extension IFG logic\index{IFG logic|(} (for ``generalized independence-friendly logic''). Her thesis is the basis for our work. With the possible exception of \propref{true for any valuation -> true for all}, none of the results in the present section are original.

Instead of focusing on which quantifiers are independent of one another, it might seem more natural to focus on which quantifiers depend on which other quantifiers. V\"a\"an\"anen\index{V\"a\"an\"anen, J.} does exactly that in his recent book \book{Dependence Logic} \cite{Vaananen:2007}.

Cylindric algebra\index{cylindric algebra} is the algebra of ordinary first-order logic in the same way that Boolean algebra\index{Boolean algebra} is the algebra of ordinary propositional logic. Cylindric algebra was first studied by Henkin\index{Henkin, L.}, Monk\index{Monk, J.~D.}, and Tarski\index{Tarski, A.} \cite{Henkin:1971, Henkin:1985}. The goal of this paper is to algebraize IFG logic in the same spirit as cylindric algebra. Our algebraization will depend heavily on the compositional semantics for independence-friendly logic put forth by Hodges\index{Hodges, W.} \cite{Hodges:1997a, Hodges:1997b}.

% Syntax
\subsection{Syntax}

In regular first-order logic, a formula $\phi$ is a string of symbols that satisfies certain conditions, and a variable $x$ is said to occur in $\phi$ if the symbol ``$x$'' appears in the string. In our version of IFG logic\index{IFG logic|)}, each formula will be a pair $\tuple{\phi, X}$ where $\phi$ is a formula in the standard sense (a string of symbols satisfying certain conditions) and $X$ is a finite set of variables. A variable is said to occur in $\tuple{\phi, X}$ if and only if it belongs to $X$. We require that any variable that appears in $\phi$ must belong to $X$. However, we will allow variables that do not appear in $\phi$ to belong to $X$. Thus every variable that appears in $\phi$ must occur in $\tuple{\phi, X}$, but a variable may occur in $\tuple{\phi, X}$ without appearing in $\phi$.

\begin{definition}
Given a first-order signature $\sigma$, an \emph{atomic IFG-formula}\index{IFG-formula!atomic|mainidx} is a pair $\tuple{\phi, X}$ where $\phi$ is an atomic first-order formula and $X$ is a finite set of variables that includes every variable that appears in $\phi$.
\end{definition}

\begin{definition}
Given a first-order signature $\sigma$, the language $\mathscr L_\mathrm{IFG}^\sigma$\index{$\mathscr L_\mathrm{IFG}^\sigma$|mainidx}\index{IFG-formula|(} is the smallest set of formulas such that:
\begin{enumerate}
	\item Every atomic IFG-formula is in $\mathscr L_\mathrm{IFG}^\sigma$.
	\item If $\tuple{\phi, Y}$ is in $\mathscr L_\mathrm{IFG}^\sigma$ and \(Y \subseteq X\), then $\tuple{\phi, X}$ is in $\mathscr L_\mathrm{IFG}^\sigma$.
	\item If $\tuple{\phi, X}$ is in $\mathscr L_\mathrm{IFG}^\sigma$, then $\tuple{\hneg\phi, X}$ is in $\mathscr L_\mathrm{IFG}^\sigma$.
	\item If $\tuple{\phi, X}$ and $\tuple{\psi, X}$ are in $\mathscr L_\mathrm{IFG}^\sigma$, and \(Y \subseteq X\), then $\tuple{\phi \hor{Y} \psi, X}$ is in $\mathscr L_\mathrm{IFG}^\sigma$.
	\item If $\tuple{\phi, X}$ is in $\mathscr L_\mathrm{IFG}^\sigma$, \(x \in X\), and \(Y \subseteq X\), then $\tuple{\hexists{x}{Y}\phi, X}$ is in $\mathscr L_\mathrm{IFG}^\sigma$.
\end{enumerate}
Above $X$ and $Y$ are finite sets of variables.
\end{definition}

From now on we will make certain assumptions about IFG-formulas that will allow us to simplify our notation. First, we will assume that the set of variables of $\mathscr L_\mathrm{IFG}^\sigma$ is $\setof{v_n}{n \in \omega}$. Second, since it does not matter much which particular variables appear in a formula, we will assume that variables with smaller indices are used before variables with larger indices. More precisely, if $\tuple{\phi, X}$ is a formula, \(v_j \in X\), and \(i \leq j\), then \(v_i \in X\). By abuse of notation, if $\tuple{\phi, X}$ is a formula and \(\abs X = N\), then we will say that $\phi$ has $N$ variables and write $\phi$ for $\tuple{\phi, X}$. As a shorthand, we will call $\phi$ an IFG$_N$-formula\index{IFG$_N$-formula}. Let \(\mathscr L^\sigma_{\mathrm{IFG}_N} = \setof{\phi \in \mathscr L^\sigma_{\mathrm{IFG}}}{\phi \text{ has $N$ variables}}\)\index{$\mathscr L^\sigma_{\mathrm{IFG}_N}$}. Third, sometimes we will write $\phi \hor{J} \psi$ instead of $\phi \hor{Y} \psi$ and $\hexists{v_n}{J}\phi$ instead of $\hexists{v_n}{Y}\phi$, where \(J = \setof{j}{v_j \in Y}\). Finally, we will use \(\phi \hand{J} \psi\) to abbreviate \(\hneg(\hneg\phi \hor{J} \hneg\psi)\) and \(\hforall{v_n}{J} \phi\) to abbreviate \(\hneg\hexists{v_n}{J}\hneg\phi\).

\begin{definition} % subformula tree
Let $\phi$ be an IFG-formula. The \emph{subformula tree} of $\phi$, denoted $\Sub(\phi)$\index{$\Sub(\phi)$ \quad subformula tree of $\phi$|mainidx}, is the smallest tree satisfying the following conditions.
\begin{enumerate}
	\item \(\tuple{\emptyset, \phi} \in \Sub(\phi)\).
	\item If \(\tuple{s, \hneg\psi} \in \Sub(\phi)\), then \(\tuple{s\concat 0, \psi} \in \Sub(\phi)\).
	\item If \(\tuple{s, \psi_1 \hor{J} \psi_2} \in \Sub(\phi)\), then \(\tuple{s\concat 1, \psi_1} \in \Sub(\phi)\) and \(\tuple{s\concat 2, \psi_2} \in \Sub(\phi)\).
	\item If \(\tuple{s, \hexists{v_n}{J}\psi} \in \Sub(\phi)\), then \(\tuple{s\concat 3, \psi} \in \Sub(\phi)\).
\end{enumerate}
For every \(\tuple{s, \psi} \in \Sub(\phi)\), \(\tuple{s, \psi} \in \Sub^+(\phi)\)\index{$\Sub(\phi)$ \quad subformula tree of $\phi$!$\Sub^+(\phi)$ \quad positive subformula tree of $\phi$|mainidx} if $s$ contains an even number of 0s, and \(\tuple{s, \psi} \in \Sub^-(\phi)\)\index{$\Sub(\phi)$ \quad subformula tree of $\phi$!$\Sub^-(\phi)$ \quad negative subformula tree of $\phi$|mainidx} if $s$ contains an odd number of 0s.
\end{definition}

From now on, we will assume that all subformulas are indexed by their position in the subformula tree. This will allow us to distinguish between multiple instances of the same formula that may occur as subformulas of $\phi$. For example, if $\phi$ is \(v_0 = v_1 \hor{v_0} v_0 = v_1\) we will distinguish between the left and right disjuncts. \index{IFG-formula|)}

% Game Semantics
\subsection{Game semantics}

Dechesne\index{Dechesne, F.} defines her semantic games in extensive form\index{games in extensive form} \cite{Dechesne:2005}, which is standard practice in game theory. It turns out that games in extensive form are more general than is necessary for our purposes. We modify her definition in order to focus on those aspects of the games that are relevant to the present discussion.  

\begin{definition} % suitable structure
Given a first-order signature $\sigma$ and a formula \(\tuple{\phi, X} \in \mathscr L^\sigma_\mathrm{IFG}\), a structure $\A$ is called \emph{suitable}\index{suitable structure|mainidx} for $\phi$ if $\A$ has an interpretation for every non-logical symbol in $\sigma$. 
\end{definition}

\begin{definition}
If $\tuple{\phi, X}$ is a formula and $\A$ is a suitable structure, then a \emph{valuation}\index{valuation|mainidx} for $\tuple{\phi, X}$ over $\A$ is a function from $X$ to $A$. Since we are assuming that $X$ has the form $\set{v_0, \ldots, v_{N-1}}$, we will identify valuations with sequences of individuals in $A$, denoted \(\vec a \in \^NA\). A set of valuations \(V \subseteq \^NA\) is called a \emph{team}\index{team|mainidx}.
\end{definition}

\begin{definition} % agree outside of $J$
Let \(\vec a, \vec b \in\, \^NA\) be two valuations, and let \(J \subseteq N\). We say that $\vec a$ and $\vec b$ \emph{agree outside of $J$}, denoted \(\vec a \approx_J \vec b\)\index{\(\vec a \approx_J \vec b\)\quad $\vec a$ and $\vec b$ agree outside of $J$|mainidx}, if \(\vec a \restrict (N\setminus J) = \vec b \restrict(N\setminus J)\).
\end{definition}

Note that $\approx_J$ is an equivalence relation on $\^NA$. Also note that $\approx_\emptyset$ is the identity relation and $\approx_N$ is the total relation on $\^NA$.

\begin{lemma} % \approx_J implies \approx_K
Let \(\vec a, \vec b \in\, \^NA\), and let \(J \subseteq K \subseteq N\). Then \(\vec a \approx_J \vec b\) implies \(\vec a \approx_K \vec b\).
\end{lemma}
\begin{proof}
If \(\vec a \approx_J \vec b\), then \(\vec a \restrict (N\setminus J) = \vec b \restrict (N\setminus J)\), which implies \(\vec a \restrict (N\setminus K) = \vec b \restrict (N\setminus K)\). Hence \(\vec a \approx_K \vec b\).
\end{proof}

\begin{definition} % $n$-variant
If \(\vec a \in\, \^NA\), \(b \in A\), and \(n < N\), define $\vec a(n:b)$\index{$\vec a(n:b)$|mainidx} to be the valuation that is like $\vec a$ except that $v_n$ is assigned the value $b$ instead of $a_n$. In other words,
\[
\vec a(n:b) = \vec a\restrict(N\setminus \set{n}) \cup \set{\tuple{n,b}}.
\]
We call $\vec a(n:b)$ an \emph{$n$-variant of $\vec a$}\index{variant|mainidx}.
\end{definition}

\begin{definition} % $n$-variation
If \(V \subseteq\, \^NA\) is a team and \(b \in A\), define 
\[
V(n:b) = \setof{\vec a(n:b)}{\vec a \in V}.\index{$V(n:b)$|mainidx}
\]

Furthermore, if \(B \subseteq A\) define
\[
V(n:B) = \setof{\vec a(n:b)}{\vec a \in V,\ b \in B}. \index{$V(n:B)$|mainidx}
\]
A set \(V' \subseteq V(n:A)\) is called an \emph{$n$-variation of $V$} if for every \(\vec a \in V\) there is at least one $n$-variant of $\vec a$ in $V'$.
Finally if \(f\colon V \to A\), and \(V' \subseteq V\), define the \emph{$n$-variation\index{variation|mainidx} of $V'$ by $f$} to be
\[
V'(n:f) = \setof{\vec a(n:f(\vec a))}{\vec a \in V'}.\index{$V(n:f)$|mainidx}
\]
\end{definition}

\begin{definition} % semantic game
Let $\phi$ be a formula with $N$ variables, let $\A$ be a suitable structure, and let \(V \subseteq \^NA\) be a team. The \emph{semantic game}\index{semantic game|mainidx} $\game(\A, \phi, V)$\index{$\game(\A, \phi, V)$|mainidx} is defined as follows. A \emph{position}\index{position|mainidx} of the game is a triple $\tuple{\psi, \vec b, \varepsilon}$, where $\psi$ is a subformula of $\phi$, \(\vec b \in \^NA\), and \(\varepsilon \in \set{0,1}\). A \emph{terminal position}\index{position!terminal|mainidx} is a position in which $\psi$ is an atomic formula. A \emph{play}\index{play|mainidx} of the game is a sequence of positions $\seq{p_0, \ldots, p_q}$ that satisfies the following conditions.
\begin{enumerate}
	\item The initial position \(p_0 = \tuple{\phi, \vec a, 1}\), where \(\vec a \in V\).
	\item If \(p_k = \tuple{\hneg\psi, \vec b, \varepsilon}\), then \(p_{k+1} = \tuple{\psi, \vec b, 1 - \varepsilon}\).
	\item If \(p_k = \tuple{\psi_1 \hor{J} \psi_2,\, \vec b,\, \varepsilon}\), then \(p_{k+1} = \tuple{\psi_1, \vec b, \varepsilon}\) or \(p_{k+1} = \tuple{\psi_2, \vec b, \varepsilon}\).
	\item If \(p_k = \tuple{\hexists{v_n}{J}\psi,\, \vec b,\, \varepsilon}\), then \(p_{k+1} = \tuple{\psi,\, \vec b(n:c),\, \varepsilon}\) for some \(c \in A\).
	\item The final position $p_q$ is a terminal position, and $p_q$ is the only terminal position in the play.
\end{enumerate}
A \emph{partial play}\index{play!partial|mainidx} $\seq{p_0, \ldots, p_\ell}$ is an initial segment of a play. A partial play that is not a play is called a \emph{proper partial play}.

For a given play of the game with final position \(p_q = \tuple{\psi, \vec b, \varepsilon}\) where $\psi$ is an atomic formula, player $\varepsilon$ \emph{wins}\index{win|mainidx} if \(\A \models \psi[\vec b]\), and player $1 - \varepsilon$ \emph{wins} if \(\A \not\models \psi[\vec b]\). In a given position \(\tuple{\psi, \vec b, \varepsilon}\), player $\varepsilon$ is called the \emph{verifier}\index{verifier|mainidx} and player $1 - \varepsilon$ is called the \emph{falsifier}\index{falsifier|mainidx}. The game $\game(\A, \phi, V)$ is the set of all possible plays. We will use $\game(\A, \phi)$\index{$\game(\A, \phi)$|mainidx} to abbreviate $\game(\A, \phi, \^NA)$.
\end{definition}

\begin{definition} % strategy
A \emph{strategy}\index{strategy|mainidx} for player $\varepsilon$ for the game $\game(\A, \phi, V)$ is a function $S$ from the set of all non-terminal positions of the game in which player $\varepsilon$ is the verifier to the set of all positions. A strategy is \emph{legal}\index{strategy!legal|mainidx} if for every proper partial play $\seq{p_0, \ldots, p_\ell}$ where player $\varepsilon$ is the verifier in $p_\ell$, the sequence $\seq{p_0, \ldots, p_\ell, p_{\ell+1}}$ is a partial play, where \(p_{\ell+1} = S(p_\ell)\), and 
\begin{enumerate}
	\item if \(p_\ell = \tuple{\psi_1 \hor{J} \psi_2,\, \vec a,\, \varepsilon}\) and \(p'_\ell = \tuple{\psi_1 \hor{J} \psi_2,\, \vec b,\, \varepsilon}\), where \(\vec a \approx_J \vec b\), then \(S(p_\ell) = S(p'_\ell)\);
	\item if \(p_\ell = \tuple{\hexists{v_n}{J}\psi,\, \vec a,\, \varepsilon}\) and \(p'_\ell = \tuple{\hexists{v_n}{J}\psi,\, \vec b,\, \varepsilon}\), where \(\vec a \approx_J \vec b\), then \(S(p_\ell) = S(p'_\ell)\).
\end{enumerate}
Given a play \(p = \seq{p_0, \ldots, p_q}\) and a strategy $S$ for player $\varepsilon$, player $\varepsilon$ is said to \emph{follow $S$ in $p$} if for every non-terminal position $p_k$ in which player $\varepsilon$ is the verifier, \(p_{k+l} = S(p_k)\).
A strategy for player $\varepsilon$ is \emph{winning}\index{strategy!winning|mainidx} if it is legal and player $\varepsilon$ wins every play in which he or she follows $S$.
\end{definition}

Observe that if $V$ is empty the game $\game(\A, \phi, V)$ has no positions nor plays. Hence the empty strategy \(\emptyset\colon \emptyset \to \emptyset\) is a winning strategy for both players. As we will see later when we define a Tarski-style satisfaction relation for IFG-formulas, this apparent defect is actually a feature.

\begin{definition} % true and false
We say that $\phi$ is \emph{true in $\A$ relative to $V$} if player 1 has a winning strategy for the semantic game $\game(\A, \phi, V)$, and that $\phi$ is \emph{false in $\A$ relative to $V$} if player 0 has a winning strategy for $\game(\A, \phi, V)$. In the first case, we call $V$ a \emph{winning team}\index{team!winning|mainidx} (or \emph{trump}\index{trump|see{team, winning}}) for $\phi$ in $\A$. In the second case, we call $V$ a \emph{losing team}\index{team!losing|mainidx} (or \emph{cotrump}\index{cotrump|see{team, losing}}) for $\phi$ in $\A$. We say that $\phi$ is \emph{true}\index{true|mainidx} in $\A$ if it is true in $\A$ relative to $\^NA$, and that $\phi$ is \emph{false}\index{false|mainidx} in $\A$ if it is false relative to $\^NA$.
\end{definition}

Thus $\phi$ is true in $\A$ if and only if player 1 has a winning strategy for the game $\game(\A, \phi)$, and $\phi$ is false in $\A$ if and only if player 0 has a winning strategy for $\game(\A, \phi)$. It is important to realize that restricting the information available to the players at different moves does not change the set of possible plays of the game $\game(\A, \phi, V)$. It only restricts the strategies the players are allowed to use. 

\begin{definition} % dual position, dual strategy
Let \(p_k = \tuple{\psi, \vec b, \varepsilon}\) be a position of the game \(\game(\A, \phi, V)\). The \emph{dual position}\index{position!dual|mainidx} of $p_k$ is \(\widetilde p_k = \tuple{\psi, \vec b, 1 - \varepsilon}\). If $S_1$ is a strategy for player 1 for $\game(\A, \phi, V)$, the \emph{dual strategy}\index{strategy!dual|mainidx} $\widetilde S_1$ for player 0 for $\game(\A, \hneg\phi, V)$ is defined by \(\widetilde S_1(\widetilde p_k) = S_1(p_k)\) for all \(p_k \in \dom(S_1)\). If $S_0$ is a strategy for player 0 for $\game(\A, \phi, V)$, the \emph{dual strategy} $\widetilde S_0$ for player 1 for $\game(\A, \hneg\phi, V)$ is defined by \(\widetilde S_0(\tuple{\hneg\phi, \vec a, 1}) = S_0(\tuple{\phi, \vec a, 0})\) for all \(\vec a \in V\), and \(\widetilde S_0(\widetilde p_k) = S_0(p_k)\) for all \(p_k \in \dom(S_0)\). 
\end{definition}

\begin{lemma} % dual strategies
\label{dual strategies}
Let $\phi$ be an IFG$_N$-formula, let $\A$ be a suitable structure, and let \(V \subseteq \^NA\). Then $S$ is a winning strategy for player $\varepsilon$ for the game $\game(\A, \phi, V)$ if and only if $\widetilde S$ is a winning strategy for player $1 - \varepsilon$ for $\game(\A, \hneg\phi, V)$.
\end{lemma}

\begin{proof}
Suppose $S_1$ is a winning strategy for player $1$ for \(\game(\A, \phi, V)\). Then $\widetilde S_1$ is a legal strategy for player 0 for $\game(\A, {\hneg\phi}, V)$. To show that $\widetilde S_1$ is a winning strategy, let \(\widetilde p = \seq{\widetilde p_0, \widetilde p_1, \ldots, \widetilde p_q}\) be a play of \(\game(\A, {\hneg\phi}, V)\) in which player 0 follows $\widetilde S_1$, and let \(\widetilde p_q = \tuple{\psi, \vec b, 1- \varepsilon}\). Then the corresponding play \(p = \seq{p_1, \ldots, p_q}\) of \(\game(\A, \phi, V)\) is a play in which player 1 follows $S_1$ and \(p_q = \tuple{\psi, \vec b, \varepsilon}\). By hypothesis, player 1 wins $p$. Hence \(\A \models \psi[\vec b]\) if \(\varepsilon = 1\) and \(\A \not\models \psi[\vec b]\) if \(\varepsilon = 0\). Thus \(\A \models \psi[\vec b]\) if \(1 - \varepsilon = 0\) and \(\A \not\models \psi[\vec b]\) if \(1 - \varepsilon = 1\). In either case, player 0 wins $\widetilde p$.

Conversely, suppose $\widetilde S_1$ is a winning strategy for player 0 for \(\game(\A, {\hneg\phi}, V)\). Then $S_1$ is a legal strategy for player 1 for $\game(\A, \psi, V)$. To show $S_1$ is a winning strategy, let \(p = \seq{p_1, \ldots, p_q}\) be a play of $\game(\A, \psi, V)$ in which player 1 follows $S_1$, where \(p_1 = \tuple{\phi, \vec a, 1}\) and \(p_q = \tuple{\psi, \vec b, \varepsilon}\). Then \(\widetilde p = \seq{\widetilde p_0, \widetilde p_1, \ldots, \widetilde p_q}\) is a play of \(\game(\A, \phi, V)\), where \(\widetilde p_0 = \tuple{{\hneg\phi}, \vec a, 1}\) and \(\widetilde p_q = \tuple{\psi, \vec b, 1 - \varepsilon}\). By hypothesis, player 0 wins $\widetilde p$. Hence \(\A \models \psi[\vec b]\) if \(1 - \varepsilon = 0\) and \(\A \not\models \psi[\vec b]\) if \(1 - \varepsilon = 1\). Thus \(\A \models \psi[\vec b]\) if \(\varepsilon = 1\) and \(\A \not\models \psi[\vec b]\) if \(\varepsilon = 0\). In either case, player 1 wins.

Similarly, $S_0$ is a winning strategy for player 0 for \(\game(\A, \phi, V)\) if and only if $\widetilde S_0$ is a winning strategy for player 1 for \(\game(\A, {\hneg\phi}, V)\). 
\end{proof}

\begin{proposition} % negation
\label{negation}
Let $\phi$ be an IFG$_N$-formula, let $\A$ be a suitable structure, and let \(V \subseteq \^NA\). Then $\phi$ is true in $\A$ relative to $V$ if and only if $\hneg\phi$ is false in $\A$ relative to $V$, and vice versa.
\end{proposition}
\begin{proof}
By the previous lemma, \eloise\ has a winning strategy for \(\game(\A, \phi, V)\) if and only if \abelard\ has a winning strategy for \(\game(\A, {\hneg\phi}, V)\), and \eloise\ has a winning strategy for \(\game(\A, \hneg\phi, V)\) if and only if \abelard\ has a winning strategy for \(\game(\A, \phi, V)\).
\end{proof}

\begin{proposition} % game double negation
\label{game double negation}
Let $\phi$ be an IFG$_N$-formula, let $\A$ be a suitable structure, and let \(V \subseteq \^NA\). Then $\phi$ is true in $\A$ relative to $V$ if and only if $\hneg(\hneg\phi)$ is true in $\A$ relative to $V$.
\end{proposition}
\begin{proof}
By \propref{negation}, $\phi$ is true in $\A$ relative to $V$ if and only if $\hneg\phi$ is false in $\A$ relative to $V$ if and only if $\hneg(\hneg\phi)$ is true in $\A$ relative to $V$.
\end{proof}

\begin{proposition}
\label{true for any valuation -> true for all}
Let $\phi$ be an IFG$_N$-sentence. If $\phi$ is true in $\A$ relative to some nonempty \(V \subseteq \^NA\), then $\phi$ is true in $\A$.
\end{proposition}

\begin{proof}
Suppose $S$ is a winning strategy for \eloise\ for \(\game(\A, \phi, V)\). We will construct a winning strategy for \eloise\ for \(\game(\A, \phi, \^NA)\). To do so, we will need to keep track of which variables the players have had the opportunity to modify during the play of the game. For each subformula $\psi$ of $\phi$, define a set of indices $J_\psi$ of those variables of $\psi$ that have been unbound:
\begin{enumerate}
	\item \(J_\phi = \emptyset\).
	\item If $\psi$ is a subformula of $\phi$ of the form $\hneg\chi$, then \(J_\chi = J_\psi\).
	\item If $\psi$ is a subformula of $\phi$ of the form $\chi_1 \hor{K} \chi_2$, then \(J_{\chi_1} = J_\psi\) and \(J_{\chi_2} = J_\psi\).
	\item If $\psi$ is a subformula of $\phi$ of the form $\hexists{v_n}{K}\chi\), then \(J_\chi = J_\psi \cup \set{n}\).
\end{enumerate}

Fix \(\vec a \in V\). For every position \(p = \tuple{\psi, \vec b, \varepsilon}\) of \(\game(\A, \phi, \^NA)\) define 
\[
f(p) = (\vec a \restrict N \setminus J_\psi) \cup (\vec b \restrict J_\psi)
\]
and \(F(p) = \tuple{\psi, f(p), \varepsilon}\). Observe that $F(p)$ is a position of \(\game(\A, \phi, V)\). Define a strategy $T$ for \eloise\ for \(\game(\A, \phi, \^NA)\) as follows:
\begin{enumerate}
	\item If \(p = \tuple{\hneg\chi, \vec b, \varepsilon}\), then \(T(p) = \tuple{\chi, \vec b, 1 - \varepsilon}\).
	\item If \(p = \tuple{\chi_1 \hor{K} \chi_2, \vec b, \varepsilon}\) and \(S(F(p)) = \tuple{\chi_i, f(p), \varepsilon}\), then \(T(p) = \tuple{\chi_i, \vec b, \varepsilon}\).
	\item If \(p = \tuple{\hexists{v_n}{K}\chi, \vec b, \varepsilon}\) and \(S(F(p)) = \tuple{\chi, f(p)(n:c), \varepsilon}\), then \(T(p) = \tuple{\chi, \vec b(n:c), \varepsilon}\).
\end{enumerate}
To show that $T$ is a legal strategy, it suffices to observe that if \(p = \tuple{\psi, \vec b, \varepsilon}\) and \(p' = \tuple{\psi, \vec b', \varepsilon}\), where \(\vec b \approx_K \vec b'\), then 
\[
f(p) = (\vec a \restrict N \setminus J_\psi) \cup (\vec b \restrict J_\psi) \approx_K (\vec a \restrict N \setminus J_\psi) \cup (\vec b' \restrict J_\psi) = f(p').
\]
Hence \(S(F(p)) = S(F(p'))\) because $S$ is a legal strategy. Thus \(T(p) = T(p')\).

To show that $T$ is a winning strategy, let \(\tuple{p_0, \ldots, p_q}\) be a play of \(\game(\A, \phi, \^NA)\) in which \eloise\ follows $T$. Then \(\tuple{F(p_0), \ldots, F(p_q)}\) is a play of \(\game(\A, \phi, V)\) in which \eloise\ follows $S$. Let \(p_q = \tuple{\psi, \vec b, \varepsilon}\) and \(F(p_q) = \tuple{\psi, f(p_q), \varepsilon}\). Note that $\vec b$ and $f(p_q)$ agree on the free variables of $\psi$ because $\phi$ was a sentence. Thus \(\A \models \psi[\vec b]\) if and only if \(\A \models \psi[f(p_q)]\) if and only if \(\varepsilon = 1\) because $S$ is a winning strategy for \eloise\ for \(\game(\A, \phi, V)\).
\end{proof}

 \renewcommand{\theenumi}{\alph{enumi}}
 \renewcommand{\labelenumi}{(\theenumi)}

\begin{definition} % reachable position
A position $\tuple{\psi, \vec b, \varepsilon}$ of the game \(\game(\A, \phi, V)\) is \emph{reachable}\index{position!reachable|mainidx} if it occurs in some play of \(\game(\A, \phi, V)\). Otherwise, it is \emph{unreachable}.
\end{definition}

For example, let $\phi$ be \(\exists v_1(v_0 = v_1)\), and let $\A$ be the equality structure with universe $\set{0,1}$. Then \(\tuple{\phi, 00, 0}\) and \(\tuple{v_0 = v_1, 10, 1}\) are both unreachable positions of the game \(\game(\A, \phi, \set{00, 01})\). The position \(\tuple{\phi, 00, 0}\) is unreachable because player 1 is always the initial verifier;  the position \(\tuple{v_0 = v_1, 11, 1}\) is unreachable because in any play of the game the initial valuation is either 00 or 01, and the players never have the opportunity to modify $v_0$.

\begin{lemma}
Let $\phi$ be an IFG-formula, and let $\tuple{\psi, \vec b, \varepsilon}$ be a reachable position of the game \(\game(\A, \phi, V)\). Then 
\begin{enumerate}
	\item \(\psi \in \Sub^+(\phi)\) if and only if \(\varepsilon = 1\).
	\item \(\psi \in \Sub^-(\phi)\) if and only if \(\varepsilon = 0\).
\end{enumerate}
\end{lemma}

\begin{theorem} % conservative extension of first-order logic
\label{conservative extension}
Let $\phi$ be a first-order formula with $N$ variables. We can treat $\phi$ as an IFG$_N$-formula if we interpret $\neg$ as $\hneg$, $\lor$ as $\hor{\emptyset}$, and $\exists v_n$ as $\hexists{v_n}{\emptyset}$. If we do so, then for every suitable structure $\A$ and team \(V \subseteq \^NA\), 
\begin{enumerate}
	\item  $\phi$ is true in $\A$ relative to $V$ if and only if \(\A \models \phi[\vec a]\) for all \(\vec a \in V\),
	\item  $\phi$ is false in $\A$ relative to $V$ if and only if \(\A \not\models \phi[\vec a]\) for all \(\vec a \in V\).
\end{enumerate}
\end{theorem}

\begin{proof}
By two simultaneous inductions on the complexity of $\phi$. A full proof can be found in \cite{Mann:2007}.
\end{proof}

\begin{corollary} % conservative extension on sentences
\label{conservative extension on sentences}
Let $\phi$ be a first-order sentence, and let $\A$ be a suitable structure. Then 
\begin{enumerate}
	\item  $\phi$ is true in $\A$ if and only if \(\A \models \phi\),
	\item  $\phi$ is false in $\A$ if and only if \(\A \not\models \phi\).
\end{enumerate}
\end{corollary}

% Trump Semantics
\subsection{Trump semantics}

The main definition and theorem of this section are due to Hodges\index{Hodges, W.} \cite{Hodges:1997a, Hodges:1997b}. See also related work by Cameron\index{Cameron, P.} and Hodges\index{Hodges, W.} \cite{Cameron:2001}, Caicedo\index{Caicedo, X.} and Krynicki\index{Krynicki, M.} \cite{Caicedo:1999}, V\"a\"an\"anen\index{V\"a\"an\"anen, J.} \cite{Vaananen:2002}, and Dechesne\index{Dechesne, F.} \cite{Dechesne:2005}.

\begin{definition} % cover
Given any set $V$, a \emph{cover}\index{cover|mainidx} of $V$ is a collection of sets $\mathscr U$ such that \(V = \bigcup \mathscr U\). A \emph{disjoint cover} of $V$ is a cover of $V$ whose members are pairwise disjoint. A \emph{partition}\index{partition|mainidx} of $V$ is a disjoint cover of $V$ whose members are all nonempty. Sometimes the members of a partition are called \emph{cells}\index{cell|mainidx}.
\end{definition}

\begin{definition} % J-saturated
Let \(V \subseteq\, \^NA\) be a team, \(J \subseteq N\), and $\mathscr U$ a cover of $V$. Call the cover \emph{$J$-saturated}\index{cover!$J$-saturated|mainidx} if every \(U \in \mathscr U\) is closed under $\approx_J$. That is, for every \(\vec a, \vec b \in V\), if \(\vec a \approx_J \vec b\) and \(\vec a \in U \in \mathscr U\), then \(\vec b \in U\).
\end{definition}

Note that every member of a $J$-saturated cover of $V$ is a union of equivalence classes under $\approx_J$. Also note that every cover of $V$ is $\emptyset$-saturated, and that the only $N$-saturated covers of $V$ are $\set{V}$ and $\set{\emptyset, V}$.

\begin{definition} % bigcup_J
Define a partial operation $\bigcup_J$\index{$\bigcup_J \mathscr U$|mainidx} on collections of sets of valuations by setting \(\bigcup_J \mathscr U = \bigcup \mathscr U\) whenever $\mathscr U$ is a $J$-saturated disjoint cover of $\bigcup \mathscr U$ and letting $\bigcup_J \mathscr U$ be undefined otherwise. Thus the formula \(V = \bigcup_J \mathscr U\) asserts that $\mathscr U$ is a $J$-saturated disjoint cover of $V$. We will use the notation \(V_1 \cup_J V_2\)\index{$V_1 \cup_J V_2$|mainidx} to abbreviate \(\bigcup_J \set{V_1, V_2}\), the notation \(V_1 \cup_J V_2 \cup_J V_3\) to abbreviate \(\bigcup_J \set{V_1, V_2, V_3}\), et cetera.
\end{definition}

\begin{lemma} % disjoint cover
\label{disjoint cover}
Let \(V \subseteq \^NA\) and \(J \subseteq N\). If $\mathscr U$ is a $J$-saturated cover of $V$, then there is a $J$-saturated disjoint cover $\mathscr U'$ of $V$ such that every cell in $\mathscr U'$ is contained in some cell in $\mathscr U$.
\end{lemma}

\begin{proof}
Let \(\seqof{U_\alpha}{\alpha < \mu}\) be an enumeration of $\mathscr U$. Define \(U'_\alpha = U_\alpha \setminus \bigcup_{\beta < \alpha} U_\beta\), and let \(\mathscr U' = \setof{U'_\alpha}{\alpha < \mu}\). By construction, $\mathscr U'$ is a disjoint cover of $V$ such that \(U'_\alpha \subseteq U_\alpha\). To show $\mathscr U'$ is $J$-saturated, suppose \(\vec a, \vec b \in V\). If \(\vec a \approx_J \vec b\) and \(\vec a \in U'_\alpha = U_\alpha \setminus \bigcup_{\beta < \alpha} U_\beta\), then \(\vec b \in U_\alpha \setminus \bigcup_{\beta < \alpha} U_\beta = U'_\alpha\) because $\mathscr U$ is $J$-saturated.
\end{proof}

\begin{lemma} % K saturated J
\label{K saturated J}
Let \(V \subseteq\, \^NA\) and  \(J \subseteq K \subseteq N\). If \(V = \bigcup_K \mathscr U\), then \(V = \bigcup_J \mathscr U\).
\end{lemma}

\begin{proof}
Suppose \(V = \bigcup_K \mathscr U\) and \(\vec a, \vec b \in V\). If \(\vec a \approx_J \vec b\) and \(\vec a \in U \in \mathscr U\), then \(\vec a \approx_K \vec b\), so \(\vec b \in U\). Thus \(V = \bigcup_J \mathscr U\).
\end{proof}

\begin{lemma} % restricted cover preserves saturation
\label{restricted cover preserves saturation}
Let \(V' \subseteq V \subseteq \^NA\). If \(V = \bigcup_J \mathscr U\), then \(V' = \bigcup_J \mathscr U'\), where 
\[
\mathscr U' = \setof{U \cap V'}{U \in \mathscr U}.
\]
\end{lemma}

\begin{proof}
Suppose \(\vec a, \vec b \in V'\). If \(\vec a \approx_J \vec b\), and \(\vec a \in U \cap V' \in \mathscr U'\), then we have \(\vec b \in U \cap V'\) because \(V = \bigcup_J \mathscr U\).
\end{proof}

\begin{lemma} % saturated union preserves saturation
\label{saturated union preserves saturation}
If \(V = V_1 \cup_J V_2\) and \(V_1 = \bigcup_J \mathscr U_1\), \(V_2 = \bigcup_J \mathscr U_2\), then \(V = \bigcup_J (\mathscr U_1 \cup \mathscr U_2)\).
\end{lemma}

\begin{proof}
Suppose \(\vec a, \vec b \in V\), \(\vec a \approx_J \vec b\), and \(\vec a \in U \in \mathscr U_1 \cup \mathscr U_2\). Without loss of generality, we may assume \(U \in \mathscr U_1\), which implies that \(\vec a \in V_1\), in which case \(\vec b \in V_1\) because \(V = V_1 \cup_J V_2\). Hence \(\vec b \in U\) because $V_1 = \bigcup_J \mathscr U_1$.
\end{proof}

\begin{definition} % independent of J
Let \(V \subseteq\, \^NA\) and \(J \subseteq N\). A function \(f\colon V \to A\) is \emph{independent of $J$}, denoted \(f\colon V \toind{J} A\)\index{$f\colon V \toind{J} A$ \quad $f$ is independent of $J$|mainidx},  if \(f(\vec{a}) = f(\vec{b})\) whenever \(\vec{a} \approx_J \vec{b}\).
\end{definition}

Note that any function \(f\colon V \to A\) is independent of $\emptyset$, and that \(f\colon V \to A\) is independent $N$ if and only if $f$ is a constant function.

\begin{lemma} % K independent J
\label{K independent J}
Let \(V \subseteq\, \^NA\) and \(J \subseteq K \subseteq N\). If \(f\colon V \toind{K} A\), then \(f\colon V \toind{J} A\).
\end{lemma}

\begin{proof}
Suppose \(f\colon V \toind{K} A\). If \(\vec a \approx_J \vec b\), then \(\vec a \approx_K \vec b\), so \(f(\vec a) = f(\vec b)\). Thus \(f\colon V \toind{J} A\).
\end{proof}

\begin{lemma} % f cup g independent J
\label{f cup g independent J}
Let \(V \subseteq \^NA\) and \(J \subseteq N\). If \(V = V_1 \cup_J V_2\) and \(f\colon V_1 \toind{J} A\), \(g\colon V_2 \toind{J} A\), then \((f\cup g)\colon V \toind{J} A\).
\end{lemma}
\begin{proof}
Suppose \(\vec a \approx_J \vec b\). Then because \(V = V_1 \cup_J V_2\), either \(\vec a, \vec b \in V_1\) or \(\vec a, \vec b \in V_2\). In the first case, \((f \cup g)(\vec a) = f(\vec a) = f(\vec b) = (f \cup g)(\vec b)\), and in the second case, \((f \cup g)(\vec a) = g(\vec a) = g(\vec b) = (f \cup g)(\vec b)\), because $f$ and $g$ are both independent of $J$. Therefore $(f \cup g)$ is independent of $J$.
\end{proof}

\begin{lemma} % V1(n:f) cup V2(n:f)
\label{V1(n:f) cup V2(n:f)}
If \(f \colon V \toind{J} A\), \(V = V_1 \cup_K V_2\), and \(n \in K\), then \(V(n:f) = V_1(n:f) \cup_K V_2(n:f)\).
\end{lemma}

\begin{proof}
Suppose \(f\colon V \toind{J} A\), \(V = V_1 \cup_K V_2\), and \(n \in K\). Then $V_1(n:f)$ and $V_2(n:f)$ are both subsets of $V(n:f)$, so \(V_1(n:f) \cup V_2(n:f) \subseteq V(n:f)\). Conversely, suppose \(\vec a(n:f(\vec a)) \in \linebreak[1] V(n:f)\), where \(\vec a \in V\). If \(\vec a \in V_1\), then \(\vec a(n: f(\vec a)) \in V_1(n:f)\), and if \(\vec a \in V_2\), then \(\vec a(n: f(\vec a)) \in V_2(n:f)\). Hence \(\vec a(n:f(\vec a)) \in V_1(n:f) \cup V_2(n:f)\). Thus \(V(n:f) \subseteq V_1(n:f) \cup V_2(n:f)\).

To show $V_1(n:f)$ and $V_2(n:f)$ are disjoint, suppose \(\vec a(n:f(\vec a)) \in V_1(n:f)\) and \(\vec b(n:f(\vec b)) \in V_2(n:f)\), where \(\vec a \in V_1\) and \(\vec b \in V_2\). If \(\vec a(n:f(\vec a)) = \vec b(n:f(\vec b))\), then \(\vec a \approx_{\set{n}} \vec b\), hence \(\vec a \approx_K \vec b\) because \(n \in K\), contradicting \(V = V_1 \cup_K V_2\).

Finally suppose \(\vec a(n:f(\vec a)) \approx_K \vec b(n:f(\vec b))\) and \(\vec a(n:f(\vec a)) \in V_i(n:f)\), where \(\vec a \in V_i\). Then \(\vec a \approx_K \vec b\) because \(n \in K\). Hence \(\vec b \in V_i\). Thus \(\vec b(n:f(\vec b)) \in V_i(n:f(\vec b))\).
\end{proof}

\begin{lemma} % composition independent of K
\label{composition independent of K}
Suppose \(f\colon V \toind{J} A\) and \(g\colon V(n:f) \toind{K} A\). Then there is a function \(h\colon V \longtoind{J \cap K} A\) such that \(V(n:f)(n:g) = V(n:h)\). If \(n \in K\), then $h$ is independent of $K$.
\end{lemma}

\begin{proof}
Define \(h\colon V \to A\) by \(h(\vec a) = g(\vec a(n:f(\vec a)))\). To show that $h$ is independent of $J \cap K$, suppose \(\vec a \approx_{J \cap K} \vec b\). Then \(\vec a \approx_{J} \vec b\), so \(f(\vec a) = f(\vec b)\). Hence 
\[
\vec a(n: f(\vec a)) \approx_{J \cap K} \vec b(n:f(\vec b)),
\]
which implies \(\vec a(n: f(\vec a)) \approx_K \vec b(n:f(\vec b))\). Thus 
\[
h(\vec a) = g(\vec a(n: f(\vec a))) = g(\vec b(n:f(\vec b))) = h(\vec b).
\]
Now suppose \(n \in K\) and \(\vec a \approx_{K} \vec b\). Then \(\vec a(n: f(\vec a)) \approx_K \vec b(n:f(\vec b))\), so 
\[
h(\vec a) = g(\vec a(n: f(\vec a))) = g(\vec b(n:f(\vec b))) = h(\vec b).
\]

Let \(f^{(n)}\colon V \to V(n:f)\) be the function that maps $\vec a$ to $\vec a(n:f(\vec a))$. Let $g^{(n)}$ and $h^{(n)}$ be defined similarly. Then \(h^{(n)} = g^{(n)} \circ f^{(n)}\) and 
\[
V(n:h) = h^{(n)}(V) = g^{(n)} \circ f^{(n)}(V) = g^{(n)}(f^{(n)}(V)) = V(n:f)(n:g). \qedhere
\]
\end{proof}

% interpolation independent of K
\begin{lemma}
\label{interpolation independent of K}
Given two functions \(f\colon V \toind{J} A\) and \(h\colon V \toind{K} A\), with \(n \in K\), there is a function \(g\colon V(n:f) \toind{K} A\) such that \(V(n:f)(n:g) = V(n:h)\).
\end{lemma}

\begin{proof}
Define \(g\colon V(n:f) \toind{K} A\) by \(g(\vec a(n:f(\vec a)) = h(\vec a)\). To show that $g$ is well defined and independent of $K$, suppose \(\vec a, \vec b \in V\) and \(\vec a(n:f(\vec a)) \approx_K \vec b(n:f(\vec b))\). Then \(\vec a \approx_K \vec b\) because \(n \in K\). Hence \(g(\vec a(n:f(\vec a)) = h(\vec a) = h(\vec b) = g(\vec b(n:f(\vec b))\). By construction \(h^{(n)} = g^{(n)} \circ f^{(n)}\), so \(V(n:h) = V(n:f)(n:g)\).
\end{proof}

\begin{lemma} % functions commute
\label{functions commute}
Suppose \(f\colon V \toind{J} A\) and \(g\colon V(m:f) \toind{K} A\). If \(m \in K\) and \(n \in J\), where \(m \not= n\), then there exist functions \(G\colon V \toind{K} A\) and \(F\colon V(n:G) \toind{J} A\) such that \(V(m:f)(n:g) = V(n:G)(m:F)\).
\end{lemma}

\begin{proof}
Define \(G\colon V \toind{K} A\) by \(G(\vec a) = g(\vec a(m:f(\vec a)))\) and \(F\colon V(n:G) \toind{J} A\) by \(F(\vec a(n:G(\vec a))) = f(\vec a)\). Observe that $G$ is independent of $K$ because \(m \in K\) and $g$ is independent of $K$. Also observe that $F$ is well defined and independent of $J$ because \(n \in J\) and $f$ is independent of $J$.

To show \(V(m:f)(n:g) = V(n:G)(m:F)\) it suffices to show that \(g^{(n)} \circ f^{(m)} = F^{(m)} \circ G^{(n)}\). So let \(\vec a \in V\). Then 
\begin{align*}
g^{(n)} \circ f^{(m)}(\vec a) &= \vec a(m:f(\vec a))(n:g(\vec a(m:f(\vec a)))) \\
&= \vec a(m: f(\vec a))(n:G(\vec a)) \\
&= \vec a(n:G(\vec a))(m:F(\vec a(n:G(\vec a)))) \\
&= F^{(m)} \circ G^{(n)}(\vec a). \qedhere
\end{align*}
\end{proof}

\begin{lemma} % V(n:A) disjoint cover
\label{V(n:A) disjoint cover}
Suppose \(n \in J\).
\begin{enumerate}
	\item If \(V = V_1 \cup_J V_2\), then \(V(n:A) = V_1(n:A) \cup_J V_2(n:A)\).
	\item If \(V(n:A) = V_1 \cup_J V_2\), then \(V_1 = V_1(n:A)\) and \(V_2 = V_2(n:A)\).
\end{enumerate}
\end{lemma}

\begin{proof}
(a)
Suppose \(V = V_1 \cup_J V_2\). If \(\vec a (n:c) \in V(n:A)\), then either \(\vec a \in V_1\) or \(\vec a \in V_2\). In the first case \(\vec a(n:c) \in V_1(n:A)\), and in the second case \(\vec a(n:c) \in V_2(n:A)\). Hence \(V(n:A) = V_1(n:A) \cup V_2(n:A)\). To show $J$-saturation, suppose \(\vec a, \vec b \in V\), \(\vec a(n:c) \approx_J \vec b(n:d)\), and \(\vec a(n:c) \in V_i(n:A)\). Then \(\vec a \approx_J \vec b\) because \(n \in J\), and \(\vec a \in V_i\), so \(\vec b \in V_i\). Hence \(\vec b(n:d) \in V_i(n:A)\). To show $V_1(n:A)$ and $V_2(n:A)$ are disjoint, suppose \(\vec c \in V_1(n:A) \cap V_2(n:A)\). Then there exist \(\vec a \in V_1\), \(\vec b \in V_2\), and \(d \in A\) such that \(\vec a(n:d) = \vec c = \vec b(n:d)\). But then \(\vec a \approx_J \vec b\), which implies \(\vec a, \vec b \in V_1 \cap V_2\), contradicting the fact that $V_1$ and $V_2$ are disjoint.

(b)
Suppose \(V(n:A) = V_1 \cup_J V_2\). By definition, \(V_i \subseteq V_i(n:A)\). If \(\vec a(n:c) \in V_i(n:A)\), then \(\vec a \approx_J \vec a(n:c)\) and \(\vec a \in V_i\), so \(\vec a(n:c) \in V_i\). Hence \(V_i(n:A) \subseteq V_i\). 
\end{proof}

We now define an Tarski-style satisfaction relation for IFG-formulas.

% (Dechesne Definition 5.3.3)
\begin{definition} % recursive satisfaction for IFG-formulas 
\index{$\A \modeltf \phi[V]$|mainidx}
Let $\phi$ be an IFG$_N$-formula, let $\A$ be a suitable structure, and let \(V \subseteq \^NA\).
\begin{itemize}
	\item If $\phi$ is atomic, then  
	\begin{itemize}
		\item[(+)] \(\A \modelt \phi[V]\) if and only if for every \(\vec a \in V\), \(\A \models \phi[\vec a]\),
		\item[($-$)] \(\A \modelf \phi[V]\) if and only if for every \(\vec a \in V\), \(\A \not\models \phi[\vec a]\).	
	\end{itemize}
	\item If $\phi$ is ${\hneg\psi}$, then 
	\begin{itemize}
		\item[(+)] \(\A \modelt {\hneg\psi[V]}\) if and only if \(\A \modelf \psi[V]\),
		\item[($-$)] \(\A \modelf {\hneg\psi[V]}\) if and only if \(\A \modelt \psi[V]\).
	\end{itemize}
	\item If \(J \subseteq N\) and $\phi$ is $\psi_1 \hor{J} \psi_2$, then 
	\begin{itemize}
		\item[(+)] \(\A \modelt \psi_1 \hor{J} \psi_2[V]\) if and only if there is a $J$-saturated disjoint cover \(V = V_1 \cup_J V_2\) such that \(\A \modelt \psi_1[V_1]\) and \( \A \modelt \psi_2[V_2]\),
		\item[($-$)] \(\A \modelf \psi_1 \hor{J} \psi_2[V]\) if and only if \(\A \modelf \psi_1[V]\) and \( \A \modelf \psi_2[V]\).
	\end{itemize}
	\item If \(J \subseteq N\) and $\phi$ is $\hexists{v_n}{J}\psi$, then 
	\begin{itemize}
		\item[(+)] \(\A \modelt \hexists{v_n}{J}\psi[V]\) if and only if there is a function \(f\colon V \toind{J} A\) independent of $J$ such that \(\A \modelt \psi[V(n:f)]\),
		\item[($-$)] \(\A \modelf \hexists{v_n}{J}\psi[V]\) if and only if \(\A \modelf \psi[V(n:A)]\).
	\end{itemize}
\end{itemize}
\end{definition}

\noindent It is easy to check that the abbreviations $\hand{J}$ and $\hforall{v_n}{J}$ behave as expected.

\begin{itemize}
	\item If \(J \subseteq N\) and $\phi$ is $\psi_1 \hand{J} \psi_2$, then 
	\begin{itemize}
		\item[($+$)] \(\A \modelt \psi_1 \hand{J} \psi_2[V]\) if and only if \(\A \modelt \psi_1[V]\) and \( \A \modelt \psi_2[V]\),
		\item[($-$)] \(\A \modelf \psi_1 \hand{J} \psi_2[V]\) if and only if there is a $J$-saturated disjoint cover \(V = V_1 \cup_J V_2\) such that \(\A \modelf \psi_1[V_1]\) and \( \A \modelf \psi_2[V_2]\).
	\end{itemize}
	\item If \(J \subseteq N\) and $\phi$ is $\hforall{v_n}{J}\psi$, then 
	\begin{itemize}
		\item[($+$)] \(\A \modelt \hforall{v_n}{J}\psi[V]\) if and only if \(\A \modelt \psi[V(n:A)]\),
		\item[($-$)] \(\A \modelf \hforall{v_n}{J}\psi[V]\) if and only if there is a function \(f\colon V \toind{J} A\) independent of $J$ such that \(\A \modelf \psi[V(n:f)]\).
	\end{itemize}
\end{itemize}

\begin{definition} % satisfaction of sentences
Let $\phi$ be an IFG$_N$-formula, and let $\A$ be a suitable structure. Define \(\A \modeltf \phi\)\index{$\A \modeltf \phi$|mainidx} if and only if \(\A \modeltf \phi[\^NA]\).
\end{definition}
 
\begin{lemma} % satisfaction by emptyset
\label{satisfaction by emptyset}
Let $\phi$ be an IFG$_N$-formula, and let $\A$ be a suitable structure. Then \(\A \modeltf \phi[\emptyset]\).
\end{lemma}
\begin{proof}
We proceed by induction on the complexity of $\phi$. If $\phi$ is atomic the lemma holds vacuously. 

Suppose $\phi$ is ${\hneg\psi}$. Then \(\A \modeltf {\hneg\psi[\emptyset]}\) if and only if \(\A \modelft \psi[\emptyset]\), which holds by inductive hypothesis.

Suppose \(J \subseteq N\) and $\phi$ is $\psi_1 \hor{J} \psi_2$. By inductive hypothesis, \(\A \modeltf \psi_1[\emptyset]\) and \(\A \modeltf \psi_2[\emptyset]\). Therefore \(\A \modeltf \psi_1 \hor{J} \psi_2[\emptyset]\).

Suppose \(J \subseteq N\) and $\phi$ is $\hexists{v_n}{J}\psi$. By inductive hypothesis, \(\A \modeltf \psi[\emptyset]\). Let $f$ be the empty function from $\emptyset$ to $A$. Then $f$ is vacuously independent of $J$ and \(\emptyset(n:f) = \emptyset\).  Therefore \(\A \modelt \hexists{v_n}{J} \psi[\emptyset]\). Also \(\emptyset(n:A) = \emptyset\), so \(\A \modelf \hexists{v_n}{J}\psi[\emptyset]\).
\end{proof}

The previous result may seems anomalous, but it is necessary for technical reasons. Specifically, in the positive disjunction clause we wish to allow \(V_1 = V\) and \(V_2 = \emptyset\) or vice versa, which corresponds to the situation in the game $\game(\A, \phi, V)$ when the verifier always wishes to choose the same disjunct. Later we will see that the empty team is the only team that can be winning and losing for the same formula.

The next lemma records the fact that if a player has a winning strategy, given that he or she knows the current valuation belongs to $V$, then that strategy should win given that he or she knows the current valuation belongs to a subset of $V$.
 
% (Dechesne Lemma 6.2.1)
\begin{lemma} % downward closure
\label{downward closure}
Let $\phi$ be an IFG$_N$-formula, let $\A$ be a suitable structure, and let \(V' \subseteq V \subseteq \^NA\). Then \(\A \modeltf \phi[V]\) implies \(\A \modeltf \phi[V']\).
\end{lemma}

\begin{proof}
We proceed by induction on the complexity of $\phi$. The atomic case follows immediately from the definition.

Suppose $\phi$ is ${\hneg\psi}$. Then \(\A \modeltf {\hneg\psi[V]}\) if and only if \(\A \modelft \psi[V]\), which by inductive hypothesis implies \(\A \modelft \psi[V']\), which holds if and only if \(\A \modeltf {\hneg\psi[V']}\).

Suppose \(J \subseteq N\) and $\phi$ is $\psi_1 \hor{J} \psi_2$. If \(\A \modelt \psi_1 \hor{J} \psi_2[V]\), then there is a $J$-saturated disjoint cover \(V = V_1 \cup_J V_2\) such that \(\A \modelt \psi_1[V_1]\) and \(\A \modelt \psi_2[V_2]\). Let \(V'_1 = V_1 \cap V'\) and \(V'_2 = V_2 \cap V'\). Then by \lemref{restricted cover preserves saturation}, \(V' = V'_1 \cup_J V'_2\) is a $J$-saturated disjoint cover, and by inductive hypothesis \(\A \modelt \psi_1[V'_1]\) and \(\A \modelt \psi_2[V'_2]\). Thus \(\A \modelt \psi_1 \hor{J} \psi_2[V']\). 

If \(\A \modelf \psi_1 \hor{J} \psi_2[V]\), then \(\A \modelf \psi_1[V]\) and \(\A \modelf \psi_2[V]\). By inductive hypothesis, \(\A \modelf \psi_1[V']\) and \(\A \modelf \psi_2[V']\). Thus \(\A \modelf \psi_1 \hor{J} \psi_2[V']\).

Suppose \(J \subseteq N\) and $\phi$ is $\hexists{v_n}{J}\psi$. If \(\A \modelt \hexists{v_n}{J}\psi[V]\), then there is a function \(f\colon V \toind{J} A\) independent of $J$ such that \(\A \modelt \psi[V(n:f)]\). The restriction of $f$ to $V'$ is independent of $J$, and by inductive hypothesis \(\A \modelt \psi[V'(n:f)]\). Thus \(\A \modelt \hexists{v_n}{J}\psi[V']\).

If \(\A \modelf \hexists{v_n}{J}\psi[V]\), then \(\A \modelf \psi[V(n:A)]\). By inductive hypothesis we have \(\A \modelf \psi[V'(n:A)]\). Thus \(\A \modelf \hexists{v_n}{J}\psi[V']\).
\end{proof}

Next we show that a given team cannot be winning and losing for the same formula, which reflects the fact that \abelard\ and \eloise\ cannot both have winning strategies for the same semantic game. 

\begin{lemma} % noncontradiction
\label{noncontradiction}
Let $\phi$ be an IFG$_N$-formula, and let $\A$ be a suitable structure. If \(V \subseteq \^NA\) is a nonempty team, then we cannot have
\[
\A \modelt \phi[V] \quad \text{and} \quad \A \modelf \phi[V].
\]
\end{lemma}

\begin{proof}
We proceed by induction on the complexity of $\phi$. Suppose \(\A \modelt \phi[V]\). If $\phi$ is atomic, then because $V$ is nonempty there is an \(\vec a \in V\) such that \(\A \models \phi[\vec a]\). Hence \(\A \not\modelf \phi[V]\).

Suppose $\phi$ is ${\hneg\psi}$. Then \(\A \modelf \psi[V]\), so by inductive hypothesis \(\A \not\modelt \psi[V]\). Hence \(\A \not\modelf {\hneg\psi[V]}\).

Suppose \(J \subseteq N\) and $\phi$ is $\psi_1 \hor{J} \psi_2$. Then there is a $J$-saturated disjoint cover \(V = V_1 \cup_J V_2\) such that \(\A \modelt \psi_1[V_1]\) and \(\A \modelt \psi_2[V_2]\). By inductive hypothesis, \(\A \not\modelf \psi_1[V_1]\) and \(\A \not\modelf \psi_2[V_2]\). It follows from \lemref{downward closure} that \(\A \not\modelf \psi_1[V]\) and \(\A \not\modelf \psi_2[V]\). Hence \(\A \not\modelf \psi_1 \hor{J} \psi_2[V]\).

Suppose \(n < N\), \(J \subseteq N\), and $\phi$ is \(\hexists{v_n}{J}\psi\). Then there is an \(f\colon V \toind{J} A\) such that \(\A \modelt \psi[V(n:f)]\). By inductive hypothesis, \(\A \not\modelf \psi[V(n:f)]\). It follows from \lemref{downward closure} that \(\A \not\modelf \psi[V(n:A)]\). Hence \(\A \not\modelf \hexists{v_n}{J}\psi[V]\).
\end{proof}

The main result of this section is that trump semantics and game semantics are equivalent. This is significant because of Hintikka's\index{Hintikka, J.} claim in \book{The Principles of Mathematics Revisited} that independence-friendly logic does not have a compositional semantics.

% (Dechesne Theorem 5.3.5)
\begin{theorem} % game semantics = compositional semantics
\label{game semantics = compositional semantics}
Let $\phi$ be an IFG$_N$-formula, let $\A$ be a suitable structure, and let \(V \subseteq \^NA\). Then 
\begin{enumerate}
	\item \(\A \modelt \phi[V]\) if and only if $\phi$ is true in $\A$ relative to $V$;
	\item \(\A \modelf \phi[V]\) if and only if $\phi$ is false in $\A$ relative to $V$.
\end{enumerate}
\end{theorem}
\begin{proof}
By two simultaneous inductions on the complexity of $\phi$. A full proof can be found in \cite{Mann:2007}.
\end{proof}

\begin{corollary}
\begin{enumerate}
	\item \(\A \modelt \phi\) if and only if $\phi$ is true in $\A$.
	\item \(\A \modelf \phi\) if and only if $\phi$ is false in $\A$.
\end{enumerate}
\end{corollary}

The next two results show that IFG-formulas that differ only in their number of variables have essentially the same meaning, as long as we do not encode any information in the extra variables. The phenomenon of one player using extra variables to send (otherwise forbidden) information to him or herself is called \emph{signaling}\index{signaling}. For example, let $\phi$ be the formula 
\[
v_0 = v_1 \hor{\set{0,1}} v_0 \not= v_1,
\]
and let $\A$ be the equality structure over $\set{0,1}$. If \(V = \set{00, 01, 10, 11}\), then \eloise\ does not have a winning strategy for \(\game(\A, \phi, V)\). However, if $\phi'$ is the 3-variable version of $\phi$, and \(V' = \set{001, 010, 100, 111}\), then \eloise\ does have a winning strategy for \(\game(\A, \phi', V')\) because she can use the value of $v_2$ to signal whether the values of $v_0$ and $v_1$ are equal.

% Dechesne Lemma 6.2.2
\begin{lemma} % expansion by Cartesian products
\label{expansion by Cartesian products}
Let $\phi$ be an IFG-formula with $M$ variables, and let $\phi'$ be the same formula but with $N$ variables,  where \(M < N\). Let $\A$ be a suitable structure, \(V \subseteq \^NA\), \(\emptyset \not= W \subseteq \^{N \setminus M}A\), and \(V' = \setof{\vec a \cup \vec b}{\vec a \in V,\ \vec b \in W}\). Then \(\A \modeltf \phi[V]\) if and only if \(\A \modeltf \phi'[V']\).
\end{lemma}
\begin{proof}
By induction on the complexity of $\phi$. A full proof can be found in \cite{Mann:2007}.
\end{proof}
\begin{corollary} % evaluation of sentences
\label{evaluation of sentences}
Let $\phi$ be an IFG-formula with $M$ variables, let $\phi'$ be the same formula but with $N$ variables,  where \(M < N\), and let $\A$ be a suitable structure. Then \(\A \modeltf \phi\) if and only if  \(\A \modeltf \phi'\).
\end{corollary}
\begin{proof}
By \lemref{expansion by Cartesian products}, \(\A \modeltf \phi\) if and only if \(\A \modeltf \phi[\^MA]\) if and only if \(\A \modeltf \phi'[\^NA]\) if and only if \(\A \modeltf \phi'\).
\end{proof}

Henceforth, we will use the same symbol for formulas that differ only in the number of variables they have. For example, we will write $\phi[V']$ instead of $\phi'[V']$.

%\input{2-IFG_Algebras}
%LaTeX by Allen L. Mann

% Independence-Friendly Cylindric Set Algebras
\section{Independence-Friendly Cylindric Set Algebras}

\begin{definition} % independence-friendly cylindric power set algebra
An \emph{independence-friendly cylindric power set algebra}\index{independence-friendly cylindric power set algebra|mainidx} is an algebra whose universe is \(\powerset(\powerset(\^NA)) \times \powerset(\powerset(\^NA))\), where $A$ is a set and $N$ is a natural number. The set $A$ is called the \emph{base set}\index{base set|mainidx}, and the number $N$ is called the \emph{dimension}\index{dimension|mainidx} of the algebra. Since each  element $X$ is an ordered pair, we will use the notation $X^+$\index{$X^+$ \quad truth coordinate of $X$|mainidx} to refer to the first coordinate of the pair, and $X^-$\index{$X^-$ \quad falsity coordinate of $X$|mainidx} to refer to the second coordinate. There are a finite number of operations: 

\begin{itemize}
	\item the constant \(0 = \tuple{\set{\emptyset}, \powerset(\^NA)}\);
	\item the constant \(1 = \tuple{\powerset(\^NA), \set{\emptyset}}\);
	\item for all \(i,j < N\), the constant $D_{ij}$\index{$D_{ij}$ \quad diagonal element|mainidx} is defined by 
	\begin{itemize}
		\item[(+)] \(D_{ij}^+ = \powerset(\setof{\vec a \in\, \^NA}{a_i = a_j})\),
		\item[($-$)] \(D_{ij}^- = \powerset(\setof{\vec a \in\, \^NA}{a_i \not= a_j})\);
	\end{itemize}
	\item if \(X = \tuple{X^+,\, X^-}\), then \(\n{X}\index{$\n{X}$ \quad negation of $X$|mainidx} = \tuple{X^-, X^+}\);
	\item for every \(J \subseteq N\), the binary operation $+_J$\index{$X +_J Y$|mainidx} is defined by
	\begin{itemize}
		\item[(+)] \(V \in (X +_J Y)^+\) if and only if \(V = V_1 \cup_J V_2\) for some
				 \(V_1 \in X^+\) and \(V_2 \in Y^+\),
		\item[($-$)] \((X +_J Y)^- = X^- \cap Y^-\);
	\end{itemize}
	\item for every \(J \subseteq N\), the binary operation $\cdot_J$\index{$X \cdot_J Y$|mainidx} is defined by
	\begin{itemize}
		\item[(+)] \((X \cdot_J Y)^+ = X^+ \cap Y^+\),
		\item[($-$)] \(W \in (X \cdot_J Y)^-\) if and only if \(W = W_1 \cup_J W_2\) for some
				 \(W_1 \in X^-\) and \(W_2 \in Y^-\);
	\end{itemize}
	\item for every \(n < N\) and \(J \subseteq N\), the unary operation $C_{n,J}$\index{$C_{n,J}(X)$ \quad cylindrification of $X$|mainidx} is defined by
		\begin{itemize}
			\item[(+)] \(V \in C_{n,J}(X)^+\) if and only if \(V(n:f) \in X^+\) for some 
					\(f\colon V \toind{J} A\),
			\item[($-$)] \(W \in C_{n,J}(X)^- \) if and only if \(W(n:A) \in X^-\).
		\end{itemize}
\end{itemize}
\end{definition}

\begin{definition} % independence-friendly cylindric set algebra
An \emph{independence-friendly cylindric set algebra}\index{independence-friendly cylindric set algebra|mainidx} (or \emph{IFG-algebra}, for short) is any subalgebra of an independence-friendly cylindric power set algebra. An \emph{IFG$_N$-cylindric set algebra}\index{IFG$_N$-cylindric set algebra|mainidx} (or \emph{IFG$_N$-algebra}) is an independence-friendly cylindric set algebra of dimension $N$.
\end{definition}

We will use the notation \(X =^+ Y\) to abbreviate \(X^+ = Y^+\) and \(X =^- Y\) to abbreviate \(X^- = Y^-\). Furthermore, \( X \leq^+ Y\) abbreviates \(X^+ \subseteq Y^+\), and \(X \leq^- Y\) abbreviates \(X^- \subseteq Y^-\). Define \(X \leq Y\) if and only if \(X \leq^+ Y\) and \(Y \leq^- X\). It follows immediately from the definitions that \(X \leq^\pm Y\) if and only if \(\n{X} \leq^\mp \n{Y}\). Hence, \(X \leq Y\) if and only if \(\n{Y} \leq \n{X}\).

% Duality
\subsection{Duality}

We include the operations $\cdot_J$ in the signature of independence-friendly cylindric set algebras for the sake of compatibility with De Morgan algebra\index{De Morgan algebra}. We could have omitted $\cdot_J$ from the definition and instead defined it in terms of $+_J$ and $\n{ }$.

\begin{lemma} % + and times
\label{+ and times}
\(X \cdot_J Y = \n{(\n{X} +_J \n{Y})}\).
\end{lemma}

\begin{proof} It follows immediately from the definition that \(\n{(\n X)} = X\) and \((\n X)^\pm = X^\mp\). Hence
\begin{align*}
(X \cdot_J Y)^+ &= X^+ \cap Y^+ \\
&= (\n{X})^- \cap (\n{Y})^- \\
&= (\n{X} +_J \n{Y})^- \\
&= (\n{(\n{X} +_J \n{Y})})^+.
\end{align*}

Also \(W \in (X \cdot_J Y)^-\) if and only if \(W = W_1 \cup_J W_2\) for some \(W_1 \in X^- = (\n X)^+\) and \(W_2 \in Y^- = (\n Y)^+\) if and only if \(W \in (\n X +_J \n Y)^+ = (\n{(\n X +_J \n Y)})^-\). 
\end{proof}

The cylindrifications $C_{n,J}$ also have their corresponding dual operations, defined by \[\Cd{n,J}(X) = \n{C_{n,J}(\n{X})}.\] We do not include the dual cylindrifications in the signature of IFG$_N$-cylindric set algebras for the sake of compatibility with cylindric algebra.

Every IFG$_N$-cylindric set algebra 
\[
\C = \tuple{U;\, 0,\, 1,\, D_{ij},\, \n{\,},\, +_J,\, \cdot_J,\, C_{n,J}}
\]
has a dual algebra
\[
\C ^\partial\index{$\C ^\partial$ \quad dual algebra} = \tuple{U;\, 1,\, 0,\, \n{D_{ij}},\, \n{\,},\, \cdot_J,\, +_J,\, \Cd{n,J}}.
\]
Furthermore, the operation \!$\n{\,}$ is an isomorphism from $\C$ to $\C^\partial$. Therefore any algebraic statement about independence-friendly cylindric set algebras can be dualized by systematically exchanging the symbols 0 with 1, $D_{ij}$ with $\n{D_{ij}}$, $+_J$ with $\cdot_J$, and $C_{n,J}$ with $\Cd{n,J}$. Statements involving the superscripts $^+$ and $^-$ can be dualized by their exchange, even though the symbols $^+$ and $^-$ do not belong to the signature. In particular, the dual of $\leq$ is $\geq$.
Henceforth we will often state a theorem and its dual together, but we will not prove the dual.

% Fixed points
\subsection{Fixed points}

An element $X$ of an IFG-algebra is a \emph{fixed point}\index{fixed point} if \(\n{X} = X\). Clearly $X$ is a fixed point if and only if \(X^+ = X^-\). Two fixed points in particular deserve special mention. In a given IFG$_N$-algebra, let \(\Omega\index{$\Omega$} = \tuple{\set{\emptyset}, \set{\emptyset}}\) and \(\mho\index{$\mho$} = \tuple{\powerset(\^NA), \powerset(\^NA)}\). Neither the symbol $\Omega$ nor the symbol $\mho$ belong to the signature of IFG-algebras. However, we will see that $\Omega$ is definable in most IFG-algebras, and it is often present even when it is not definable.

\begin{lemma} % Omega1
\label{Omega1}
In any IFG$_N$-algebra $\C$, if there exists an element \(X \in \C\) such that \(X^+ \cap X^- = \set{\emptyset}\) and \(\^NA \notin X^+ \cup X^-\), then \(\Omega \in \C\).
\end{lemma}

\begin{proof}
Suppose $X$ is such an element, and consider \(C_{0,\emptyset} \ldots C_{N-1,\emptyset}(X \cdot_N \n{X})\). Observe that \(V \in C_{0,\emptyset} \ldots C_{N-1,\emptyset}(X \cdot_N \n{X})^+\) if and only if the exist functions \(f_0, \ldots, f_{N-1}\) such that \(V(0:f_0)\ldots(N-1:f_{N-1}) \in (X \cdot_N \n{X})^+ = X^+ \cap X^- = \set{\emptyset}\) if and only if \(V = \emptyset\). Also, \(W \in C_{0,\emptyset} \ldots C_{N-1,\emptyset}(X \cdot_N \n{X})^-\) if and only if \(W(0:A) \ldots (N-1:A) \linebreak[1] \in (X \cdot_N \n{X})^-\). Note that if \(W = \emptyset\) then \(W(0:A) \ldots (N-1:A) = \emptyset\), and if \(W \not= \emptyset\) then \(W(0:A) \ldots (N-1:A) = \^NA\). Also note that \(\emptyset \in (X \cdot_N \n{X})^-\) because \(\emptyset \in X^+ \cap X^-\), while \(\^NA \notin (X \cdot_N \n{X})^-\) because if it did we would have \(\^NA \in X^+\) or \(\^NA \in X^-\). Thus \(W \in C_{N-1,\emptyset} \ldots C_{0,\emptyset}(X \cdot_N \n{X})^-\) if and only if \(W = \emptyset\).
\end{proof}

In any IFG-algebra of dimension at least two and whose base set  has at least two elements we have \(D_{01}^+ \cap D_{01}^- = \emptyset\) and \(\^NA \notin D_{01}^+ \cup D_{01}^-\). Hence 
\[
C_{0,\emptyset} \ldots C_{N-1,\emptyset}(D_{01} \cdot_N \n{D_{01}}) = \Omega.
\]
In the next section we will show that if \(\abs{A} = 1\) or \(N = 1\) then there are examples of IFG-algebras that do not include $\Omega$.

% $\Cyls_{IFG_N}(\A)$
\subsection{$\Cyls_{IFG_N}(\A)$}

The motivating examples of IFG-algebras are those whose elements are the meanings of IFG-formulas.

\begin{definition} % \norm{\phi}
Given a structure $\A$ and an IFG$_N$-formula $\phi$, define 
\begin{align*}
\trump{\phi}_\A &= \setof{V \subseteq\, \^NA}{\A \modelt \phi[V]}, \\
\cotrump{\phi}_\A &= \setof{W \subseteq\, \^NA}{\A \modelf \phi[W]}, \\
\norm{\phi}_\A &= \tuple{\trump{\phi}_\A, \cotrump{\phi}_\A}.
\end{align*}
\index{$\trump{\phi}_\A$ \quad winning teams for $\phi$ |mainidx}
\index{$\cotrump{\phi}_\A$ \quad losing teams for $\phi$ |mainidx}
\index{$\norm{\phi}_\A$ \quad meaning of $\phi$ |mainidx}
\index{IFG-formula!meaning of|mainidx}
By \thmref{game semantics = compositional semantics}, $\trump{\phi}_\A$ is the set of all winning teams for $\phi$ in $\A$, and $\cotrump{\phi}_\A$ is the set of all losing teams. We will call $\norm{\phi}_\A$ the \emph{meaning} of $\phi$ in $\A$. If the structure is clear from context we will omit the subscript and simply write $\trump{\phi}$, $\cotrump{\phi}$, and $\norm{\phi}$.
\end{definition}

\begin{definition} % IFG_N-cylindric set algebra over A
Given an $\sigma$-structure $\A$, the \emph{IFG$_N$-algebra over $\A$}\index{IFG$_N$-cylindric set algebra over $\A$|mainidx}, denoted $\Cyls_{IFG_N}(\A)$\index{$\Cyls_{IFG_N}(\A)$ \quad IFG$_N$-cylindric set algebra over $\A$|mainidx}, is the IFG$_N$-algebra whose universe is the set \(\setof{\norm{\phi}_\A}{\phi \in \mathscr L^\sigma_{{IFG}_N}}\). Observe that if \(X = \norm{\phi}\) and \(Y = \norm{\psi}\), then 
\begin{align*}
0 &= \norm{v_0 \not= v_0}, \\
1 &= \norm{v_0 = v_0}, \\
D_{ij} &= \norm{v_i = v_j}, \\
\n{X} &= \norm{\hneg\phi}, \\
X +_J Y &= \norm{\phi \hor{J} \psi}, \\
X \cdot_J Y &= \norm{\phi \hand{J} \psi}, \\
C_{n,J}(X) &= \norm{\hexists{v_n}{J} \phi}.
\end{align*}
\end{definition}

\begin{lemma} % emptyset
\label{emptyset}
For any formula $\phi$ and suitable structure $\A$ we have 
\[
\trump{\phi}_\A \cap \cotrump{\phi}_\A = \set{\emptyset}.
\]
\end{lemma}

\begin{proof}
By \lemref{satisfaction by emptyset} and \lemref{noncontradiction}.
\end{proof}

\begin{lemma} % algebraic downward closure
\label{algebraic downward closure}
For any formula $\phi$ and suitable structure $\A$, if \(V' \subseteq V \in \norm{\phi}^\pm_\A\), then \(V' \in \norm{\phi}^\pm_\A\).
\end{lemma}

\begin{proof}
By \lemref{downward closure}.
\end{proof}

In ordinary first-order logic a sentence can have one of only two possible truth-values: \emph{true} or \emph{false}. Thus if $\phi$ is a sentence either \(\phi^\A = 0\) or \(\phi^\A = 1\). In IFG logic a sentence $\phi$ may be neither true nor false, so $\norm{\phi}_\A$ can take values other than  0 or 1. Interestingly, the only other possible value is $\Omega$.

\begin{proposition} % three truth-values
\label{three truth-values}
If $\phi$ is an IFG$_N$-sentence and $\A$ is a suitable structure, then \(\norm{\phi}_\A \in \set{0, \Omega, 1}\).
\end{proposition}

\begin{proof}
If \(\norm{\phi}_\A \not= \Omega\), then there is a nonempty \(V \subseteq \^NA\) such that \(V \in \trump{\phi}_\A\) or \(V \in \cotrump{\phi}_\A\). By \propref{true for any valuation -> true for all}, in the first case \(\norm{\phi}_\A = 1\), and in the second case \(\norm{\phi}_\A = 0\).
\end{proof}

\begin{lemma} % A = 1
\label{A = 1}
If \(\abs{\A} = 1\), then for every IFG$_N$-formula $\phi$, either \(\A \modelt \phi\) or \(\A \modelf \phi\).
\end{lemma}

\begin{proof}
If \(\abs{\A} = 1\), then \(\abs{\^NA} = 1\). Let \(\^NA = \set{\vec a}\). 

Suppose $\phi$ is atomic. Then either \(\A \models \phi[\vec a]\) or \(\A \not\models \phi[\vec a]\). In the first case \(\A \modelt \phi[\set{\vec a}]\), and in the second case \(\A \modelf \phi[\set{\vec a}]\).

Suppose $\phi$ is $\hneg\psi$. Then by inductive hypothesis \(\A \modelt \psi[\set{\vec a}]\) or \(\A \modelf \psi[\set{\vec a}]\), in which case \(\A \modelf \hneg\psi[\set{\vec a}]\) or \(\A \modelt \hneg\psi[\set{\vec a}]\).

Suppose $\phi$ is \(\psi_1 \hor{J} \psi_2\). If \(\A \modelt \psi_1[\set{\vec a}]\), then \(\A \modelt \psi_1 \hor{J} \psi_2[\set{\vec a}]\) because \(\set{\vec a} = \set{\vec a} \cup_J \emptyset\) and it is always the case that \(\A \modelt \psi_2[\emptyset]\). Similarly, if \(\A \modelt \psi_2[\set{\vec a}]\), then \(\A \modelt \psi_1 \hor{J} \psi_2[\set{\vec a}]\). If \(\A \modelf \psi_1[\set{\vec a}]\) and \(\A \modelf \psi_2[\set{\vec a}]\), then \(\A \modelf \psi_1 \hor{J} \psi_2[\set{\vec a}]\).

Suppose $\phi$ is \(\hexists{v_n}{J}\psi\). If \(\A \modelt \psi[\set{\vec a}]\), let \(f\colon \set{\vec a} \toind{J} A\) be the function that sends $\vec a$ to $a_0$. Then \(\set{\vec a}(n:f) = \set{\vec a}\), so \(\A \modelt \psi[\set{\vec a}(n:f)]\). Hence \(\A \modelt \hexists{v_n}{J}\psi[\set{\vec a}]\). If \(\modelf \psi[\set{\vec a}]\), then \(\set{\vec a}(n:A) = \set{\vec a}\), so \(\A \modelf \psi[\set{\vec a}(n:A)]\). Hence \(\A \modelf \hexists{v_n}{J}\psi[\set{\vec a}]\).
\end{proof}

Thus if \(\abs{\A} = 1\), \(\Cyls_{\mathrm{IFG}_{N}}(\A)\) is essentially the same structure as the two-element Boolean algebra.

\begin{lemma} % Omega in Cs(A) 1
\label{Omega in Cs(A) 1}
Let $\A$ be a structure. The element \(\Omega \in \Cyls_{\mathrm{IFG}_{N}}(\A)\) if and only if there is an atomic IFG$_N$-formula $\phi$ such that \(\A \not\modelt \phi\) and \(\A \not\modelf \phi\).
\end{lemma}

\begin{proof}
Suppose \(\A \not\modelt \phi[\^NA]\) and \(\A \not\modelf \phi[\^NA]\). Then \(\trump{\phi}_\A \cap \cotrump{\phi}_\A = \set{\emptyset}\) by \lemref{emptyset}, and \(\^NA \notin \trump{\phi}_\A \cup \cotrump{\phi}_\A\) by hypothesis. Hence \(\Omega \in \Cyls_{\mathrm{IFG}_{N}}(\A)\) by \lemref{Omega1}. Conversely, suppose that for every atomic IFG$_N$-formula, either \(\A \modelt \phi[\^NA]\) or \(\A \modelf \phi[\^NA]\). We wish to show that for every IFG$_N$-formula, either \(\A \modelt \phi[\^NA]\) or \(\A \modelf \phi[\^NA]\).

Suppose $\phi$ is $\hneg\psi$. By inductive hypothesis, either \(\A \modelt \psi[\^NA]\) or \(\A \modelf \psi[\^NA]\). In the first case \(\A \modelf \hneg\psi[\^NA]\), and in the second case \(\A \modelt \hneg\psi[\^NA]\).

Suppose $\phi$ is \(\psi_1 \hor{J} \psi_2\). Then by inductive hypothesis \(\A \modeltf \psi_1[\^NA]\) and \(\A \modeltf \psi_2[\^NA]\). If \(\A \modelt \psi_1[\^NA]\) or \(\A \modelt \psi_2[\^NA]\),  then \(\A \modelt \psi_1 \hor{J} \psi_2[\^NA]\). If \(\A \modelf \psi_1[\^NA]\) and \(\A \modelf \psi_2[\^NA]\), then \(\A \modelf \psi_1 \hor{J} \psi_2[\^NA]\).

Suppose $\phi$ is \(\hexists{v_n}{J}\psi\). Then by inductive hypothesis \(\A \modeltf \psi[\^NA]\). If \(\A \modelt \psi[\^NA]\), let \(f\colon \^NA \toind{J} A\) be any function independent of $J$ (e.g., a constant function). Then \(\^NA(n:f) \subseteq \^NA\), so \(\A \modelt \psi[\^NA(n:f)]\). Hence \(\A \modelt \hexists{v_n}{J}\psi[\^NA]\). If \(\A \modelf \psi[\^NA]\), then \(\A \modelf \psi[\^NA(n:A)]\) because \(\^NA(n:A) = \^NA\). Hence \(\A \modelf \hexists{v_n}{J}\psi[\^NA]\).
\end{proof}

\begin{proposition} % Omega in Cs(A) 2
\label{Omega in Cs(A) 2}
Let $\A$ be a nonempty structure. Then \(\Omega \in \Cyls_{\mathrm{IFG}_{N}}(\A)\) if and only if at least one of the following conditions is satisfied:
\begin{enumerate}
	\item \(\abs{\A} \geq 2\) and \(N \geq 2\);
	\item \(\abs{\A} \geq 2\), \(N = 1\), and there exist terms $s$ and $t$ and elements \(x,y \in \A\) such that \(s^\A(x) = t^\A(x)\) and \(s^\A(y) \not= t^\A(y)\);
	\item \(\abs{\A} \geq 2\), \(N = 1\), and there exist an $m$-ary relation symbol $R$, terms \(\seqof{s_i}{i < m}\), and elements \(x,y \in \A\) such that 
	\[
	\tuple{s_0^\A(x), \ldots, s_{m-1}^\A(x)} \in R^\A \quad \text{and} \quad \tuple{s_0^\A(y), \ldots, s_{m-1}^\A(y)} \notin R^\A.
	\]
\end{enumerate}
\end{proposition}

\begin{proof}
If (a) is true, then \(\A \not\modelt v_0 = v_1[\^NA]\) and \(\A \not\modelf v_0 = v_1[\^NA]\). If (b) is true, then \(\A \not\modelt s = t[A]\) and \(\A \not\modelf s = t[A]\). If (c) is true, then \(\A \not\modelt Rs_0\ldots s_{m-1}[A]\) and \(\A \not\modelf Rs_0\ldots s_{m-1}[A]\).

Conversely, suppose (a), (b), and (c) are all false. If \(\abs{\A} = 1\) then \(\Omega \notin \Cyls_{\mathrm{IFG}_{N}}(\A)\) by \lemref{A = 1}. Suppose \(\abs{\A} \geq 2\) and \(N = 0\). If $\A$ is a relational structure then \(\Cyls_{\mathrm{IFG}_{0}}(\A) = \emptyset\) because there are no IFG$_0$-formulas in a purely relational language. Otherwise every atomic IFG$_0$-formula has the form \(s = t\) or \(Rs_0\ldots s_{m-1}\), where all of the terms involved are closed. In the first case either \(\A \modelt s = t[\set{\emptyset}]\) or \(\A \modelf s = t[\set{\emptyset}]\), and in the second case either \(\A \modelt Rs_0\ldots s_{m-1}[\set{\emptyset}]\) or \(\A \modelf Rs_0\ldots s_{m-1}[\set{\emptyset}]\). 

Suppose \(\abs{\A} \geq 2\), \(N = 1\), every pair of terms $s$ and $t$ has the property that either for every \(x \in \A\), \(s^\A(x) = t^\A(x)\), or for every \(x \in \A\), \(s^\A(x) \not= t^\A(x)\), and every relation symbol $R$ and sequence of terms \(\seqof{s_i}{i < m}\) has the property that either for every \(x \in \A\), \(\tuple{s_0^\A(x), \ldots, s_{m-1}^\A(x)} \in R^\A\) or for every \(x \in \A\), \(\tuple{s_0^\A(x), \ldots, s_{m-1}^\A(x)} \notin R^\A\). Every atomic IFG$_1$-formula has the form \(s = t\) or \(Rs_0\ldots s_{m-1}\), where all of the terms involve at most one variable. In the first case, either \(\A \modelt s = t[A]\) or \(\A \modelf s = t[A]\). In the second case, \(\A \modelt Rs_0\ldots s_{m-1}[A]\) or \(\A \modelf Rs_0\ldots s_{m-1}[A]\).
\end{proof}
 
% Rooted, suited, and double-suited IFG-algebras
\subsection{Rooted, suited, and double-suited IFG-algebras}

\lemref{emptyset} and \lemref{algebraic downward closure} inspire the following definitions.

\begin{definition} % rooted
A set of teams $X^*$ is \emph{rooted}\index{rooted element|mainidx} if \(\emptyset \in X^*\). A pair \(\tuple{X^+, X^-}\) of sets of teams is \emph{rooted} if both of its coordinates are rooted. An IFG-algebra is \emph{rooted}\index{IFG-cylindric set algebra!rooted|mainidx} if all of its elements are rooted.
\end{definition}

\begin{proposition} % subalgebra generated rooted
\label{subalgebra generated rooted}
The subalgebra of an IFG-algebra generated by a set of rooted elements is a rooted IFG-algebra.
\end{proposition}

\begin{proof}
The constant elements $0$, $1$, and $D_{ij}$ are all rooted, and if $X$ is rooted, then $\n{X}$ is.

Suppose $X$ and $Y$ are rooted. Then \(\emptyset \in (X +_J Y)^+\) because \(\emptyset = \emptyset \cup_J \emptyset\) and \(\emptyset \in X^+\), \(\emptyset \in Y^+\). Also, \(\emptyset \in X^- \cap Y^- = (X +_J Y)^-\).

Suppose $X$ is rooted. Then \(\emptyset \in C_{n,J}(X)^+\) because the empty function $f$ from  $\emptyset$ to $A$ is vacuously independent of $J$ and \(\emptyset(n:f) = \emptyset \in X^+\). Also, \(\emptyset \in C_{n,J}(X)^-\) because \(\emptyset(n:A) = \emptyset \in X^-\).
\end{proof}

Given a set $A$, let $\Root_N(A)$\index{$\Root_N(A)$ \quad rooted IFG$_N$-cylindric set algebra over $A$} denote the IFG$_N$-cylindric set algebra whose universe is the set of all rooted elements in \(\powerset(\powerset(\^NA)) \times \powerset(\powerset(\^NA))\).

\begin{definition} % suit 							%\cite{Cameron:2001}
A nonempty subset of $\powerset(\^NA)$ is called a \emph{suit}\index{suit|mainidx} if it is closed under subsets. That is, a suit $X^*$ is a nonempty collection of subsets of $\^NA$ such that \(V' \subseteq V \in X^*\) implies \(V' \in X^*\). A \emph{double suit}\index{double suit|mainidx} is a pair $\tuple{X^+, X^-}$ of suits such that \(X^+ \cap X^- = \set{\emptyset}\).
\end{definition}

\begin{definition} % suited
An independence-friendly cylindric set algebra is \emph{suited}\index{IFG-cylindric set algebra!suited|mainidx} if all of its elements are pairs of suits. It is \emph{double-suited}\index{IFG-cylindric set algebra!double-suited|mainidx} if all of its elements are double suits.
\end{definition}

Note that suits, double suits, and suited independence-friendly cylindric set algebras are all rooted.

\begin{proposition} % subalgebra generated suited
\label{subalgebra generated suited}
The subalgebra of an IFG$_N$-algebra generated by a set of pairs of suits is a suited IFG$_N$-algebra.
\end{proposition}

\begin{proof}
The constant elements $0$, $1$, and $D_{ij}$ are all pairs of suits, and if $X$ is a pair of suits, then so is $\n{X}$.

Suppose $X$ and $Y$ are pairs of suits. If \(V' \subseteq V \in (X +_J Y)^+\), then \(V = V_1 \cup_J V_2\) for some \(V_1 \in X^+\), \(V_2 \in Y^+\). Let \(V_1' = V_1 \cap V'\) and \(V_2' = V_2 \cap V'\). By \lemref{restricted cover preserves saturation}, \(V' = V_1' \cup_J V_2'\), and \(V_1' \in X^+\), \(V_2' \in Y^+\). Hence \(V' \in (X +_J Y)^+\). Thus $(X +_J Y)^+$ is a suit. If \(W' \subseteq W \in (X +_J Y)^- = X^- \cap Y^-\), then \(W' \subseteq W \in X^-\) and \(W' \subseteq W \in Y^-\). Hence \(W' \in X^- \cap Y^- = (X +_J Y)^-\). Thus $(X +_J Y)^-$ is a suit. 

Suppose $X$ is a pair of suits. If \(V' \subseteq V \in C_{n,J}(X)^+\), then \(V(n:f) \in X^+\) for some \(f\colon V \toind{J} A\). Let \(g = f\restrict V'\). Then \(g\colon V' \toind{J} A\), and \(V'(n:g) \subseteq V(n:f)\), so \(V'(n:g) \in X^+\). Hence \(V' \in C_{n,J}(X)^+\). Thus \(C_{n,J}(X)^+\) is a suit. If \(W' \subseteq W \in C_{n,J}(X)^-\), then \(W'(n:A) \subseteq W(n:A) \in X^-\), so \(W'(n:A) \in X^-\). Hence $C_{n,J}(X)^-$ is a suit. 
\end{proof}

\begin{proposition} % subalgebra generated double-suited
\label{subalgebra generated double-suited}
The subalgebra of an IFG$_N$-algebra generated by a set of double suits is a double-suited IFG$_N$-algebra.
\end{proposition}

\begin{proof}
The constant elements $0$, $1$, and $D_{ij}$ are all double suits, and if $X$ is a double suit, then so is $\n{X}$.

Suppose $X$ and $Y$ are double suits. We already know by the previous proposition that $X +_J Y$ is a pair of suits. To show that \((X +_J Y)^+ \cap (X +_J Y)^- = \set{\emptyset}\), suppose \(V \in (X +_J Y)^+ \cap (X +_J Y)^- = (X +_J Y)^+ \cap X^- \cap Y^-\). Then \(V = V_1 \cup_J V_2\) for some \(V_1 \in X^+\), \(V_2 \in Y^+\). Since \(V_1 \subseteq V \in X^-\) and \(V_2 \subseteq V \in Y^-\), we have \(V_1 \in X^+ \cap X^-\) and \(V_2 \in Y^+ \cap Y^-\). Hence \(V_1 = \emptyset = V_2\). Thus \(V = \emptyset\). Therefore \(X +_J Y\) is a double suit.

Suppose $X$ is a double suit. We already know by the previous proposition that $C_{n,J}(X)$ is a pair of suits. To show that \(C_{n,J}(X)^+ \cap C_{n,J}(X)^- = \set{\emptyset}\), suppose \(V \in C_{n,J}(X)^+ \cap C_{n,J}(X)^-\). Then \(V(n:f) \in X^+\) for some \(f\colon V \toind{J} A\), and \(V(n:A) \in X^-\). But \(V(n:f) \subseteq V(n:A)\), so \(V(n:f) \in X^+ \cap X^-\). Hence \(V(n:f) = \emptyset\), which implies \(V = \emptyset\). Therefore \(C_{n,J}(X)\) is a double suit.
\end{proof}

Given a set $A$, let $\Suit_N(A)$\index{$\Suit_N(A)$ \quad suited IFG$_N$-cylindric set algebra over $A$} denote the IFG$_N$-algebra whose universe is the set of all pairs of suits in \(\powerset(\powerset(\^NA)) \times \powerset(\powerset(\^NA))\). Let $\DSuit_N(A)$\index{$\DSuit_N(A)$ \quad double-suited IFG$_N$-cylindric set algebra over $A$} denote the IFG$_N$-algebra whose universe is the set of all double suits in \(\powerset(\powerset(\^NA)) \times \powerset(\powerset(\^NA))\).

% IFG$_N$-algebras and De Morgan algebra
\subsection{IFG$_N$-algebras and De Morgan algebra}
\label{subsec:IFG-algebras and De Morgan algebra}

\begin{definition} % bounded distributive lattice
A \emph{bounded distributive lattice}\index{bounded distributive lattice|mainidx} is an algebra \(\mathfrak L = \seq{L; 0,1,\join,\meet}\) such that 
\begin{alignat*}{2}
(x \join y) \join z &= x \join (y \join z),    &\qquad  
(x \meet y) \meet z  &= x \meet (y \meet z), \\
x \join y &= y \join x,		&\qquad
x \meet y &= y \meet x, \\
x \join (x \meet y) &= x,	&\qquad
x \meet (x \join y) &= x, \\
x \join (y \meet z) &= (x \join y) \meet (x \join z), &\qquad
x \meet (y \join z) &= (x \meet y) \join (x \meet z), \\
0 \join x &= x,		&\qquad
1 \meet x &= x.
\end{alignat*}
\end{definition}

We will refer to each pair of axioms except the last as associativity, commutativity, absorption, and distributivity, respectively. We can define a partial order $\leq$ on $\mathfrak L$ by \(x \leq y\) if and only if \(x \join y = y\) (if and only if \(x \meet y = x\)). Hence the last pair of axioms asserts that \(0 \leq x \leq 1\).

\begin{definition} % De Morgan algebra
A \emph{De Morgan algebra}\index{De Morgan algebra|mainidx} \(\A = \seq{A; 0, 1, {\hneg}\,, \join, \meet}\) is a bounded distributive lattice with an additional unary operation $\hneg\ $ that satisfies \(\hneg{\hneg x} = x\) and \(\hneg(x \join y) = {\hneg x} \meet {\hneg y}\).
\end{definition}

\begin{proposition} % De Morgan algebra properties
\label{De Morgan algebra properties}
Let $\A$ be a De Morgan algebra. Then $\A$ satisfies
\begin{enumerate}
	\item \(\hneg(x \meet y) = {\hneg x} \join {\hneg y}\),
	\item \(\hneg 0 = 1\),
	\item \(x \leq y\) if and only if \({\hneg y} \leq {\hneg x}\).
\end{enumerate}
\end{proposition}

\begin{proof}
(a)
\(
\hneg(x \meet y) = {\hneg({\hneg{\hneg x}} \meet {\hneg{\hneg y}})} = {\hneg{\hneg({\hneg x} \join {\hneg y})}} = {\hneg x} \join {\hneg y}.
\)

(b)
\(
\hneg 0 = {\hneg 0} \meet 1 = {\hneg 0} \meet {\hneg{\hneg 1}} = {\hneg(0 \join {\hneg 1})} = {\hneg{\hneg 1}} = 1.
\)

(c)
\(x \leq y\) if and only if \(x \join y = y\) if and only if \({\hneg x} \meet {\hneg y} = {\hneg y}\) if and only if \({\hneg y} \leq {\hneg x}\).
\end{proof}

\begin{definition} % Kleene algebra
A \emph{Kleene algebra}\index{Kleene algebra|mainidx} \(\A = \seq{A; 0, 1, \hneg\,, \join, \meet}\) is a De Morgan algebra that satisfies the additional axiom \(x \meet {\hneg x} \leq y \join {\hneg y}\).
\end{definition}

\begin{definition} % Boolean algebra
A \emph{Boolean algebra}\index{Boolean algebra|mainidx} \(\A = \seq{A; 0, 1, -, \join, \meet}\) is a De Morgan algebra that satisfies the complementation axioms \(x \meet -x = 0\) and \(x \join -x = 1\).
\end{definition}

We refer the reader to \cite{Balbes:1974} for the elementary theory of De Morgan algebras. 

Unlike cylindric set algebras, independence-friendly cylindric set algebras do not have an underlying Boolean algebra structure. The complementation axioms fail. However, the reduct of a rooted IFG$_N$-algebra to the signature 
\[
\tuple{0, 1, \n{}, +_N, \cdot_N}
\]
is a De Morgan algebra, and the same reduct of a double-suited IFG$_N$-algebra is a Kleene algebra.

% Associativity and commutativity
\subsubsection{Associativity and commutativity}

Let $\C$ be an IFG$_N$-algebra with base set $A$. Let \(X, Y, Z\) be elements in $\C$, let \(i,j,k, \ell, m, n < N\), let \(J, K, L \subseteq N\), and let \(V, W \subseteq \^NA\).

\begin{proposition} % commutativity
\label{commutativity}
\(X +_J Y = Y +_J X\) and \(X \cdot_J Y = Y \cdot_J X\).
\end{proposition}

\begin{proof}
\(V \in (X +_J Y)^+\) if and only if \(V = V_1 \cup_J V_2\) for some \(V_1 \in X^+\) and \(V_2 \in Y^+\) if and only if \(V = V_2 \cup_J V_1\) for some \(V_2 \in Y^+\) and \(V_1 \in X^+\) if and only if \(V \in (Y +_J X)^+\). Also, \((X +_J Y)^- = X^- \cap Y^- = Y^- \cap X^- = (Y +_J X)^-\).
\end{proof}

\begin{proposition} % associativity
\label{associativityJJ}
\((X +_J Y) +_J Z = X +_J (Y +_J Z)\) and \((X \cdot_J Y) \cdot_J Z = X \cdot_J (Y \cdot_J Z)\).
\end{proposition}

\begin{proof}
Suppose \(V \in ((X +_J Y) +_J Z)^+\). Then \(V = (V_1 \cup_J V_2) \cup_J V_3\) for some \(V_1 \in X^+\), \(V_2 \in Y^+\), and \(V_3 \in Z^+\). By \lemref{saturated union preserves saturation}, 
\begin{align*}
V &= (V_1 \cup_J V_2) \cup_J V_3 \\
&= V_1 \cup_J V_2 \cup_J V_3 \\
&= V_1 \cup_J (V_2 \cup_J V_3).
\end{align*}
Hence \(V \in (X +_J (Y +_J Z))^+\). By the same argument, if \(V \in (X +_J (Y +_J Z))^+\), then \(V \in ((X +_J Y) +_J Z)^+\). Therefore, \(((X +_J Y) +_J Z)^+ = (X +_J (Y +_J Z))^+\).

In addition, 
\begin{align*}
((X +_J Y) +_J Z)^- &= (X^- \cap Y^-) \cap Z^- \\
&= X^- \cap (Y^- \cap Z^-) \\
&= (X +_J (Y +_J Z))^-. \qedhere
\end{align*}
\end{proof}

Associativity can fail if the two operations are not the same. For example, let $\A$ be the structure with universe \(A = \set{0,1,2}\) in which each element is named by a constant symbol. Consider the IFG$_1$-algebra over $\A$ (that is, $\Cyls_{IFG_1}(\A)$). Let \(X = \norm{v_0 = 0}\), \(Y = \norm{v_0 = 1}\), and \(Z  = \norm{v_0 = 2}\). Then
\begin{align*}
X &= \tuple{\set{\emptyset, \set{0}}, \set{\emptyset, \set{1}, \set{2}, \set{1,2}}}, \\
Y &= \tuple{\set{\emptyset, \set{1}}, \set{\emptyset, \set{0}, \set{2}, \set{0,2}}}, \\
Z &= \tuple{\set{\emptyset, \set{2}}, \set{\emptyset, \set{0}, \set{1}, \set{0,1}}}.
\end{align*}

We claim that \((X +_\emptyset Y) +_N Z \not= X +_\emptyset (Y +_N Z)\). On the left, 
\begin{align*}
(X +_\emptyset Y)^+ &= \set{\emptyset, \set{0}, \set{1}, \set{0,1}}, \\
((X +_\emptyset Y) +_N Z)^+ &= \set{\emptyset, \set{0}, \set{1}, \set{0,1}, \set{2}},
\end{align*}
while on the right, 
\begin{align*}
(Y +_N Z)^+ &= \set{\emptyset, \set{1}, \set{2}}, \\
(X +_\emptyset (Y +_N Z))^+ &= \set{\emptyset, \set{0}, \set{1}, \set{2}, \set{0,1}, \set{0,2}}.
\end{align*}
Of course, on the falsity-axis, \(((X +_\emptyset Y) + _N Z)^- = X^- \cap Y^- \cap Z^- = (X +_\emptyset (Y +_N Z))^-\).

Notice that \(((X +_\emptyset Y) +_N Z)^+ \subseteq (X +_\emptyset (Y +_N Z))^+\). This suggests the following strengthening of \propref{associativityJJ}.

\begin{lemma} % associativityJK
\label{associativityJK}
Let \(J \subseteq K\). Then
\begin{enumerate}
	\item  \((X +_J Y) +_K Z \leq^+ X +_J (Y +_K Z)\) and \((X +_J Y) +_K Z =^- X +_J (Y +_K Z)\);
	\item  \((X \cdot_J Y) \cdot_K Z =^+ X \cdot_J (Y \cdot_K Z)\) and \((X \cdot_J Y) \cdot_K Z \leq^- X \cdot_J (Y \cdot_K Z)\).
\end{enumerate}
\end{lemma}

\begin{proof}
(a)
Suppose \(V \in ((X +_J Y) +_K Z)^+\). Then \(V = (V_1 \cup_J V_2) \cup_K V_3\) for some \(V_1 \in X^+\), \(V_2 \in Y^+\), and \(V_3 \in Z^+\). By \lemref{K saturated J} and \lemref{saturated union preserves saturation},
\begin{align*}
V &= (V_1 \cup_J V_2) \cup_K V_3 \\
&= (V_1 \cup_J V_2) \cup_J V_3 \\
&= V_1 \cup_J V_2 \cup_J V_3 \\
&= V_1 \cup_J (V_2 \cup_J V_3).
\end{align*}
Let \(V_5 = V_2 \cup V_3\). Since $V_1$, $V_2$, and $V_3$ are pairwise disjoint, we have by \lemref{restricted cover preserves saturation} that \(V_5 = V \cap V_5 = ((V_1 \cup_J V_2) \cap V_5) \cup_K (V_3 \cap V_5) = V_2 \cup_K V_3\). Hence \(V_5 \in (Y +_K Z)^+\). Thus \(V \in (X +_J (Y +_K Z))^+\).

(b)
This is just the dual of part (a).
\end{proof}

% The elements 0, $\Omega$, $\mho$, and 1
\subsubsection{The elements 0, $\Omega$, $\mho$, and 1}

\begin{proposition} % 0, Omega, mho, 1
\label{0, Omega, mho, 1}
\begin{enumerate}
	\item \(X +_J 0 = X = X \cdot_J 1\).
	\item \(X \cdot_J 0 = 0\) and \(X +_J 1 = 1\) if and only if $X$ is rooted.
	\item \(\n{\Omega} = \Omega +_J \Omega = \Omega \cdot_J \Omega = \Omega\).
	\item \(\n{\mho} = \mho +_J \mho = \mho \cdot_J \mho = \mho\).
\end{enumerate}
\end{proposition}
\begin{proof}
(a)
Suppose \(V \in (X +_J 0)^+\). Then \(V = V_1 \cup_J V_2\) for some \(V_1 \in X^+\) and \(V_2 \in 0^+\). But then \(V_2 = \emptyset\), so \(V_1 = V\). Thus \(V \in X^+\). Conversely, suppose \(V \in X^+\). Then \(V = V \cup_J \emptyset\), \(V \in X^+\), and \(\emptyset \in 0^+\). Thus \(V \in (X +_J 0)^+\). Therefore \((X +_J 0)^+ = X^+\). Also, \((X +_J 0)^- = X^- \cap\, 0^- = X^-\).

(b)
Suppose \(X \cdot_J 0 = 0\) and \(X +_J 1 = 1\). Then \(\emptyset \in (X \cdot_J 0)^+ = X^+ \cap 0^+ \subseteq X^+\), and \(\emptyset \in (X +_J 1)^- = X^- \cap 1^- \subseteq X^-\). Hence $X$ is rooted. Conversely, suppose $X$ is rooted. Then \((X \cdot_J 0)^+ = X^+ \cap 0^+ = \set{\emptyset}\), and for every \(V \subseteq \^NA\), \(V = \emptyset \cup_J V\), where \(\emptyset \in X^-\) and \(V \in 0^-\). Hence \(V \in (X \cdot_J 0)^-\). Thus \((X \cdot_J 0)^- = \powerset(\^NA)\). Therefore \(X \cdot_J 0 = 0\). Similarly, \(X +_J 1 = 1\).

(c)
That \(\n{\Omega} = \Omega\) follows immediately from the definitions. To show that \(\Omega +_J \Omega = \Omega\), observe that \(\emptyset = \emptyset \cup_J \emptyset\). Hence \(\emptyset \in (\Omega +_J \Omega)^+\). Conversely, if \(V \in (\Omega +_J \Omega)^+\), then \(V = V_1 \cup_J V_2\) for some \(V_1, V_2 \in \Omega^+\). Hence \(V = V_1 = V_2 = \emptyset\). Thus \((\Omega +_J \Omega)^+ = \Omega^+\). Also \((\Omega +_J \Omega)^- = \Omega^- \cap \Omega^- = \Omega^-\).

(d)
That \(\n{\mho} = \mho\) follows immediately from the definitions. To show that \(\mho +_J \mho = \mho\), observe that for any \(V \subseteq \^NA\), \(V = V \cup_J \emptyset\). Hence \(V \in (\mho +_J \mho)^+\). Thus \((\mho +_J \mho)^+ = \mho^+\). Also \((\mho +_J \mho)^- = \mho^- \cap \mho^- = \mho^-\).
\end{proof}

\begin{proposition} % X<Omega<Y
\label{X<Omega<Y}
If $X$ and $Y$ are double suits, and \(X \leq \Omega \leq Y\), then \(X \cdot_J Y = X\) and \(X +_J Y = Y\).
\end{proposition}

\begin{proof}
Suppose \(X \leq \Omega \leq Y\). Then \(X^+ = \set{\emptyset} = Y^-\), so \((X \cdot_J Y)^+ = X^+ \cap Y^+ = \set{\emptyset} = X^- \cap Y^- = (X +_J Y)^-\). Also, \(V \in (X +_J Y)^+\) if and only if \(V \in Y^+\), and \(W \in (X \cdot_J Y)^-\) if and only if \(W \in X^-\). Therefore \(X \cdot_J Y = X\) and \(X +_J Y = Y\).
\end{proof}

% Absorption
\subsubsection{Absorption}

The absorption axioms do not hold in general, but they do hold in important special cases. For example, if $X$ is not rooted then \((X +_N X)^+ = \emptyset = (X \cdot_N X)^-\). Hence \(X +_N (X \cdot_K 1) \not= X\) whenever $X$ is a nonempty, non-rooted element. In contrast, if $X$ and $Y$ are rooted, then the absorption axioms hold partially for every pair of addition and multiplication operations, and they hold fully when the ``outside'' operation is $+_N$ or $\cdot_N$.

\begin{lemma} % rooted +
\label{rooted +}
If $Y$ is rooted, then \(X^+ \subseteq (X +_J Y)^+\) and \(X ^- \subseteq (X \cdot_J Y)^-\).
\end{lemma}

\begin{proof}
Suppose \(V \in X^+\). Then \(V = V \cup_J \emptyset\), where \(V \in X^+\) and \(\emptyset \in Y^+\). Hence \(V \in (X +_J Y)^+\).
\end{proof}

\begin{lemma} % +_N and \cdot_N
If $X$ and $Y$ are rooted, then \((X +_N Y)^+ = X^+ \cup Y^+\) and \((X \cdot_N Y)^- = X^- \cup Y^-\).
\end{lemma}

\begin{proof}
By definition, \(V \in (X +_N Y)^+\) if and only if \(V = V_1 \cup_N V_2\) for some \(V_1 \in X^+\) and \(V_2 \in Y^+\), which holds if and only if \(V_1 = V\) and \(V_2 = \emptyset\) or vice versa. Hence \(V \in (X +_N Y)^+\) if and only if \(V \in X^+\) or \(V \in Y^+\).
\end{proof}

\begin{lemma} % absorption_JK
\label{absorption_JK}
Suppose $X$ and $Y$ are rooted. Then 
\begin{enumerate}
	\item \(X \leq^+ X +_J (X \cdot_K Y)\) and \(X = ^- X +_J (X \cdot_K Y)\);
	\item \(X =^+ X \cdot_J (X +_K Y)\) and \(X \leq^- X \cdot_J (X +_K Y)\).
\end{enumerate}
Thus \(X \cdot_J (X +_K Y) \leq X \leq X +_J (X \cdot_K Y)\).
\end{lemma}

\begin{proof}
(a)
Suppose \(V \in X^+\). Then \(V = V \cup_J \emptyset\) where \(V \in X^+\) and \(\emptyset \in X^+ \cap Y^+ = (X \cdot_K Y)^+\). Hence \(V \in (X +_J (X \cdot_K Y))^+\).

Suppose \(W \in X^-\). Then \(W = W \cup_K \emptyset\), where \(W \in X^-\) and \(\emptyset \in Y^-\). Hence \(W \in X^- \cap (X \cdot_K Y)^- = (X +_J (X \cdot_K Y))^-\). Conversely, suppose \(W \in (X +_J (X \cdot_K Y))^-\). Then \(W \in X^- \cap (X \cdot_K Y )^- \subseteq X^-\).
\end{proof}

To show that absorption can fail even when $X$ and $Y$ are rooted, consider the IFG$_2$-cylindric set algebra over the equality structure $\A$ whose universe is \(A = \set{0,1}\). Let \(X = D_{01} +_N \n{D_{01}}\). Then 
\begin{align*}
X^+ &= \set{\emptyset, \set{00}, \set{01}, \set{10}, \set{11}, \set{00, 11}, \set{01, 10}}, \\
X^- &= \set{\emptyset}.
\end{align*}
However, \(\set{00, 01, 10, 11} \in (X +_\emptyset (X +_N X))^+\) because \(\set{00, 01, 10, 11} = \set{00, 11} \cup_\emptyset \set{01, 10}\), where \(\set{00, 11} \in X^+\) and \(\set{01, 10} \in X^+ = (X +_N X)^+\). Thus 
\[
(X +_\emptyset (X +_N X))^+ \supset X^+.
\]

One can obtain a similar example where \((X \cdot_\emptyset (X +_N X))^- \supset X^-\) by taking \(X = D_{01} \cdot_N \n{D_{01}}\).

Notice that \(D_{01} +_N \n{D_{01}} \not= 1\), which demonstrates that the complementation axioms fail in any IFG-algebra of dimension greater than 1.

\begin{lemma} % flat absorption
\label{flat absorption}
If $X$ is flat, then \(X +_J (X \cdot_K Y) = X\). If $\n{X}$ is flat, then \(X \cdot_J (X +_K Y) = X\).
\end{lemma}

\begin{proof}
Suppose \(X^+ = \powerset(V)\). If \(V' \in (X +_J (X \cdot_K Y))^+\), then \(V' = V_1 \cup_J V_2\) for some \(V_1 \in X^+\) and \(V_2 \in (X \cdot_K Y)^+ = X^+ \cap Y^+\). Hence \(V' = V_1 \cup_J V_2 \subseteq V\). Thus \(V' \in \powerset(V) = X^+\). 
\end{proof}

\begin{lemma} % absorption_NJ
\label{absorption_NJ}
Suppose $X$ and $Y$ are rooted. Then
\begin{enumerate}
	\item \(X +_N (X \cdot_J Y) = X\);
	\item \(X \cdot_N (X +_J Y) = X\).
\end{enumerate}
\end{lemma}

\begin{proof}
(a) 
First, observe that \((X +_N (X \cdot_J Y))^+ = X^+ \cup (X^+ \cap Y^+) = X^+\). Next, suppose \(W \in (X +_N (X \cdot_J Y))^-\). Then \(W \in X^- \cap (X \cdot_J Y)^- \subseteq X^-\). Conversely, suppose \(W \in X^-\). Then \(W = W \cup_J \emptyset\) and \(W \in X^-\), \(\emptyset \in Y^-\). Hence \(W \in (X \cdot_J Y)^-\). Thus \(W \in X^- \cap (X \cdot_J Y)^- = (X +_N (X \cdot_J Y))^-\).
\end{proof}

Thus, the reduct of a rooted IFG$_N$-algebra to the signature $\tuple{+_N, \cdot_N}$ is a lattice. As with all lattices, we can define a partial order $\leq'$ by declaring \(X \leq' Y\) if and only if \(X +_N Y = Y\) (or equivalently, \(X \cdot_N Y = X\)). The next proposition shows that in a rooted IFG$_N$-algebra, our two partial orders $\leq$ and $\leq'$ agree.

\begin{proposition} % less than
\label{less than}
Suppose $X$ and $Y$ are rooted. Then \(X \leq Y\) if and only if \(X +_N Y = Y\) if and only if \(X \cdot_N Y = X\).
\end{proposition}
\begin{proof}
First, \(X \leq Y\) if and only if \(X^+ \subseteq Y^+\) and \(Y^- \subseteq X^-\) if and only if \(X^+ \cup Y^+ = Y^+\) and \(X^- \cap Y^- = Y^-\) if and only if \((X +_N Y)^+ = Y^+\) and \((X +_N Y)^- = Y^-\).

Second, \(X \leq Y\) if and only if \(X^+ \subseteq Y^+\) and \(Y^- \subseteq X^-\) if and only if \(X^+ \cap Y^+ = X^+\) and \(X^- \cup Y^- = X^-\) if and only if \((X \cdot_N Y)^+ = X^+\) and \((X \cdot_N Y)^- = X^-\).
\end{proof}

\begin{proposition} % 0 < X < 1
\label{0 < X < 1}
If $X$ is rooted, then \(0 \leq X \leq 1\).
\end{proposition}
\begin{proof}
By \propref{0, Omega, mho, 1} and \propref{less than}.
\end{proof}

\begin{proposition} % +_K less than +_J
\label{+_K less than +_J}
If \(J \subseteq K\), then \(X +_K Y \leq X +_J Y\) and \(X \cdot_J Y \leq X \cdot_K Y\).
\end{proposition}
\begin{proof}
Suppose \(J \subseteq K\). If \(V \in (X +_K Y)^+\), then \(V = V_1 \cup_K V_2\) for some \(V_1 \in X^+\) and \(V_2 \in Y^+\). But \(V = V_1 \cup_K V_2\) implies \(V = V_1 \cup_J V_2\), so \(V \in (X +_J Y)^+\). Thus \((X +_K Y)^+ \subseteq (X +_J Y)^+\). In addition, \((X +_J Y)^- = X^- \cap Y^- = (X +_K Y)^-\).
\end{proof}

Therefore, when \(X\) is rooted, we have the following string of inequalities:
\[
X \cdot_\emptyset X \leq X \cdot_J X \leq X \cdot_N X = X = X +_N X \leq X +_J X \leq X +_\emptyset X.
\]

\begin{proposition} % + preserves order
\label{+ preserves order}
If \(X \leq X'\) and \(Y \leq Y'\), then \(X +_J Y \leq X' +_J Y'\) and \(X \cdot_J Y \leq X' \cdot_J Y'\).
\end{proposition}
\begin{proof}
Suppose \(X \leq X'\) and \(Y \leq Y'\). If \(V \in (X +_J Y)^+\), then \(V = V_1 \cup_J V_2\) for some \(V_1 \in X^+ \subseteq (X')^+\) and \(V_2 \in Y^+ \subseteq (Y')^+\). Hence \(V \in (X' +_J Y')^+\). Also, \((X' +_J Y')^- = (X')^- \cap (Y')^- \subseteq X^- \cap Y^- = (X +_J Y)^-\). Thus \(X +_J Y \leq X' +_J Y'\).
\end{proof}

\subsubsection{Distributivity}

\begin{lemma} % distributivity_JK
\label{distributivity_JK}
Suppose $X$ is a double suit. Then 
\begin{enumerate}
	\item \((X \cdot_J (Y +_K Z))^\pm \subseteq ((X \cdot_J Y) +_K (X \cdot_J Z))^\pm\);
	\item \((X +_J (Y \cdot_K Z))^\pm \subseteq ((X +_J Y) \cdot_K (X +_J Z))^\pm\).
\end{enumerate}
\end{lemma}

\begin{proof}
(a)
Suppose \(V \in (X \cdot_J (Y +_K Z))^+\). Then \(V \in X^+ \cap (Y +_K Z)^+\), which implies that \(V = V_1 \cup_K V_2\) for some \(V_1 \in Y^+\) and \(V_2 \in Z^+\). Observe that since \(V_1, V_2 \subseteq V \in X^+\) we have \(V_1, V_2 \in X^+\). Hence \(V_1 \in X^+ \cap Y^+ =  (X \cdot_J Y)^+\) and \(V_2 \in X^+ \cap Z^+ = (X \cdot_J Z)^+\). Thus \(V \in ((X \cdot_J Y) +_K (X \cdot_J Z))^+\). 

Now suppose \(W \in (X \cdot_J (Y +_K Z))^-\). Then \(W = W_1 \cup_J W_2\) for some \(W_1 \in X^-\) and \(W_2 \in Y^- \cap Z^-\). It follows that \(W \in (X \cdot_J Y)^- \cap (X \cdot_J Z)^- = ((X \cdot_J Y) +_K (X \cdot_J Z))^-\).
\end{proof}

To show that distributivity can fail, consider the IFG$_2$-cylindric set algebra over the structure $\A$ whose universe is \(A = \set{0, 1}\) and in which each element is named by a constant. Let \(X = \norm{v_0 = 0} +_N \norm{v_0 = 1}\), and let 
\begin{align*}
V &= \set{00, 01, 10, 11}, \\
V_1 &= \set{00, 01}, \\
V_2 &= \set{10, 11}.
\end{align*}
Observe that \(V = V_1 \cup_{\set{1}} V_2\), \(V_1 \in X^+\) and \(V_2 \in X^+\). Hence
\[
V \in (X +_{\set{1}} X)^+ = ((X \cdot_J 1) +_{\set{1}} (X \cdot_J 1))^+.
\]
However, \(V \notin X^+ = (X \cdot_J (1 +_{\set{1}} 1))^+\).

\begin{lemma} % distributivity_JN
\label{distributivity_JN}
Suppose $X$, $Y$, and $Z$ are rooted.
\begin{enumerate}
	\item \(X \cdot_J (Y +_N Z) =^+ (X \cdot_J Y) +_N (X \cdot_J Z)\).
	\item \(X \cdot_N (Y +_K Z) =^- (X \cdot_N Y) +_K (X \cdot_N Z)\).
	\item \(X +_N (Y \cdot_K Z) =^+ (X +_N Y) \cdot_K (X +_N Z)\).
	\item \(X +_J (Y \cdot_N Z) =^- (X +_J Y) \cdot_N (X +_J Z)\).
\end{enumerate}
\end{lemma}

\begin{proof}
\begin{align*}
(X \cdot_J (Y +_N Z))^+ &= X^+ \cap (Y^+ \cup Z^+) \\
&= (X^+ \cap Y^+) \cup (X^+ \cap Z^+) \\
&= ((X \cdot_J Y) +_N (X \cdot_J Z))^+. \\
\\
(X \cdot_N (Y +_K Z))^- &= X^- \cup (Y^- \cap Z^-) \\
&= (X^- \cup Y^-) \cap (X^- \cup Z^-) \\
&= ((X \cdot_N Y) +_K (X \cdot_N Z))^-. \\
\\
(X +_N (Y \cdot_K Z))^+ &= X^+ \cup (Y^+ \cap Z^+) \\
&= (X^+ \cup Y^+) \cap (X^+ \cup Z^+) \\
&= ((X +_N Y) \cdot_K (X +_N Z))^+. \\
\\
(X +_J (Y \cdot_N Z))^- &= X^- \cap (Y^- \cup Z^-) \\
&= (X^- \cap Y^-) \cup (X^- \cap Z^-) \\
&= ((X +_J Y) \cdot_N (X +_J Z))^-. \qedhere
\end{align*}
\end{proof}

The previous lemmas and propositions combine to yield the following theorem.

\begin{theorem} % rooted -> De Morgan
\label{rooted -> De Morgan}
The reduct of a rooted IFG$_N$-algebra to the signature $\seq{0, 1, \n{}, +_N, \cdot_N}$ is a De Morgan algebra\index{De Morgan algebra}.
\end{theorem}
\index{IFG-cylindric set algebra!rooted}

\subsubsection{Complementation}

We have already remarked that the complementation axioms fail in independence-friendly cylindric set algebras. The failure of the complementation axioms stems from the failure of the law of excluded middle in IFG logic. However sentences of the form \(\phi \hand{J} {\hneg \phi}\) are never true, while sentences of the form \(\phi \hor{K} {\hneg\phi}\) are never false. 

\begin{proposition} % X*-X < Y+-Y
\label{X*-X < Y+-Y}
If $X$ and $Y$ are double suits, then \(X \cdot_J \n{X} \leq Y +_K \n{Y}\).
\end{proposition}

\begin{proof}
Suppose $X$ and $Y$ are double suits. Then \((X \cdot_J \n{X})^+ = X^+ \cap X^- = \set{\emptyset} \subseteq (Y +_K \n{Y})^+\), and 
\((Y +_K \n{Y})^- = Y^+ \cap Y^- = \set{\emptyset} \subseteq (X +_J \n{X})^-\).
\end{proof}

Thus, if a double-suited IFG$_N$-cylindric set algebra includes $\Omega$, then 
\[
X \cdot_J \n{X} \leq \Omega \leq X +_K \n{X}.
\]

\begin{theorem} % double-suited -> Kleene algebra
\label{double-suited -> Kleene algebra}
\index{IFG-cylindric set algebra!double-suited}
\index{Kleene algebra}
The reduct of a double-suited IFG$_N$-algebra to the signature $\seq{0, 1, \n{}, +_N, \cdot_N}$ is a Kleene algebra.
\end{theorem}

Even though the complementation axioms do not hold universally, it might still be possible for particular elements to satisfy them. In a double-suited IFG$_N$-algebra, the only elements that satisfy the complementation axioms are 0 and 1. 

\begin{lemma} % complementation
\label{complemetation}
\(X +_N \n{X} = 1\) (equivalently, \(X \cdot_N \n{X} = 0\)) if and only if \(X^+ \cup X^- = \powerset(\^NA)\) and \(X^+ \cap X^- = \set{\emptyset}\).
\end{lemma}

\begin{proof}
Suppose \(X +_N \n{X} = 1\). Then \(X^- \cap X^+ = (X +_N \n{X})^- = \set{\emptyset}\). Hence $X$ is rooted. Thus \(X^+ \cup X^- = (X +_N \n{X})^+ = \powerset(\^NA)\).

Conversely, suppose  \(X^+ \cup X^- = \powerset(\^NA)\) and \(X^+ \cap X^- = \set{\emptyset}\). Then $X$ is rooted, so \((X +_N \n{X}) = X^+ \cup X^- = 1^+\) and \((X +_N \n{X})^- = X^- \cap X^+ = 1^-\).
\end{proof}

\begin{lemma} % suited complementation
\label{suited complemetation}
Suppose $X$ is a double suit. Then \(X +_N \n X = 1\) if and only if \(X = 0\) or \(X = 1\).
\end{lemma}

\begin{proof}
Suppose \(X +_N \n X = 1\). Then, by the previous lemma, \((X +_N \n X)^+ = X^+ \cup X^- = \powerset(\^NA)\). In particular, \(\^NA \in X^+ \cup X^-\). If \(\^NA \in X^+\) then, since $X^+$ is a suit, \(X^+ = \powerset(\^NA)\) and \(X^- = \set{\emptyset}\). Hence \(X = 1\). Similarly, if \(\^NA \in X^-\), then \(X = 0\). 

Conversely, suppose \(X = 0\) or \(X = 1\). In either case, \(X +_N \n X = 0 +_N 1 = 1\).
\end{proof}

\begin{definition} % complement
Let $L$ be a bounded lattice. Two elements \(x,y \in L\) are \emph{complements}\index{complements|mainidx} if \(x \join y = 1\) and \(x \meet y = 0\). An element \(x \in L\) is \emph{complemented}\index{complemented element|mainidx} if it has a complement.
\end{definition}

For example, $\Omega$ and $\mho$ are complements in any IFG-algebra that includes both. 
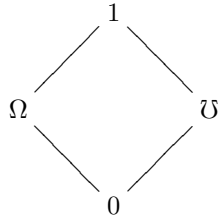
\[
\xymatrix{
					& 1 \ar@{-}[dl] \ar@{-}[dr]	&					\\
	\Omega \ar@{-}[dr] 		& 					& \mho \ar@{-}[dl]			\\
					& 0  &			\\
}
\]
In a double-suited IFG-algebra the only complemented elements are 0 and 1.
\vspace{0 pt}

\begin{lemma} % not complemented
\label{not complemented}
Suppose \(X +_N Y = 1\) and \(X \cdot_N  Y = 0\). Then  \(X^+ \cup Y^+ = \powerset(\^NA) = X^- \cup Y^-\) and \(X^+ \cap Y^+ = \set{\emptyset} = X^- \cap Y^-\).
\end{lemma}

\begin{proof}
First, \(X^+ \cap Y^+ = (X \cdot_N Y)^+ = \set{\emptyset} = (X +_N Y)^- = X^- \cap Y^-\). Hence $X$ and $Y$ are both rooted. Thus \(X^+ \cup Y^+ = (X +_N Y)^+ = \powerset(\^NA) = (X \cdot_N Y)^- = X^- \cup Y^-\).
\end{proof}

\begin{proposition} % suited not complemented
\label{suited not complemented}
Suppose $X$ and $Y$ are double suits. Then \(X +_N Y = 1\) and \(X \cdot_N  Y = 0\) if and only if \(X = 1\) and \(Y = 0\), or vice versa.
\end{proposition}

\begin{proof}
Suppose \(X +_N Y = 1\) and \(X \cdot_N Y = 0\). By the previous lemma, \(X^+ \cup Y^+ = \powerset(\^NA) = X^- \cup Y^-\) and \(X^+ \cap Y^+ = \set{\emptyset} = X^- \cap Y^-\). In particular, \(\^NA \in X^+ \cup Y^+\). Since $X$ and $Y$ are both double suits, if \(\^NA \in X^+\) then \(X = 1\), and if \(\^NA \in Y^+\) then \(Y = 1\). In the first case, \(X^+ \cap Y^+ = \set{\emptyset}\) and \(X^- \cup Y^- = \powerset(\^NA)\) imply that \(Y = 0\), and in the second case, that \(Y = 1\).

Conversely, suppose \(X = 1\) and \(Y = 0\), or vice versa. In either case, \(X +_N Y = 0 +_N 1 = 1\) and \(X \cdot_N Y = 0 \cdot_N 1 = 0\).
\end{proof}

%% IFG$_N$-algebras and cylindric algebra
\subsection{IFG$_N$-algebras and cylindric algebra}

In addition to the axioms of Boolean algebra, the axioms of cylindric algebra\index{cylindric algebra} \cite{Henkin:1971} are:

\begin{itemize}
	\item[(C1)] \(c_n(0) = 0\).
	\item[(C2)] \(x \cdot c_n(x) = x\).
	\item[(C3)] \(c_n(x \cdot c_n(y)) = c_n(x) \cdot c_n(y)\).
	\item[(C4)] \(c_m c_n(x) = c_n c_m(x)\).
	\item[(C5)] \(c_i(d_{ij}) = 1\).
	\item[(C6)] \(c_k(d_{ik} \cdot d_{kj}) = d_{ij}\) provided \(k \not= i,j\).
	\item[(C7)] \(c_i(d_{ij} \cdot x) \cdot c_i(d_{ij} \cdot -x) = 0\) provided \(i \not= j\).
\end{itemize}

We will show that analogs of the axioms (C1)--(C6) hold in all double-suited independence-friendly cylindric set algebras. Axiom (C7) does not hold in general, even in double-suited independence-friendly cylindric set algebras.

\subsubsection{Axiom C1}

\begin{proposition} % C(0) = 0
\label{C(0) = 0}
\(C_{n,J}(0)  = 0\) and \(\Cd{n,J}(1) = 1\).
\end{proposition}
\begin{proof}
Suppose \(V \in C_{n,J}(0)^+\). Then \(V(n:f) \in 0^+\) for some \(f\colon V \toind{J} A\). But then \(V(n:f) = \emptyset\), which holds if and only if \(V = \emptyset\). Thus \(V \in 0^+\). Conversely, suppose \(V \in 0^+\). Then \(V = \emptyset\). Observe that the empty function $f$ from $\emptyset$ to $A$ is vacuously independent of $J$ and \(\emptyset(n:f) = \emptyset \in 0^+\). Thus \(\emptyset \in C_{n,J}(0)^+\). Therefore \(C_{n,J}(0)^+ = 0^+\).

\(C_{n,J}(0)^- \subseteq \powerset(\^NA) = 0^-\) is immediate. To show \(0^- \subseteq C_{n,J}(0)^-\), suppose \(W \in 0^-\). Then \(W \subseteq \^NA\), so \(W(n:A) \subseteq \^NA\). Hence \(W(n:A) \in 0^-\). Thus \(W \in C_{n,J}(0)^-\). Therefore \(C_{n,J}(0)^- = 0^-\).
\end{proof}

\begin{proposition} % C(X) = 0 iff X = 0
\label{C(X) = 0 iff X = 0}
If $X$ is a double suit, then \(C_{n,J}(X) = 0\) if and only if \(X = 0\).
\end{proposition}

\begin{proof}
Observe that \(\^NA \in C_{n,J}(X)^-\) if and only if \(\^NA = \^NA(n:A) \in X^-\). Thus if \(X \not= 0\), then \(\^NA \notin X^-\), so \(\^NA \notin C_{n,J}(X)^-\), hence \(C_{n,J}(X) \not= 0\).
\end{proof}

\begin{proposition} % C(Omega) = Omega
\label{C(Omega) = Omega}
\(C_{n,J}(\Omega) = \Omega = \Cd{n,J}(\Omega)\) and \(C_{n,J}(\mho) = \mho = \Cd{n,J}(\mho)\).
\end{proposition}

\begin{proof}
\(\emptyset \in C_{n,J}(\Omega)^+\) because $\Omega$ is rooted. Conversely, suppose \(V \in C_{n,J}(\Omega)^+\). Then \(V(n:f) \in \Omega^+\) for some \(f\colon V \toind{J} A\). Hence \(V(n:f) = \emptyset\), which implies \(V = \emptyset\). Thus \(C_{n,J}(\Omega)^+ = \Omega^+\). Also, \(W \in C_{n,J}(X)^-\) if and only if \(W(n:A) \in \Omega^-\) if and only if \(W = \emptyset\). Hence \(C_{n,J}(X)^- = \Omega^-\).

We know that $C_{n,J}(\mho)$ is a pair of suits because $\mho$ is. Thus it suffices to show that \(\^NA \in C_{n,J}(\mho)^+\) and \(\^NA \in C_{n,J}(\mho)^-\). Observe that for any \(f\colon \^NA \toind{J} A\) we have \(\^NA(n:f) \subseteq \^NA \in \mho^+\). Hence \(\^NA \in C_{n,J}(\mho)^+\). Also, \(\^NA(n:A) = \^NA \in \mho^-\). Hence \(\^NA \in C_{n,J}(\mho)^-\).
\end{proof}

\begin{proposition} % C(1) = 1
\label{C(1) = 1}
\(C_{n,J}(1) = 1\) and \(\Cd{n,J}(0) = 0\).
\end{proposition}
\begin{proof}
That \(C_{n,J}(1)^+ \subseteq 1^+\) is immediate. To show \(1^+ \subseteq C_{n,J}(1)^+\), suppose \(V \in 1^+\). Define \(f\colon V \toind{J} A\) by \(f(\vec a) = c\) for some arbitrary \(c \in A\). Then \(V(n:f) \in 1^+\), so \(V \in C_{n,J}(1)^+\). Therefore \(C_{n,J}(1)^+ = 1^+\).

That \(1^- \subseteq C_{n,J}(1)^-\) is immediate. To show \(C_{n,J}(1)^- \subseteq 1^-\), suppose \(W \in C_{n,J}(1)^-\). Then \(W(n:A) \in 1^-\). Hence \(W(n:A) = \emptyset\), which implies \(W = \emptyset\). Thus \(W \in 1^-\). Therefore \(C_{n,J}(1)^- = 1^-\).
\end{proof}

\subsubsection{Axiom C2}

\begin{proposition} % X < C(X)
\begin{enumerate}
	\item If \(n \notin J\), then \(X \cdot_K C_{n,J}(X) =^+ X\).
	\item If $X^-$ is a suit, then \(X \cdot_N C_{n,J}(X) =^- X\). 
\end{enumerate}
\noindent In particular, if $X$ is a double suit then \(X \cdot_N C_{n,\emptyset}(X) = X\).
\end{proposition}

\begin{proof}
(a)
Observe that \((X \cdot_K C_{n,J}(X))^+ = X^+ \cap C_{n,J}(X)^+ \subseteq X^+\). Conversely, suppose \(V \in X^+\). If \(n \notin J\), the projection \(\pr_n\colon V \to A\) that maps $\vec a$ to $a_n$ is independent of $J$, and \(V(n:\pr_n) = V \in X^+\). Hence \(V \in X^+ \cap C_{n,J}(X)^+ = (X \cdot_K C_{n,J}(X))^+\).

(b)
Suppose $X^-$ is a suit and \(W \in (X \cdot_N C_{n,J}(X))^-\). Then \(W = W \cup_N \emptyset\), where \(W \in X^-\) and \(\emptyset \in C_{n,J}(X)^-\), or vice versa. If \(W \in X^-\) we are done, so suppose \(W \in C_{n,J}(X)^-\). Then \(W \subseteq W(n:A) \in X^-\); hence \(W \in X^-\). Conversely, suppose \(W \in X^-\). Since $X^-$ is a suit, $C_{n,J}(X)^-$ is also a suit. In particular, \(\emptyset \in C_{n,J}(X)^-\). Hence \(W = W \cup_N \emptyset\), where \(W \in X^-\) and \(\emptyset \in C_{n,J}(X)^-\). Thus \(W \in (X \cdot_N C_{n,J}(X))^-\).
\end{proof}

To give an example where \(n \in J\) and \(X \cdot_K C_{n,J}(X) \not=^+ X\), let $\A$ be the equality structure with universe \(A = \set{0,1}\), and consider \(\Cyls_{\mathrm{IFG}_2}(\A)\). Observe that 
\begin{align*}
D_{01} &= \tuple{\powerset(\set{00, 11}), \powerset(\set{01, 10})}, \\
C_{0,N}(D_{01}) &= \tuple{\set{\emptyset, \set{00}, \set{11}}, \set{\emptyset}}, \\
(D_{01} \cdot_K C_{0,N}(D_{01}))^+ &= \set{\emptyset, \set{00}, \set{11}} \not = D_{01}^+.
\end{align*}

In section \ref{subsec:IFG-algebras and De Morgan algebra} we showed that the partial order $\leq$ interacts nicely with the operations $\n{\,}$, $+_J$, and $\cdot_J$ (see \propref{+_K less than +_J} and \propref{+ preserves order}). It interacts equally well with the operations $C_{n,J}$ and $\Cd{n,J}$. In particular, the previous proposition shows that in a double-suited IFG$_N$-algebra, \(\Cd{n,J}(X) \leq X \leq C_{n,J}(X)\) whenever \(n \notin J\).

% Axiom C3
\subsubsection{Axiom C3}

\begin{proposition} % C_K < C_J
\label{C_K < C_J}
If \(J \subseteq K\), then \(C_{n,K}(X) \leq C_{n,J}(X)\) and \(\Cd{n,J}(X) \leq \Cd{n,K}(X)\).
\end{proposition}

\begin{proof}
Suppose \(J \subseteq K\) and \(V \in C_{n,K}(X)^+\). Then \(V(n:f) \in X^+\) for some function \(f\colon V \toind{K} A\). But then \(f\colon V \toind{J} A\); hence \(V \in C_{n,J}(X)^+\). Also \(C_{n,J}(X)^- = C_{n,K}(X)^-\) by definition.
\end{proof}

\begin{proposition} % C preserves <
\label{C preserves <}
If \(X \leq Y\) then \(C_{n,J}(X) \leq C_{n,J}(Y)\) and \(\Cd{n,J}(X) \leq \Cd{n,J}(Y)\).
\end{proposition}

\begin{proof}
Suppose \(X \leq Y\) and \(V \in C_{n,J}(X)^+\). Then \(V(n:f) \in X^+ \subseteq Y^+\) for some \(f\colon V \toind{J} A\). Hence \(V \in C_{n,J}(Y)^+\). Now suppose \(W \in C_{n,J}(Y)^-\). Then \(W(n:A) \in Y^- \subseteq X^-\). Hence \(W \in C_{n,J}(X)^-\).
\end{proof}

\begin{proposition} % C(X C(Y)) = C(X) C(Y)
\label{C(X C(Y)) = C(X) C(Y)}
\begin{enumerate}
	\item If \(J \subseteq K\), then \(C_{n,J}(X \cdot_L C_{n,K}(Y)) \leq^+ C_{n,J}(X) \cdot_L C_{n,J}(Y)\).
	\item If \(n \in K\), then \(C_{n,J}(X \cdot_L C_{n,K}(Y)) =^+ C_{n,J}(X) \cdot_L C_{n,K}(Y)\).
	\item If \(n \in L\), then 
		\begin{align*}
			C_{n,J}(X) \cdot_L C_{n,K}(Y) &\leq^- C_{n,J}(X \cdot_L C_{n,K}(Y)), \\
			C_{n,J}(X \cdot_L C_{n,K}(Y)) &\leq^- C_{n,P}(C_{n,J}(X) \cdot_L C_{n,K}(Y)).
		\end{align*}
	\item If \(n \in L\) and \((C_{n,J}(X) \cdot_L C_{n,K}(Y))^-\) is a suit, then 
	\[
	C_{n,J}(X \cdot_L C_{n,K}(Y)) =^- C_{n,J}(X) \cdot_L C_{n,K}(Y).
	\]
\end{enumerate}
\noindent Thus, if  $X$ and $Y$ are double suits and \(n \in K \cap L\), then 
\[
C_{n,J}(X \cdot_L C_{n,K}(Y)) = C_{n,J}(X) \cdot_L C_{n,K}(Y).
\]
\end{proposition}

\begin{proof}
(a)
Suppose \(V \in (C_{n,J}(X \cdot_L C_{n,K}(Y)))^+\). Then there is an \(f\colon V \toind{J} A\) such that 
\[
V(n:f) \in (X \cdot_L C_{n,K}(Y))^+ = X^+ \cap C_{n,K}(Y)^+,
\] 
which implies that there is a \(g\colon V(n:f) \toind{K} A\) such that \(V(n:f)(n:g) \in Y^+\). By \lemref{composition independent of K}, if \(J \subseteq K\) there is an \(h\colon V \toind{J} A\) such that \(V(n:f)(n:g) = V(n:h)\). Hence \(V \in C_{n,J}(X)^+ \cap C_{n,J}(Y)^+ = (C_{n,J}(X) \cdot_L C_{n,J}(Y))^+\).

(b)
Suppose \(V \in C_{n,J}(X \cdot_L C_{n,K}(Y))^+\). Then there is a function \(f\colon V \toind{J} A\) such that 
\[
V(n:f) \in (X \cdot_L C_{n,K}(Y))^+ = X^+ \cap C_{n,K}(Y)^+,
\]
which implies that there is a function \(g\colon V(n:f) \toind{K} A\) such that \(V(n:f)(n:g) \in Y^+\). If \(n \in K\), then by \lemref{composition independent of K} there is a function \(h\colon V \toind{K} A\) such that \(V(n:f)(n:g) = V(n:h)\). Hence \(V \in C_{n,J}(X)^+ \cap C_{n,K}(Y)^+ = (C_{n,J}(X) \cdot_L C_{n,K}(Y))^+\).

Conversely, suppose \(V \in (C_{n,J}(X) \cdot_L C_{n,K}(Y))^+ = C_{n,J}(X)^+ \cap C_{n,J}(Y)^+\). Then there exist \(f\colon V \toind{J} A\) and \(h\colon V \toind{J} A\) such that \(V(n:f) \in X^+\) and \(V(n:h) \in Y^+\). If \(n \in K\), then by \lemref{interpolation independent of K} there is a \(g\colon V(n:f) \toind{K} A\) such that \(V(n:f)(n:g)\linebreak[1] = V(n:h)\). It follows that \(V(n:f) \in X^+ \cap C_{n,K}(Y)^+ = (X \cdot_L C_{n,K}(Y))^+\). Hence \(V \in  (C_{n,J}(X \cdot_L C_{n,K}(Y)))^+\).

(c)
Suppose \(n \in L\) and \(W \in (C_{n,J}(X) \cdot_L C_{n,K}(Y))^-\). Then \(W = W_1 \cup_L W_2\) for some \(W_1 \in C_{n,J}(X)^-\) and \(W_2 \in C_{n,K}(Y)^-\), which implies \(W_1(n:A) \in X^-\) and \(W_2(n:A) \in Y^-\). By \lemref{V(n:A) disjoint cover}, \(W(n:A) = W_1(n:A) \cup_L W_2(n:A)\). Also, \(W_2(n:A)(n:A) = W_2(n:A) \in Y^-\), \linebreak[1] so \(W_2(n:A) \in C_{n,K}(Y)^-\). Thus, \(W(n:A) \in (X \cdot_L C_{n,K}(Y))^-\). Therefore we have \(W \in (C_{n,J}(X \cdot_L C_{n,K}(Y)))^-\).

Suppose \(W \in (C_{n,J}(X \cdot_L C_{n,K}(Y)))^-\). Then \(W(n:A) \in (X \cdot_L C_{n,K}(Y))^-\), which means \(W(n:A) = W_1 \cup_L W_2\) for some \(W_1 \in X^-\) and \(W_2 \in C_{n,K}(Y)^-\). By \lemref{V(n:A) disjoint cover}, \(W_1 = W_1(n:A)\), hence \(W_1 \in C_{n,J}(X)^-\). Thus \(W(n:A) \in (C_{n,J}(X) \cdot_L C_{n,K}(Y))^-\). Therefore \(W \in C_{n,P}(C_{n,J}(X) \cdot_L C_{n,K}(Y))^-\).

(d)
If \((C_{n,J}(X) \cdot_L C_{n,K}(Y))^-\) is a suit, then \(W \subseteq W(n:A) \in (C_{n,J}(X) \cdot_L C_{n,K}(Y))^-\) implies \(W \in (C_{n,J}(X) \cdot_L C_{n,K}(Y))^-\).
\end{proof}

To give an example where (b) fails, again let $\A$ be the equality structure with universe \(A = \set{0,1}\), and consider $\Cyls_{\mathrm{IFG}_{2}}(\A)$. Setting \(X = 1\) and \(Y = D_{01}\), observe that 
\[
D_{01} \leq C_{1,\set{0}}(D_{01})= C_{1,N}(1) \cdot_L C_{1,\set{0}}(D_{01}).
\]
In particular, \(\set{00, 11} \in (C_{1,N}(1) \cdot_L C_{1,\set{0}}(D_{01}))^+\). However, 
\[
\set{00, 11} \notin C_{1,N}(1 \cdot_L C_{1,\set{0}}(D_{01}))^+.
\]
To see why, suppose the contrary. Then there is an \(f\colon \set{00,11} \toind{N} A\) such that \(\set{00, 11}(1:f) \in (1 \cdot_L C_{1,\set{0}}(D_{01}))^+ = C_{1,\set{0}}(D_{01})^+\). Hence $f$ is a constant function. Thus, either \(\set{00, 11}(1:f) \linebreak[1] = \set{00, 10} \in C_{1,\set{0}}(D_{01})^+\) or \(\set{00, 11}(1:f) = \set{01, 11} \in C_{1,\set{0}}(D_{01})^+\). In the first case, \(\set{00, 10}(1:g) \in D_{01}^+\) for some \(g\colon \set{00, 10} \toind{\set{0}} A\). Since $g$ is a constant function, either \(\set{00, 10}(1:g) = \set{00, 10}\) or \(\set{00, 10}(1:g) = \set{01, 11}\). Note that neither set belongs to $D_{01}^+$. The second case is similar.

\subsubsection{Axiom C4}

\begin{proposition}
\begin{enumerate}
	\item \(C_{n,J}(C_{n,K}(X)) \leq C_{n,J \cap K}(X)\).
	\item If \(n \in K\), then \(C_{n,J}(C_{n,K}(X)) = C_{n,K}(X)\).
\end{enumerate}
\end{proposition}
\begin{proof}
(a)
Suppose \(V \in C_{n,J}(C_{n,K}(X))^+\). Then \(V(n:f)(n:g) \in X^+\) for some \(f\colon V \toind{J} A\) and \(g\colon V(n:f) \toind{K} A\). By \lemref{composition independent of K} there is an \(h\colon V \toind{J \cap K} A\) such that \(V(n:f)(n:g) = V(n:h)\). Hence \(V \in C_{n,J \cap K}(X)^+\). Also, \(W \in C_{n,J}(C_{n,K}(X))^-\) if and only if \(W(n:A)(n:A) \in X^-\) if and only if \(W(n:A) \in X^-\) if and only if \(W \in C_{n,J \cap K}(X)^-\).

(b) If \(n \in K\), then the $h$ from above is independent of $K$, so \(V \in C_{n,K}(X)^+\). Conversely, suppose \(V \in C_{n,K}(X)^+\). Then there is a function \(h\colon V \toind{K} A\) such that \(V(n:h) \in X^+\). Let \(f\colon V \toind{J} A\) be any function independent of $J$. By \lemref{interpolation independent of K} there is a \(g\colon V(n:f) \toind{K} A\) such that \(V(n:f)(n:g) = V(n:h)\). Hence \(V \in C_{n,J}(C_{n,K}(X))^+\).
\end{proof}

\begin{proposition} % C_mC_n = C_nC_m
\label{C_mC_n = C_nC_m}
If \(m \in K\) and \(n \in J\), where \(m \not= n\), then 
\[
C_{m,J}(C_{n,K}(X)) = C_{n,K}(C_{m,J}(X)).
\]
\end{proposition}

\begin{proof}
Suppose \(V \in C_{m,J}(C_{n,K}(X))^+\). Then there exist functions  \(f\colon V \toind{J} A\) and \(g\colon V(m:f) \toind{K} A\) such that \(V(m:f)(n:g) \in X^+\). By \lemref{functions commute} there exist \(G\colon V \toind{K} A\) and \(F\colon V(n:G) \toind{J} A\) such that \(V(m:f)(n:g) = V(n:G)(m:F)\). Hence \(V \in C_{n,K}(C_{m,J}(X))^+\). Thus \(C_{m,J}(C_{n,K}(X))^+ \subseteq C_{n,K}(C_{m,J}(X))^+\). The reverse containment follows by symmetry. Therefore 
\[
C_{m,J}(C_{n,K}(X))^+ = C_{n,K}(C_{m,J}(X))^+.
\]

Observe that  \(W \in C_{m,J}(C_{n,K}(X))^-\) if and only if \(W(m:A)(n:A) \in X^-\) if and only if \(W(n:A)(m:A) \in X^-\) if and only if \(W \in C_{n,K}(C_{m,J}(X))^-\).
\end{proof}

%%% Axiom C5
\subsubsection{Axiom C5}

\begin{proposition} %CD = 1
\label{CD = 1}
If \(j \notin J\), then \(C_{i,J}(D_{ij}) = 1\).
\end{proposition}

\begin{proof}
That \(C_{i,J}(D_{ij})^+ \subseteq 1^+\) is immediate. To show \(1^+ \subseteq C_{i,J}(D_{ij})^+\), suppose \(V \in 1^+\). Define \(f\colon V \toind{J} A\) by \(f(\vec a) = a_j\). Observe that $f$ is independent of $J$ because \(j \notin J\). Also note that \(V(i:f) = \setof{\vec a(i:a_j)}{\vec a \in V} \in D_{ij}^+\). Hence \(V \in C_{i,J}(D_{ij})^+\). Therefore \(C_{i,J}(D_{ij})^+ = 1^+\).

We know \(\emptyset \in C_{i,J}(D_{ij})^-\), so to show \(C_{i,J}(D_{ij})^- = \set{\emptyset} = 1^-\) it suffices to show that any \(W \in C_{i,J}(D_{ij})^-\) must be empty. Suppose \(W \in C_{i,J}(D_{ij})^-\). Then by definition \(W(i:A) \in D_{ij}^-\). It follows that $W = \emptyset$ because if \(\vec a \in W\), then \(\vec a(i:a_j) \in W(i:A)\), which contradicts \(W(i:A) \in D_{ij}^-\). Therefore \(C_{i,J}(D_{ij})^- = 1^-\).
\end{proof}

To show that the hypothesis \(j \notin J\) is necessary, let $\A$ be the equality structure $\A$ with universe is \(A = \set{0,1}\), and consider $\Cyls_{\mathrm{IFG}_2}(\A)$. Observe that  
\[
C_{0,N}(D_{01}) = \tuple{\powerset(\set{00, 10}) \cup \powerset(\set{01, 11}), \set{\emptyset}}.
\]

\subsubsection{Axiom C6}

\begin{proposition} % C(DD) = D
\label{C(DD) = D}
If \(i \notin J\) or \(j \notin J\), and \(i \not= k \not = j\), then 
\[
C_{k,J}(D_{ik} \cdot_\emptyset D_{kj}) = D_{ij}.
\]
\end{proposition}

\begin{proof}
Without loss of generality, suppose \(i \notin J\) and \(i \not= k \not = j\). Suppose \(V \in C_{k,J}(D_{ik} \cdot_\emptyset D_{kj})^+\). Then there is an \(f\colon V \toind{J} A\) such that \(V(k:f) \in (D_{ik} \cdot_\emptyset D_{kj})^+ = D_{ik}^+ \cap D_{kj}^+\). Thus every \(\vec b \in V(k:f)\) has the property that \(b_i = b_k = b_j\). Let \(\vec a \in V\). Then \(\vec a(k:f(\vec a)) \in V(k:f)\), so \(a_i = a_j\). Therefore \(V \in D_{ij}^+\).

Conversely, suppose \(V \in D_{ij}^+\). Then every \(\vec a \in V\) has the property that \(a_i = a_j\). Since \(i \notin J\), the projection \(\pr_i\colon V \toind{J} A\) is independent of $J$, and every \(\vec b \in V(k:\pr_i)\) has the property that \(b_i = b_k = b_j\). Hence \(V(k:\pr_i) \in D_{ik}^+ \cap D_{kj}^+ = (D_{ik} \cdot_\emptyset D_{kj})^+\). Thus \(V \in C_{k,J}(D_{ik} \cdot_\emptyset D_{kj})^+\).

Now suppose \(W \in C_{k,J}(D_{ik} \cdot_\emptyset D_{kj})^-\). Then \(W(k:A) \in (D_{ik} \cdot_\emptyset D_{kj})^-\), which means \(W(k:A) = W_1 \cup_\emptyset W_2\) for some \(W_1 \in D_{ik}^-\) and \(W_2 \in D_{kj}^-\). To show that \(W \in D_{ij}^-\), suppose to the contrary that \(W \notin D_{ij}^-\). Then there is an \(\vec a \in W\) such that \(a_i = a_j\). Let \(\vec b = \vec a(k:a_i)\). Then \(\vec b \in W(k:A)\), but \(\vec b \notin W_1\) and \(\vec b \notin W_2\), which contradicts the fact that \(W(k:A) = W_1 \cup W_2\) is a disjoint cover. Therefore \(W \in D_{ij}^-\).

Conversely, suppose \(W \in D_{ij}^-\). Then every \(\vec a \in W\) has the property that \(a_i \not= a_j\). Let \(W_1 = \setof{\vec a \in W(k:A)}{a_i \not= a_k}\) and \(W_2 = \setof{\vec a \in W(k:A)}{a_i = a_k}\). Then we have \(W(k:A) = W_1 \cup_\emptyset W_2\), \(W_1 \in D_{ik}^-\), and \(W_2 \in D_{kj}^-\). Hence \(W(k:A) \in (D_{ik} \cdot_\emptyset D_{kj})^-\). Thus \(W \in C_{k,J}(D_{ik} \cdot_\emptyset D_{kj})^-\).
\end{proof}

%%% Axiom C7
\subsubsection{Axiom C7}

Axiom C7 is a complementation axiom, so it is not surprising that it fails in IFG-algebras. Nevertheless, like the Boolean complementation axioms, axiom C7 fails in a nice way.

\begin{proposition} % C(DX)C(DnX) < Omega
\label{C(DX)C(DnX) < Omega}
If $X$ is a double suit and \(i \not= j\), then 
\[
C_{i,K}(D_{ij} \cdot_J X) \cdot_L C_{i,K}(D_{ij} \cdot_J \n{X}) \leq \Omega.
\]
\end{proposition}

\begin{proof}
Suppose \(V \in (C_{i,K}(D_{ij} \cdot_J X) \cdot_L C_{i,K}(D_{ij} \cdot_J \n{X}))^+\). Then \(V \in C_{i,K}(D_{ij} \cdot_J X)^+\) and \(V \in C_{i,K}(D_{ij} \cdot_J \n{X})^+\), which means that \(V(i:f) \in (D_{ij} \cdot_J X)^+ = D_{ij}^+ \cap X^+\) and also \(V(i:g) \in (D_{ij} \cdot_J \n{X})^+ = D_{ij}^+ \cap X^-\) for some \(f\colon V \toind{K} A\) and \(g\colon V \toind{K} A\). Hence, for every \(\vec a \in V\), \(f(\vec a) = a_j = g(\vec a)\), so \(f = g\) and \(V(i:f) = V(i:g)\). Thus \(V(i:f) \in X^+ \cap X^- = \set{\emptyset}\), which implies \(V = \emptyset\).
\end{proof}

Thus an IFG-formula of the form 
\[
\hexists{v_i}{K}(v_i = v_j \hand{J} \phi) \hand{L} \hexists{v_i}{K}(v_i = v_j \hand{J} \hneg\phi)
\]
is never true. However, it might not be false. Let $\A$ be the structure whose universe is \(A = \set{0,1}\) and in which every element is named by a constant, and consider the 2-dimensional IFG-cylindric set algebra over $\A$. Let \(X = \norm{v_0 = 0} +_N \norm{v_0 = 1}\). Then 
\[
(C_{0,\emptyset}(D_{01} \cdot_\emptyset X) \cdot_\emptyset C_{0,\emptyset}(D_{01} \cdot_\emptyset \n{X}))^- \not= 0^-.
\]
To see why, suppose to the contrary that \(\^NA \in (C_{0,\emptyset}(D_{01} \cdot_\emptyset X) \cdot_\emptyset C_{0,\emptyset}(D_{01} \cdot_\emptyset \n{X}))^-\). Then \(\^NA = W_1 \cup W_2\) for some \(W_1 \in C_{0,\emptyset}(D_{01} \cdot_\emptyset X)^-\) and \(W_2 \in C_{0,\emptyset}(D_{01} \cdot_\emptyset \n{X})^-\). Consider the valuations 00 and 11. Suppose for the sake of a contradiction that \(00 \in W_1\). Then \(00 \in W_1(0:A) \in (D_{01} \cdot_\emptyset X)^-\), which means \(W_1(0:A) = W_3 \cup W_4\) for some \(W_3 \in D_{01}^-\) and \(W_4 \in X^- = \set{\emptyset}\). Hence \(00 \in W_3 \in D_{01}^-\), which contradicts the definition of $D_{01}^-$. Thus \(00 \notin W_1\). Similarly, \(11 \notin W_1\). Therefore, \(\set{00, 11} \subseteq W_2\), which means \(\set{00, 01, 10, 11} \subseteq W_2(0:A) \in (D_{01} \cdot_J \n{X})^-\). Hence \(\set{00, 01, 10, 11} = W_5 \cup W_6\) for some \(W_5 \in D_{01}^-\) and \(W_6 \in X^+\). Without loss of generality we may assume \(W_5 = \set{01, 10}\) and \(W_6 = \set{00, 11}\), which contradicts the fact that \(\set{00, 11} \notin X^+\).

% More on cylindrifications
\subsubsection{More on cylindrifications}

\begin{lemma} % C_0...C_N-1 rooted
\label{C_0...C_N-1 rooted}
Suppose $X$ is rooted.
\begin{enumerate}
	\item \(C_{0,\emptyset} \ldots C_{{N-1},\emptyset}(X)^+ = \set{\emptyset}\) if and only if \(X^+ = \set{\emptyset}\).
	\item \(C_{0,N}\ldots C_{N-1,N}(X)^+ = \powerset(\^NA)\) if and only if there is an \(\vec a \in \^NA\) such that \(\set{\vec a} \in X^+\); otherwise \(C_{0,N}\ldots C_{N-1,N}(X)^+ = \set{\emptyset}\).
	\item \(C_{0,J_0}\ldots C_{N-1,J_{N-1}}(X)^- = \powerset(\^NA)\) if and only if \(\^NA \in X^-\).
	\item \(C_{0,J_0}\ldots C_{N-1,J_{N-1}}(X)^- = \set{\emptyset}\) if and only if \(\^NA \notin X^-\).
\end{enumerate}
\end{lemma}

\begin{proof}
(a)
If \(X^+ = \set{\emptyset}\), then for any \(V \subseteq \^NA\) we have \(V \in C_{0,\emptyset} \ldots C_{{N-1},\emptyset}(X)^+\) if and only if there exist \(f_0,\ldots, f_{N-1}\) such that \(V(0: f_0)\ldots(N-1:f_{n-1}) = \emptyset\) if and only if \(V = \emptyset\).  Hence \(C_{0,\emptyset} \ldots C_{{N-1},\emptyset}(X)^+ = \set{\emptyset}\). Conversely, if there is a nonempty \(V \in X^+\), then \linebreak[4] \(V(0:\pr_0)\ldots(N-1:\pr_{N-1}) = V \in X^+\), so \(V \in C_{0,\emptyset} \ldots C_{{N-1},\emptyset}(X)^+\).

(b)
Suppose \(\set{\vec a} \in X^+\), and let $f_n$ be the function that takes the constant value $a_n$. Then for any \(V \subseteq \^NA\), \(V(0:f_0)\ldots({N-1},f_{N-1}) = \set{\vec a} \in X^+\). Hence \(V \in C_{0,N}\ldots C_{N-1,N}(X)^+\). Thus \(C_{0,N}\ldots C_{N-1,N}(X)^+ = \powerset(\^NA)\). Conversely, suppose \(C_{0,N}\ldots C_{N-1,N}(X)^+ = \powerset(\^NA)\). Then there exist constant functions $f_n$ such that \(\^NA(0:f_0)\ldots(N-1:f_{N-1}) \in X^+\). Let $a_n$ be the constant value taken by $f_n$. Then \(\^NA(0:f_0)\ldots(N-1:f_{N-1}) = \set{\vec a}\).

If there is no \(\vec a \in \^NA\) such that \(\set{\vec a} \in X^+\), then for any nonempty \(V \subseteq \^NA\) and constant functions \(f_0,\ldots,f_{N-1}\), \(V(0:f_0)\ldots(N-1:f_{N-1}) \notin X^+\). Hence \(V \notin C_{0,N}\ldots C_{N-1,N}(X)^+\). Thus \(C_{0,N}\ldots C_{N-1,N}(X)^+ = \set{\emptyset}\).

(c) and (d) 
Let \(W \subseteq \^NA\). By definition \(W \in C_{0,J_0}\ldots C_{N-1,J_{N-1}}(X)^-\) if and only if \(W(0:A)\ldots(N-1,A) \in X^-\). If \(W = \emptyset\), then \(W(0:A)\ldots(N-1,A) = \emptyset\), otherwise \(W(0:A)\ldots(N-1,A) = \^NA\).
\end{proof}

\begin{corollary} % C_0...C_N-1 pair of suits
\label{C_0...C_N-1 pair of suits}
If $X$ is a pair of suits, then 
\[
C_{0,N}\ldots C_{N-1,N}(X) = 
\begin{cases}
1      		& \text{if \(X \not\leq \Omega\) and \(X \not\leq \mho\)}, \\
\Omega		& \text{if \(0 < X \leq \Omega\)}, \\
\mho		& \text{if \(0 < X \leq \mho\)}, \\
0	    	& \text{if \(X = 0\)}.
\end{cases}
\]
\end{corollary}

%\begin{corollary} % C_0,N...C_N-1,N double suit
%\label{C_0,N...C_N-1,N double suit}
%If $X$ is a double suit, then 
%\[
%C_{0,N}\ldots C_{N-1,N}(X) = 
%\begin{cases}
%1	    	& \text{if \(X \not\leq \Omega\)}, \\
%\Omega		& \text{if \(0 < X \leq \Omega\)}, \\
%0      		& \text{if \(X = 0\)}.
%\end{cases}
%\]
%\end{corollary}

\begin{proposition} % C_0,J...C_N-1,J double suit
\label{C_0,J...C_N-1,J double suit}
If $X$ is a double suit, then 
\[
C_{0,J_0}\ldots C_{N-1,J_{N-1}}(X) = 
\begin{cases}
1	    	& \text{if \(X \not\leq \Omega\)}, \\
\Omega		& \text{if \(0 < X \leq \Omega\)}, \\
0      		& \text{if \(X = 0\)}.
\end{cases}
\]
\end{proposition}

\begin{proof}
By \propref{C_K < C_J} and \propref{C preserves <}, 
\[
C_{0,N}\ldots C_{N-1,N}(X) \leq C_{0,J_0}\ldots C_{N-1,J_{N-1}}(X).
\]
Thus \(C_{0,J_0}\ldots C_{N-1,J_{N-1}}(X) = 1\) if \(X \not\leq \Omega\). If \(0 < X \leq \Omega\), then \(X^+ = \set{\emptyset}\) and \(\^NA \notin X^-\). Hence \(C_{0,J_0}\ldots C_{N-1,J_{N-1}}(X) = \Omega\). Finally, \(C_{0,J_0}\ldots C_{N-1,J_{N-1}}(0) = 0\) by \propref{C(0) = 0}. 
\end{proof}

In terms of IFG logic, what \propref{C_0,J...C_N-1,J double suit} says is that the semantic game associated with \(\hexists{v_0}{J_0}\ldots\hexists{v_{N-1}}{J_{N-1}}\phi\) is like the game associated with $\phi$ except that \eloise\ is allowed to choose the initial valuation. Since \eloise\ can use constant functions to specify any initial valuation she wishes, the amount of information available to her is irrelevant.

\begin{theorem} % 3 implies Omega
\label{3 implies Omega}
If $\C$ is a double-suited IFG-algebra and \(\abs{\C} > 2\), then \(\Omega \in \C\).
\end{theorem}

\begin{proof}
Suppose \(X \in \C \setminus \set{0,1}\). Then \(\^NA \notin X^+ \cup X^-\), and there exists a nonempty team \(V \in X^+ \cup X^-\). Hence \((X \cdot_N \n{X})^+ = X^+ \cap X^- = \set{\emptyset}\), and \(V \in X^+ \cup X^- = (X \cdot_N \n{X})^-\), so \(C_{0, N}\cdots C_{N-1,N}(X\cdot_N \n{X}) = \Omega\).
\end{proof}

% The trivial algebra and IFG$_0$-algebras
\subsection{The trivial algebra and IFG$_0$-algebras}

Thus far we have neglected the case when the base set \(A = \emptyset\). If \(A = \emptyset\) and \(N > 0\), then \(\^NA = \emptyset\). Hence \(0 = 1 = D_{ij} = \Omega = \mho\). Thus, by \propref{0, Omega, mho, 1} and \propref{C(Omega) = Omega}, $\set{\Omega}$ is an IFG$_N$-algebra. We will refer to $\set{\Omega}$ as the trivial IFG$_N$-algebra\index{IFG-cylindric set algebra!trivial}.

\begin{proposition} % trivial algebra
\label{trivial algebra}
The trivial algebra $\set{\Omega}$ is the only double-suited IFG$_N$-algebra with an empty base set. In fact, it is the only rooted IFG$_N$-algebra with an empty base set. It is also the only IFG$_N$-algebra with only one element.
\end{proposition}

\begin{proof}
If \(A = \emptyset\) then \(\powerset(\powerset(\^NA)) \times \powerset(\powerset(\^NA))\) is
\(
\set{\tuple{\emptyset, \emptyset}, \tuple{\emptyset, \set{\emptyset}}, \tuple{\set{\emptyset}, \emptyset}, \tuple{\set{\emptyset}, \set{\emptyset}}}.
\)
The only rooted element is $\tuple{\set{\emptyset}, \set{\emptyset}}$, which is a double suit. 

Suppose $\C$ is an IFG$_N$-algebra such that \(\abs{\C} = 1\). Then the base set $A$ must be empty because otherwise \(0 = \tuple{\set{\emptyset}, \powerset(\^NA)} \not= \tuple{\powerset(\^NA), \set{\emptyset}} = 1\).
\end{proof}

If \(N = 0\), then \(\^NA = \set{\emptyset}\). It is important to distinguish between the empty team $\emptyset$ and the team that includes only the empty valuation $\set{\emptyset}$. To avoid confusion let \(\vec a = \emptyset\) be the empty valuation. Then \(\powerset(\powerset(\^NA)) \times \powerset(\powerset(\^NA)) = \powerset(\powerset(\set{\vec a})) \times \powerset(\powerset(\set{\vec a}))\) has sixteen elements. The only suits are $\set{\emptyset}$ and $\powerset(\set{\vec a})$, so the only pairs of suits are \(0 = \tuple{\set{\emptyset}, \powerset(\set{\vec a})}\),  \(\Omega = \tuple{\set{\emptyset}, \set{\emptyset}}\), \(\mho = \tuple{\powerset(\set{\vec a}), \powerset(\set{\vec a})}\), and \(1 = \tuple{\powerset(\set{\vec a}), \set{\emptyset}}\). Thus there are three suited IFG$_0$-algebras:
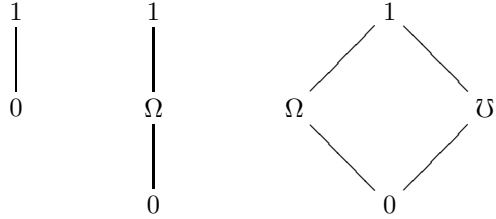
\[
\xymatrix{
	1 \ar@{-}[d] \\
	0
}
\qquad\qquad
\xymatrix{
	1 \ar@{-}[d] \\
	\Omega \ar@{-}[d] \\
	0
}
\qquad\qquad
\xymatrix{
						& 1 \ar@{-}[dl] \ar@{-}[dr] \\
\Omega \ar@{-}[dr] 	&									& \mho \ar@{-}[dl]\\
						& 0
}
\]

%\input{3-Finite_Axiomatizability}
%LaTeX by Allen L. Mann

% Finite Axiomatizability
\section{Finite Axiomatizability}

We would like to know whether the equational theory of IFG$_N$-algebras is finitely axiomatizable. Unfortunately determining the answer is beyond our current abilities, even in the one-dimensional case. To make the problem more tractable, we consider certain reducts of IFG$_1$-algebras. First we will examine reducts of IFG$_1$-algebras to the signature 
\(
\tuple{0,1,\n{},+_{\set{0}},\cdot_{\set{0}}}.
\)
Then will consider reducts to the signature 
\[
\tuple{0,1,\n{},+_{\set{0}},\cdot_{\set{0}}, C_{0,\set{0}}}.
\]

%% The De Morgan reduct
\subsection{The De Morgan reduct}

\begin{definition} % De Morgan reduct
The reduct of an IFG$_N$-cylindric set algebra $\C$ to the signature \(\tuple{0,1,\n{},+_N,\cdot_N}\) is called the \emph{De Morgan reduct}\index{De Morgan reduct|mainidx} of $\C$.
\end{definition}

\begin{definition} % fixed point
An element $x$ of a De Morgan algebra is a \emph{fixed point}\index{fixed point|mainidx} (or \emph{center}\index{center|mainidx}) if \({\hneg x} = x\). A \emph{centered De Morgan algebra}\index{De Morgan algebra!centered|mainidx} is a De Morgan algebra with a center.
\end{definition}

\begin{proposition} % centered Kleene algebra
\label{centered Kleene algebra}
The center of a centered Kleene algebra is unique.
\end{proposition}

\begin{proof}
Suppose $a$ and $b$ are both fixed points of a centered Kleene algebra. Then \(a = a \meet {\hneg a} \leq b \join {\hneg b} = b\), and \(b = b \meet {\hneg b} \leq a \join {\hneg a} = a\).
\end{proof}

Let $\mathbf{\underline{B}}$\index{$\mathbf{\underline{B}}$} be the two-element De Morgan algebra, let \(\mathbf{\underline{K}}\)\index{$\mathbf{\underline{K}}$} be the three-element De Morgan algebra with fixed point $a$, and let \(\mathbf{\underline{M}}\)\index{$\mathbf{\underline{M}}$} be the four-element De Morgan algebra with two fixed points $a$ and $b$ such that \(a \meet b = 0\) and \(a \join b = 1\):
\[
\xymatrix{
1 \ar@{-}[d]	&&	1 \ar@{-}[d]	&&					& 1 \ar@{-}[dl]\ar@{-}[dr]						\\
0				&& a \ar@{-}[d]	&& a \ar@{-}[dr] 	&								& b \ar@{-}[dl] 	\\
				&& 0				&&					& 0												\\
\textbf{\underline{B}}				&& \textbf{\underline{K}}				&&					& \textbf{\underline{M}}
}
\]

\begin{theorem}[Kalman\index{Kalman, J.~A.} \cite{Kalman:1958}] % subdirectly irreducible De Morgan algebras
\label{subdirectly irreducible De Morgan algebras}
The subdirectly irreducible De Morgan algebras are exactly \linebreak[3] $\mathbf{\underline{B}}$, $\mathbf{\underline{K}}$, and $\mathbf{\underline{M}}$.
\end{theorem}

The following corollary is an immediate consequence of \thmref{subdirectly irreducible De Morgan algebras} and Birkhoff's subdirect representation theorem (II\S8.6 in \cite{Burris:1981}).

\begin{corollary} % subdirect product of B, K, M
\label{subdirect product of B, K, M}
Every De Morgan algebra is isomorphic to a subdirect product of copies of $\mathbf{\underline{B}}$, $\mathbf{\underline{K}}$, and $\mathbf{\underline{M}}$.
\end{corollary}

\begin{corollary} % subdirect product of B, K
\label{subdirect product of B, K}
\begin{enumerate}
	\item Every Kleene algebra is isomorphic to a subdirect product of copies of $\mathbf{\underline{B}}$ and $\mathbf{\underline{K}}$.
	\item Every centered Kleene algebra is isomorphic to a subdirect product of copies of $\mathbf{\underline{K}}$.
	\item Every De Morgan algebra that has a subalgebra isomorphic to $\mathbf{\underline{M}}$ is isomorphic to a subdirect product of copies of $\mathbf{\underline{M}}$.
\end{enumerate}
\end{corollary}

\begin{proof}
(a)
Let $\A$ be a Kleene algebra. If $\mathbf{\underline{M}}$ were a factor in a subdirect representation of $\A$, then $\mathbf{\underline{M}}$ would be a homomorphic image of $\A$, and hence a Kleene algebra. 

(b)
Let $\A$ be a Kleene algebra, and let $\B$ be a subdirect representation of $\A$. If $\mathbf{\underline{B}}$ is a factor of $\B$ (say the first), then every tuple in $\B$ has the form $\tuple{0, \ldots}$ or $\tuple{1, \ldots}$, so no tuple in $\B$ is a fixed point.

(c)
Let $\A$ be a De Morgan algebra, and let $\B$ be a subdirect representation of $\A$. If $\mathbf{\underline{B}}$ is a factor of $\B$, then $\B$ would have no fixed points at all, and if $\mathbf{\underline{K}}$ is a factor of $\B$ (say the first), then every fixed point in $\B$ would have the form  $\tuple{a, \ldots}$, so no two fixed points could meet to $\tuple{0, 0, \ldots}$ or join to $\tuple{1, 1, \ldots}$.
\end{proof}

Let $\mathcal{B}$\index{$\mathcal{B}$ \quad variety of Boolean algebras} denote the variety of Boolean algebras, $\mathcal{K}$\index{$\mathcal{K}$ \quad variety of Kleene algebras} denote the variety of Kleene algebras, and $\mathcal{M}$\index{$\mathcal{M}$ \quad variety of De Morgan algebras} denote the variety of De Morgan algebras.

\begin{theorem}[Kalman\index{Kalman, J.~A.} \cite{Kalman:1958}] % subvarieties of De Morgan algebras
\label{subvarieties of De Morgan algebras}
The nontrivial subvarieties of $\mathcal{M}$ are 
\(
\mathcal{B} \subset \mathcal{K} \subset \mathcal{M}.
\)
\end{theorem}

\begin{theorem} % De Morgan variety generated
\label{De Morgan variety generated}
The variety generated by the De Morgan reducts of all suited IFG$_N$-algebras is $\mathcal M$.
The variety generated by the De Morgan reducts of all double-suited IFG$_N$-algebras is $\mathcal K$.
\end{theorem}

\begin{proof}
The De Morgan reduct of every suited IFG$_N$-algebra is a De Morgan algebra. Hence the variety generated by the De Morgan reducts of all suited IFG$_N$-algebras is contained in $\mathcal M$. On the other hand, \(\set{0, 1}\), \(\set{0, \Omega, 1}\), and \(\set{0, \Omega, \mho, 1}\) are all suited IFG$_N$-algebras, and their De Morgan reducts are isomorphic to $\mathbf{\underline{B}}$, $\mathbf{\underline{K}}$, and $\mathbf{\underline{M}}$, respectively. Thus the De Morgan  reducts of suited IFG$_N$-algebras generate all the subdirectly irreducible De Morgan algebras, which in turn generate $\mathcal M$.

The De Morgan reduct of every double-suited IFG$_N$-algebra is a Kleene algebra. Hence the variety generated by the De Morgan reducts of all double-suited IFG$_N$-algebras is contained in $\mathcal K$. On the other hand, \(\set{0, 1} \iso \mathbf{\underline{B}}\) and  \(\set{0, \Omega, 1} \iso \mathbf{\underline{K}}\) are both double-suited IFG$_N$-algebras. Thus the De Morgan  reducts of double-suited IFG$_N$-algebras generate all the subdirectly irreducible Kleene algebras, which in turn generate $\mathcal K$.
\end{proof}

It follows from \thmref{De Morgan variety generated} that the equational theory of De Morgan reducts of suited IFG$_N$-algebras is axiomatized by the axioms of De Morgan algebra, and the equational theory of De Morgan reducts of double-suited IFG$_N$-algebras is axiomatized by the axioms of Kleene algebra.

%% \input{topological_duality.tex}

%% The monadic De Morgan reduct
\subsection{Monadic De Morgan algebras}

A monadic Boolean algebra\index{monadic Boolean algebra} is a Boolean algebra with an additional unary operation called a quantifier. Monadic Boolean algebras were first studied by Halmos\index{Halmos, P.~R.} \cite{Halmos:1962} following the investigations of McKinsey\index{McKinsey, J.~C.~C.} and Tarski\index{Tarski, A.} \cite{McKinsey:1944} into the algebraic properties of closure operators. Cignoli\index{Cignoli, R.} \cite{Cignoli:1991} generalized the notion of a monadic Boolean algebra by adding quantifiers to bounded distributive lattices. Cignoli \index{Cignoli, R.} calls a bounded distributive lattice with a quantifier a $Q$-distributive lattice. Petrovich\index{Petrovich, A.} \cite{Petrovich:1996, Petrovich:1999} extended the results of Cignoli\index{Cignoli, R.} to monadic De Morgan algebras\index{monadic De Morgan algebra}. Our presentation follows \cite{Petrovich:1999}.

\begin{definition} % Q-distributive lattice
An \emph{existential quantifier}\index{quantifier!existential|mainidx} on a bounded distributive lattice is a unary operation $\nabla$\index{$\nabla x$ \quad existential quantifier|mainidx} such that:
\begin{enumerate}
	\item[(Q1)] \(\nabla 0 = 0\),
	\item[(Q2)] \(x \leq \nabla x \),
	\item[(Q3)] \(\nabla(x \join y) = \nabla x \join \nabla y\),
	\item[(Q4)] \(\nabla(x \meet \nabla y) = \nabla x \meet \nabla y\).
\end{enumerate}
A \emph{universal quantifier}\index{quantifier!universal|mainidx} is a unary operation $\Delta$\index{$\Delta x$ \quad universal quantifier|mainidx} such that:
\begin{enumerate}
	\item[(Q1$'$)] \(\Delta 1 = 1\),
	\item[(Q2$'$)] \(\Delta x \leq x \),
	\item[(Q3$'$)] \(\Delta(x \meet y) = \Delta x \meet \Delta y\),
	\item[(Q4$'$)] \(\Delta(x \join \Delta y) = \Delta x \join \Delta y\).
\end{enumerate}
The two kinds of quantifiers are dual to each other in the usual way. A \emph{quantifier}\index{quantifier|mainidx} without modification is assumed to be existential. A \emph{Q-distributive lattice}\index{$Q$-distributive lattice|mainidx} \(\tuple{L, \nabla}\) is a bounded distributive lattice $L$ equipped with a quantifier $\nabla$.
\end{definition}

\begin{definition} % monadic De Morgan algebra
An \emph{quantifier}\index{quantifier|mainidx} on a De Morgan algebra must satisfy the additional condition:
\begin{enumerate} \setcounter{enumi}{4}
	\item[(Q5)] \(\nabla({\hneg\nabla x}) = {\hneg\nabla x}\).
\end{enumerate}
A De Morgan algebra equipped with a quantifier is called a \emph{monadic De Morgan algebra}\index{monadic De Morgan algebra|mainidx}.
\end{definition}

\begin{proposition} % Q-distributive lattice properties
\label{Q-distributive lattice properties}
In any $Q$-distributive lattice,
\begin{enumerate}
	\item \(\nabla 1 = 1\),
	\item \(\nabla\nabla x = \nabla x\).
\end{enumerate}
\end{proposition}

\begin{proof}
(a)
\(1 \leq \nabla 1\) by axiom Q2, and \(\nabla 1 \leq 1\) because for all $x$, \(x \leq 1\).

(b)
\(\nabla\nabla x = \nabla(1 \meet \nabla x) = \nabla 1 \meet \nabla x = 1 \meet \nabla x = \nabla x\).
\end{proof}

\begin{lemma} % quantifier iff subalgebra
\label{quantifier iff subalgebra}
If $\A$ is a De Morgan algebra, and $\nabla$ is a quantifier on its underlying distributive lattice, then $\nabla$ is a quantifier on $\A$ if and only if the range of $\nabla$ is a subalgebra of $\A$.
\end{lemma}

\begin{proof}
Suppose $\nabla$ is a quantifier on $\A$. Then \(0 = \nabla 0\), \(1 = \nabla 1\), and \({\hneg\nabla x} = \nabla({\hneg\nabla x})\).
Also, \(\nabla x \join \nabla y = \nabla(x \join y)\), and \(\nabla x \meet \nabla y = \nabla(x \meet \nabla y)\). Thus the range of $\nabla$ is a subalgebra of $\A$. Conversely, suppose the range of $\nabla$ is a subalgebra of $\A$. Then \({\hneg\nabla x}\) is in the range of $\nabla$, so \(\nabla({\hneg\nabla x}) = {\hneg\nabla x}\).
\end{proof}

\begin{example} % type 0 quantifier
Let $\A$ be a De Morgan algebra, and let $\nabla$ be defined by 
\[
\nabla x = 
\begin{cases}
1     & \text{if \(x > 0\)}, \\
0     & \text{if \(x = 0\)}.
\end{cases}
\]
Then $\nabla$ is a quantifier on $\A$. Such a quantifier will be called a \emph{quantifier of type 0}\index{quantifier!type 0}. If $\A$ is any De Morgan algebra, then $\tuple{\A, \nabla_0}$\index{$\tuple{\A, \nabla_0}$ \quad monadic De Morgan algebra with the type 0 quantifier} will denote the monadic De Morgan algebra $\A$ with the type 0 quantifier. In particular, \(\tuple{\mathbf{\underline{B}}, \nabla_0}\), \(\tuple{\mathbf{\underline{K}}, \nabla_0}\), and \(\tuple{\mathbf{\underline{M}}, \nabla_0}\) are all monadic De Morgan algebras.
\end{example}

\begin{lemma} % type 1 lemma
\label{type 1 lemma}
If $\A$ is a centered Kleene algebra with fixed point $c$, then 
\[
\setof{x \in \A}{x \meet c = 0} = \set{0}.
\]
\end{lemma}

\begin{proof}
If $\A$ is a centered Kleene algebra with fixed point $c$, then $\A$ is isomorphic to a subdirect product $\B$ of copies of $\mathbf{\underline{K}}$. If \(x \not= 0\) in $\A$, then the corresponding tuple \(\vec y \in \B\) has at least one nonzero coordinate $y_i$. Let \(\vec a = \tuple{a, a, \ldots}\) be the fixed point of $\B$. Then \(y_i \in \set{a, 1}\), so \(y_i \meet a = a\) in $\mathbf{\underline{K}}$. Hence \(\vec y \meet \vec a \not= \vec 0\) in $\B$. Thus \(x \meet c \not= 0\) in $\A$.
\end{proof}

\begin{example} % type 1 quantifier
Let $\A$ be a centered De Morgan algebra, and let $a$ be a fixed point of $\A$. Define $\nabla$ by 
\[
\nabla x = 
\begin{cases}
1	& \text{if \(x \not\leq a\)}, \\
a	& \text{if \(0 < x \leq a\)}, \\
0	& \text{if \(x = 0\)}.
\end{cases}
\]
Then $\nabla$ is a quantifier on $\A$ if and only if \(\setof{x \in \A}{x \meet a = 0} = \set{0}\). Such a quantifier will be called a \emph{quantifier of type 1}\index{quantifier!type 1}. If $\A$ is a centered Kleene algebra with fixed point $a$, then $a$ is unique and \(\setof{x \in \A}{x \meet a = 0} = \set{0}\). Hence $\nabla$ is a quantifier on $\A$. If $\A$ is a centered Kleene algebra, let \(\tuple{\A, \nabla_1}\)\index{$\tuple{\A, \nabla_1}$ \quad centered Kleene algebra with the type 1 quantifer} denote the monadic Kleene algebra $\A$ with the type 1 quantifier. In particular, \(\tuple{\mathbf{\underline{K}}, \nabla_1}\) is a monadic Kleene algebra.
\end{example}

\begin{lemma} % type 2 lemma
\label{type 2 lemma}
Let $L$ be a bounded distributive lattice. If \(a \meet b = 0\) and \(a \join b = 1\), then
\[
\setof{x \in L}{x \meet a = 0} = [0,b] \quad\text{ and }\quad \setof{x \in L}{x \meet b = 0} = [0,a].
\]
\end{lemma}

\begin{proof}
If \(0 \leq x \leq b\), then \(0 \leq a \meet x \leq a \meet b = 0\). Hence \(x \meet a = 0\). Conversely, if \(x \meet a = 0\), then 
\[
x \join b = (x \join b) \meet (a \join b) = (x \meet a) \join b = b.
\]
Thus \(0 \leq x \leq b\).
\end{proof}

\begin{example}% type 2 quantifier
\index{quantifier!type 2(}
Let $\A$ be a centered De Morgan algebra with two fixed points $a$ and $b$ such that \(a \meet b = 0\) and \(a \join b = 1\), i.e., $\A$ contains a subalgebra isomorphic to $\mathbf{\underline{M}}$.
Define $\nabla$ by 
\[
\nabla x = 
\begin{cases}
1	& \text{if \(x \not\leq a\) and \(x \not\leq b\)}, \\
a	& \text{if \(0 < x \leq a\)}, \\
b	& \text{if \(0 < x \leq b\)}, \\
0	& \text{if \(x = 0\)}.
\end{cases}
\]
It is easy to check that $\nabla$ satisfies axioms Q1--Q3. To check Q4, suppose \(x,y \in \A\). If \(x = 0\) or \(y = 0\), then \(\nabla(x \meet \nabla y) = 0 = \nabla x \meet \nabla y\). If \(0 < x,y \leq a\), then \(\nabla(x \meet \nabla y) = a = \nabla x \meet \nabla y\). Similarly, if \(0 < x,y \leq b\), then \(\nabla(x \meet \nabla y) = b = \nabla x \meet \nabla y\). If \(x \not\leq a,b\) and \(0 < y \leq a\), then \(x \meet a \not= 0\) by \lemref{type 2 lemma}, so 
\[
\nabla(x \meet \nabla y) = \nabla(x \meet a) = a = 1 \meet a = \nabla x \meet \nabla y.
\]
Similarly, if \(x \not\leq a,b\) and \(0 < y \leq b\), then \(x \meet b \not= 0\), so 
\[
\nabla(x \meet \nabla y) = \nabla(x \meet b) = b = 1 \meet b = \nabla x \meet \nabla y.
\]
Finally, if \(y \not\leq a,b\), then \(\nabla(x \meet \nabla y) = \nabla(x \meet 1) = \nabla x = \nabla x \meet 1 = \nabla x \meet \nabla y\).

Thus $\nabla$ is a quantifier on the underlying bounded distributive lattice of $\A$, and since the range of $\nabla$ is a subalgebra of $\A$, \lemref{quantifier iff subalgebra} tells us that $\nabla$ is a quantifier on $\A$. Such a quantifier will be called a \emph{quantifier of type 2}\index{quantifier!type 2)}. If $\A$ is a De Morgan algebra with two specified fixed points $a$ and $b$, then \(\tuple{\A, \nabla_2}\)\index{$\tuple{\A, \nabla_2}$ \quad monadic De Morgan algebra with a type 2 quantifier} will denote the monadic De Morgan algebra $\A$ with the associated type 2 quantifier. In particular, \(\tuple{\mathbf{\underline{M}}, \nabla_2}\) is a monadic De Morgan algebra.
\end{example}

\begin{proposition}[Petrovich \cite{Petrovich:1999}] % subdirectly irreducible -> subdirectly irreducible
\label{subdirectly irreducible -> subdirectly irreducible}
Let \(\tuple{\A, \nabla}\) be a subdirectly irreducible monadic De Morgan algebra. Then the range of $\nabla$ is a subdirectly irreducible De Morgan algebra.
\end{proposition}

\begin{corollary}[Petrovich \cite{Petrovich:1999}] % subdirectly irreducible -> type 012
\label{subdirectly irreducible -> type 012}
Let \(\tuple{\A, \nabla}\) be a subdirectly irreducible monadic De Morgan algebra. Then $\nabla$ is a quantfier of type 0, a quantifer of type 1, or a quantifier of type 2.
\end{corollary}

The varieties of monadic Boolean, Kleene, and De Morgan algebras will be denoted by $\mathcal B^\nabla$\index{$\mathcal B^\nabla$ \quad variety of monadic Boolean algebras}, $\mathcal K^\nabla$\index{$\mathcal K^\nabla$ \quad variety of monadic Kleene algebras}, and $\mathcal M^\nabla$\index{$\mathcal M^\nabla$ \quad variety of monadic De Morgan algebras}, respectively. We can form subvarieties of $\mathcal M^\nabla$ by imposing additional conditions on $\nabla$. For example, let $\mathcal M^{\nabla_{01}}$\index{$\mathcal M^{\nabla_{01}}$ \quad variety of monadic De Morgan algebras satisfying \(\nabla x \meet {\hneg \nabla x} \leq \nabla y \join {\hneg \nabla y}\)} denote the variety of monadic De Morgan algebras that satisfy the equation equivalent to 
\[
\nabla x \meet {\hneg \nabla x} \leq \nabla y \join {\hneg \nabla y},\]
i.e., the range of $\nabla$ is a Kleene algebra, and let $\mathcal M^{\nabla_0}$\index{$\mathcal M^{\nabla_0}$ \quad variety of monadic De Morgan algebras satisfying \(\nabla x \meet {\hneg\nabla x} = 0\)} denote the variety of monadic De Morgan algebras that satisfy the equation 
\[
\nabla x \meet {\hneg\nabla x} = 0,
\]
i.e., the range of $\nabla$ is a Boolean algebra. It follows that if \(\tuple{\A, \nabla}\) is a subdirectly irreducible monadic De Morgan algebra in $\mathcal M^{\nabla_{01}}$, then $\nabla$ must be a quantifier of type 0 or 1, and if it is subdirectly irreducible in $\mathcal M^{\nabla_0}$, then $\nabla$ is a quantifier of type 0. Similarly, let $\mathcal K^{\nabla_0}$\index{$\mathcal K^{\nabla_0}$ \quad variety of monadic Kleene algebras satisfying \(\nabla x \meet {\hneg\nabla x} = 0\)} denote the variety of monadic Kleene algebras that satisfy \(\nabla x \meet {\hneg\nabla x} = 0\). In particular, we have that \(\tuple{\mathbf{\underline{M}}, \nabla_2}\) belongs to $\mathcal M^\nabla$ but not $\mathcal M^{\nabla_{01}}$, that \(\tuple{\mathbf{\underline{M}}, \nabla_0}\) belongs to $\mathcal M^{\nabla_{0}}$ but not $\mathcal K^\nabla$, and that \(\tuple{\mathbf{\underline{K}}, \nabla_1}\) belongs to $\mathcal K^{\nabla}$ but not $\mathcal M^{\nabla_{0}}$. Also, \(\mathcal M^{\nabla_0} \cap \mathcal K^\nabla = \mathcal K^{\nabla_0}\), and \(\tuple{\mathbf{\underline{K}}, \nabla_0}\) belongs to $\mathcal K^{\nabla_0}\) but not $\mathcal B^\nabla$. 

Thus we have the following sublattice of the lattice of subvarieties of $\mathcal M^\nabla$:
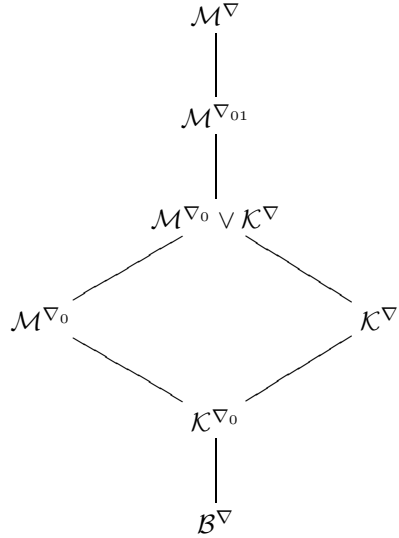
\[
\xymatrix{
										& {\mathcal M^\nabla} \ar@{-}[d]						&	\\
										& {\mathcal M^{\nabla_{01}}} \ar@{-}[d] 				&	\\
										& {\mathcal M^{\nabla_{0}} \join \mathcal K^\nabla} \ar@{-}[dl] \ar@{-}[dr]	&	\\
{\mathcal M^{\nabla_0}} \ar@{-}[dr] 	& 														& {\mathcal K^{\nabla}} \ar@{-}[dl]	\\
										& {\mathcal K^{\nabla_0}} \ar@{-}[d]					&	\\
										& {\mathcal B^{\nabla}}								&			
}
\]

To show that \(\mathcal M^{\nabla_{0}} \join \mathcal K^\nabla\) is properly contained in $\mathcal M^{\nabla_{01}}$ will require some results from universal algebra.
%which appear as II\S12.6 and IV\S6.8 in \cite{Burris:1981}.

\begin{definition} % congruence-distributive
A class of algebras is \emph{congruence-distributive}\index{congruence-distributive|mainidx} if the congruence lattice of every algebra in the class is distributive. 
\end{definition}

\begin{theorem}[J\'onsson\index{J\'onsson, B.} \cite{Jonsson:1967}] % Jonsson terms
\label{Jonsson terms}
A variety $\mathcal V$ is congruence-distributive if and only if for some integer $n \geq 2$ there exist terms \(p_0(x,y,z), \ldots, p_n(x,y,z)\) such that $\mathcal V$ satisfies
\begin{alignat*}{2}
p_0(x,y,z) &= x, \\
p_n(x,y,z) &= z, \\
p_i(x,y,x) &= x 				& \qquad & 0 \leq i \leq n, \\
p_i(x,x,y) &= p_{i+1}(x,x,y)	& \qquad & \text{for $i$ even,} \\
p_i(x,y,y) &= p_{i+1}(x,y,y)	& \qquad & \text{for $i$ odd.}
\end{alignat*}
\end{theorem}
\index{J\'onsson terms}

\begin{example}
The variety of lattices is congruence-distributive because it has J\'onsson terms\index{J\'onsson terms}
\begin{align*}
p_0(x,y,z) &= x, \\
p_1(x,y,z) &= (x \meet y) \join (y \meet z) \join (z \meet x), \\
p_2(x,y,z) &= z.
\end{align*}
It follows that any algebra that has an underlying lattice structure is congruence-distributive. In particular, the variety of De Morgan algebras and the variety of monadic De Morgan algebras are both congruence-distributive.
\end{example}

%\begin{definition} % Pu(K)
%If $K$ is a class of algebras, let $\mathsf H(K)$, $\mathsf S(K)$, $\mathsf P(K)$, and $\Pu(K)$ denote the classes consisting, repectively, of all homomorphic images, subalgebras, products, and ultraproducts of members of $K$.
%\end{definition}

%\begin{theorem}[J\'onsson \cite{Jonsson:1967}] % HSPu
%\label{HSPu}
%Let $\HSP(K)$ be a congruence-distributive variety. If $\A$ is a subdirectly irreducible algebra in $\HSP(K)$, then \(\A \in \HSPu(K)\).
%\end{theorem}

\begin{lemma}[J\'onsson\index{J\'onsson, B.} \cite{Jonsson:1967}] % V1 join V2 -> V1 cup V2
\label{V1 join V2 -> V1 cup V2}
If $\mathcal V_1$ and $\mathcal V_2$ are varieties such that \(\mathcal V_1 \join \mathcal V_2\) is congruence-distributive, then every member of \(\mathcal V_1 \join \mathcal V_2\) is isomorphic to a subdirect product of a member of $\mathcal V_1$ and a member of $\mathcal V_2$ and, in particular, every subdirectly irreducible member of $\mathcal V_1 \join \mathcal V_2$ belongs either to $\mathcal V_1$ or $\mathcal V_2$.
\end{lemma}

\begin{corollary}[J\'onsson\index{J\'onsson, B.} \cite{Jonsson:1967}] % distributive varieties
\label{distributive varieties}
If $\mathcal V$ is a congruence-distributive variety, then then lattice of all subvarieties of $\mathcal V$ is distributive.
\end{corollary}

\begin{lemma} % simple M-nabla-01 not M-nabla-0 nor K-nabla
\label{simple M-nabla-01 not M-nabla-0 nor K-nabla}
There is a simple monadic De Morgan algebra \(\tuple{\A,\nabla}\) with a quantifier of type 1 such that $\A$ is not a Kleene algebra.
\end{lemma}

\begin{proof}
Let $\A$ be the following subalgebra of \(\mathbf{\underline{K}} \times \mathbf{\underline{M}}\), where $a$ is the fixed point of $\mathbf{\underline{K}}$, and $b$, $c$ are the fixed points of $\mathbf{\underline{M}}$:
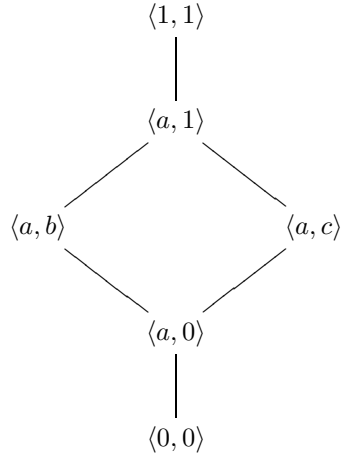
\[
\xymatrix{
										& \tuple{1,1} \ar@{-}[d]						&	\\
										& \tuple{a,1} \ar@{-}[dl] \ar@{-}[dr]	&	\\
\tuple{a,b} \ar@{-}[dr] 	& 														& \tuple{a,c} \ar@{-}[dl]	\\
										& \tuple{a,0} \ar@{-}[d]					&	\\
										& \tuple{0,0}								&			
}
\]
Consider the fixed point $\tuple{a,b}$. The only \(x \in \A\) that satisfies \(x \meet \tuple{a,b} = \tuple{0,0}\) is $\tuple{0,0}$. Therefore \(\set{\tuple{0,0}, \tuple{a,b}, \tuple{1,1}}\) is the range of a type 1 quantifier $\nabla$ on $\A$.

To show that $\tuple{\A, \nabla}$ is simple, suppose \(\tuple{0,0} \cong \tuple{a,0}\). Then \(\tuple{0,0} = \nabla\tuple{0,0} \cong \nabla\tuple{a,0} = \tuple{a,b}\), which implies \(\tuple{1,1} = {\hneg\tuple{0,0}} \cong {\hneg\tuple{a,b}} = \tuple{a,b}\). Hence \(\tuple{0,0} \cong \tuple{1,1}\), and $\cong$ is the total congruence. Now suppose \(\tuple{a,0} \cong \tuple{a,b}\). Then \(\tuple{a,1} = {\hneg\tuple{a,0}} \cong {\hneg\tuple{a,b}} = \tuple{a,b}\), which implies \(\tuple{1,1} = \nabla\tuple{a,1} \cong \nabla\tuple{a,b} = \tuple{a,b}\). Hence \(\tuple{0,0} \cong \tuple{1,1}\), and $\cong$ is the total congruence. Finally, suppose \(\tuple{a,b} \cong \tuple{a,c}\). Then \(\tuple{a,b} = \nabla\tuple{a,b} \cong \nabla\tuple{a,c} = \tuple{1,1}\). Hence $\cong$ is the total congruence.
\end{proof}

\begin{proposition} % M-nabla-0 join K-nabla < M-nabla-0
\label{M-nabla-0 join K-nabla < M-nabla-0}
\(\mathcal M^{\nabla_0} \join \mathcal K^\nabla \subset \mathcal M^{\nabla_{01}}\).
\end{proposition}

\begin{proof}
Let \(\tuple{\A, \nabla}\) be a simple monadic De Morgan algebra \(\tuple{\A,\nabla}\) with a quantifier of type 1 such that $\A$ is not a Kleene algebra. Then \(\tuple{\A, \nabla} \notin \mathcal M^{\nabla_0} \cup \mathcal K^\nabla\), and by \lemref{V1 join V2 -> V1 cup V2}, \(\tuple{\A, \nabla} \notin \mathcal M^{\nabla_0} \join \mathcal K^\nabla\).
\end{proof}

We have shown that we can separate the monadic De Morgan algebras with quantifiers of type 0 from those with quantifiers of type 1 or type 2 by means of an equation asserting that the range of the quantifier is a Boolean algebra. We can separate the monadic De Morgan algebras with quantifiers of type 0 or type 1 from those with quantifiers of type 2 by an equation asserting that the range of the quantifier is a Kleene algebra. However what we cannot yet do is separate those De Morgan algebras with type 1 quantifiers from those with type 0 or type 2 quantifiers. In particular, is there a variety that includes all monadic Kleene algebras with a quantifier of type 1 but not those with a quantifier of type 0? For starters, we should look for an equation that is true in \(\tuple{\mathbf{\underline K}, \nabla_1}\) but false in \(\tuple{\mathbf{\underline K}, \nabla_0}\). Let $a$ be the fixed point of $\mathbf{\underline K}$. Notice that for all \(x \in \mathbf{\underline K}\), 
\[
\nabla_1(x \meet {\hneg x}) \leq {\hneg\nabla_1(x \meet {\hneg x})},
\]
while 
\[
\nabla_0(a \meet {\hneg a}) \not\leq {\hneg\nabla_0(a \meet {\hneg a})}.
\]

\begin{proposition} % fixed by nabla
\label{fixed by nabla}
Let $\tuple{\A, \nabla}$ be a monadic De Morgan algebra. If $\tuple{\A, \nabla}$ satisfies 
\[
\nabla(x \meet {\hneg x}) \leq {\hneg\nabla(x \meet {\hneg x})},
\]
then every fixed point of $\A$ is fixed by $\nabla$.
\end{proposition}

\begin{proof}
First, suppose $\tuple{\A, \nabla}$ is subdirectly irreducible and satisfies \(\nabla(x \meet {\hneg x}) \leq {\hneg\nabla(x \meet {\hneg x})}\). Then $\nabla$ is a quantifier of type 0, type 1, or type 2. If $c$ is a fixed point that is not fixed by $\nabla$, then \(\nabla c = 1\) because any two distinct fixed points are incomparable. But then \(\nabla(c \meet {\hneg c}) = \nabla c = 1\) and \({\hneg\nabla(c \meet {\hneg c})} = 0\), contrary to hypothesis.

Now let $\tuple{\A, \nabla}$ be any monadic De Morgan algebra that satisfies \(\nabla(x \meet {\hneg x}) \leq {\hneg\nabla(x \meet {\hneg x})}\). Every fixed point $c$ in $\A$ corresponds to a tuple of fixed points in the subdirect representation of $\A$ in which every coordinate is fixed by $\nabla$. Hence $c$ is fixed by $\nabla$. 
\end{proof}

\begin{corollary} % subdirectly irreducible M-nabla-2
\label{subdirectly irreducible M-nabla-2}
Let \(\tuple{\A, \nabla}\) be a subdirectly irreducible monadic De Morgan algebra that satisfies \(\nabla(x \meet {\hneg x}) \leq {\hneg\nabla(x \meet {\hneg x})}\). Then either
\begin{enumerate}
	\item $\A$ is a centered De Morgan algebra with exactly two fixed points, and $\nabla$ is a quantifier of type 2,
	\item $\A$ is a centered Kleene algebra, and $\nabla$ is a quantifier of type 1, or
	\item $\A$ is a Boolean algebra, and $\nabla$ is a quantifier of type 0.
\end{enumerate}
\end{corollary}

\begin{proof}
Since $\nabla$ is a quantifier of type 0, type 1, or type 2, and every fixed point in $\A$ is fixed by $\nabla$, $\A$ can have at most two fixed points. Furthermore, if $\A$ has two distinct fixed points, then $\nabla$ must be a quantifier of type 2. 

If $\A$ has a unique fixed point $c$, then $\nabla$ must be a quantifier of type 1. Suppose for the sake of a contradiction that there exists an \(x \in \A\) such that \(x \meet {\hneg x} \not\leq c\). Then \(\nabla(x \meet {\hneg x}) = 1\) and \({\hneg \nabla(x \meet {\hneg x})} = 0\), contrary to hypothesis. Therefore, for every \(x,y \in \A\) we have \(x \meet {\hneg x} \leq c \leq y \join {\hneg y}\). Hence $\A$ is a Kleene algebra.

If $\A$ is not centered, then $\nabla$ must be a quantifier of type 0. Suppose for the sake of a contradiction that there exists an \(x \in \A\) such that \(x \meet {\hneg x} > 0\). Then \(\nabla(x \meet {\hneg x}) = 1\) and \({\hneg\nabla(x \meet {\hneg x})} = 0\), contrary to hypothesis. Therefore $\A$ satisfies the equation \(x \meet {\hneg x} = 0\). Hence $\A$ is a Boolean algebra.
\end{proof}

Let $\mathcal M^{\nabla_2}$\index{$\mathcal M^{\nabla_2}$ \quad variety of monadic De Morgan algebras satisfying \(\nabla(x \meet {\hneg x}) \leq {\hneg\nabla(x \meet {\hneg x})}\)} denote the variety of monadic De Morgan algebras that satisfy the equation equivalent to \(\nabla(x \meet {\hneg x}) \leq {\hneg\nabla(x \meet {\hneg x})}\), and let $\mathcal K^{\nabla_1}$\index{$\mathcal K^{\nabla_1}$ \quad variety of monadic Kleene algebras satisfying \(\nabla(x \meet {\hneg x}) \leq {\hneg\nabla(x \meet {\hneg x})}\)} denote the variety of monadic Kleene algebras that satisfy \(\nabla(x \meet {\hneg x}) \leq {\hneg\nabla(x \meet {\hneg x})}\). %Then we can expand our lattice of subvarieties of $\mathcal M^\nabla$ somewhat:

\begin{proposition} % sublattice of varieties
\label{sublattice of varieties}
The following diagram is a sublattice of the lattice of subvarieties of $\mathcal M^\nabla$.
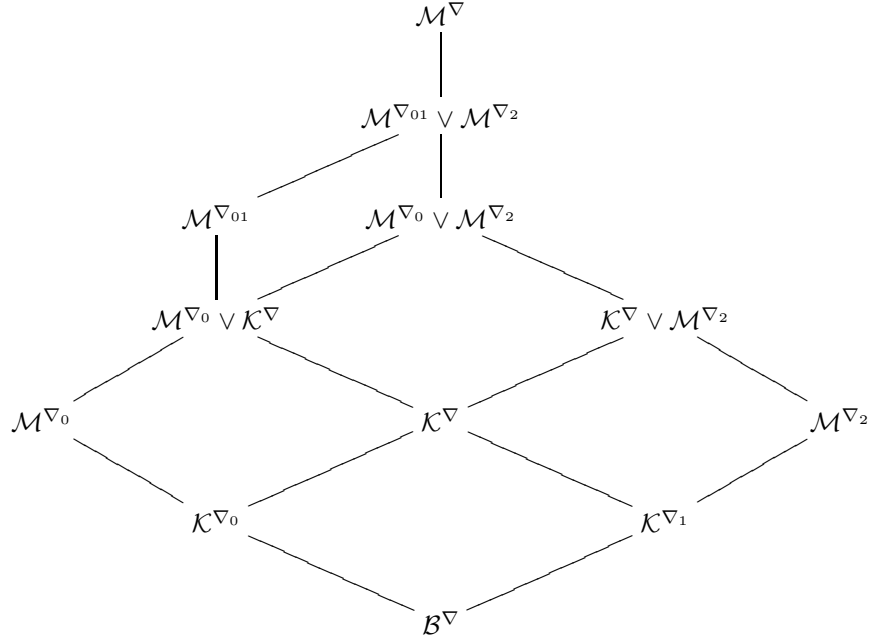
\[
\xymatrix{
&																					& {\mathcal M^\nabla} \ar@{-}[d] \\
&																					& {\mathcal M^{\nabla_{01}} \join \mathcal M^{\nabla_2}} \ar@{-}[dl] \ar@{-}[d]	 \\
&{\mathcal M^{\nabla_{01}}} \ar@{-}[d] 											& {\mathcal M^{\nabla_{0}} \join \mathcal M^{\nabla_2}} \ar@{-}[dl] \ar@{-}[dr] \\
&{\mathcal M^{\nabla_{0}} \join \mathcal K^\nabla} \ar@{-}[dl] \ar@{-}[dr]		& 					 																										& {\mathcal K^\nabla \join \mathcal M^{\nabla_2}} \ar@{-}[dl] \ar@{-}[dr] \\
{\mathcal M^{\nabla_0}} \ar@{-}[dr]		&										& {\mathcal K^{\nabla}} \ar@{-}[dl] \ar@{-}[dr]																			&& {\mathcal M^{\nabla_2}} \ar@{-}[dl] \\
&{\mathcal K^{\nabla_0}} \ar@{-}[dr]												&																															& {\mathcal K^{\nabla_1}} \ar@{-}[dl]	\\
&																					& {\mathcal B^{\nabla}}																									&			
}
\]
\end{proposition}

\begin{proof}
Starting at the bottom, \(\mathcal K^{\nabla_0} \meet \mathcal K^{\nabla_1} = \mathcal B^\nabla\) by \corref{subdirectly irreducible M-nabla-2}, and \(\mathcal K^{\nabla_0} \join \mathcal K^{\nabla_1} = \mathcal K^\nabla\) because if $\tuple{\A, \nabla}$ is a subdirectly irreducible algebra in $\mathcal K^\nabla$, then either $\nabla$ is of type 0 or 1. In the first case \(\tuple{\A, \nabla} \in \mathcal K^{\nabla_0}\). In the second case, let $a$ be the fixed point of $\A$. Then for any \(x \in \A\), \(x \meet {\hneg x} \leq a\), so \(\nabla(x \meet {\hneg x}) \leq a \leq {\hneg\nabla(x \meet {\hneg x})}\). Hence \(\tuple{\A, \nabla} \in \mathcal K^{\nabla_1}\). 

Moving up and to the right, \(\mathcal K^\nabla \meet \mathcal M^{\nabla_2} = \mathcal K^{\nabla_1}\) by definition. The monadic Kleene algebra \(\tuple{\mathbf{\underline K}, \nabla_0}\) belongs to \(\mathcal K^\nabla \setminus \mathcal M^{\nabla_2}\), while \(\tuple{\mathbf{\underline M}, \nabla_2}\) belongs to \(\mathcal M^{\nabla_2} \setminus K^\nabla\), so neither \(\mathcal K^\nabla \subseteq \mathcal M^{\nabla_2}\) nor \(\mathcal M^{\nabla_2} \subseteq \mathcal K^\nabla\).

Moving left, \(\mathcal M^{\nabla_0} \meet \mathcal K^\nabla = \mathcal K^{\nabla_0}\) by definition. Also, \(\tuple{\mathbf{\underline M}, \nabla_0} \in \mathcal M^{\nabla_0} \setminus \mathcal K^\nabla\), and \(\tuple{\mathbf{\underline K}, \nabla_1} \in \mathcal K^\nabla \setminus \mathcal M^{\nabla_0}\). However, \(\mathcal M^{\nabla_0} \join \mathcal K^\nabla \subset \mathcal M^{\nabla_{01}}\) by \propref{M-nabla-0 join K-nabla < M-nabla-0}.

Before moving right it will be useful to note that by \corref{subdirectly irreducible M-nabla-2}, \(\mathcal M^{\nabla_0} \meet \mathcal M^{\nabla_2} = \mathcal B^\nabla\) and \(\mathcal M^{\nabla_{01}} \meet \mathcal M^{\nabla_2} = \mathcal K^{\nabla_1}\). Also, \(\tuple{\mathbf{\underline M}, \nabla_0} \in \mathcal M^{\nabla_0} \setminus \mathcal M^{\nabla_2}\), while \(\tuple{\mathbf{\underline M}, \nabla_2} \in \mathcal M^{\nabla_2} \setminus \mathcal M^{\nabla_{01}}\). Thus 
\begin{align*}
(\mathcal M^{\nabla_0} \join \mathcal K^\nabla) \join (\mathcal K^\nabla \join \mathcal M^{\nabla_2}) &= \mathcal M^{\nabla_0} \join \mathcal K^\nabla \join \mathcal M^{\nabla_2} \\
&= \mathcal M^{\nabla_0} \join (\mathcal K^{\nabla_0} \join \mathcal K^{\nabla_1}) \join \mathcal M^{\nabla_2} \\
&= (\mathcal M^{\nabla_0} \join \mathcal K^{\nabla_0}) \join (\mathcal K^{\nabla_1} \join \mathcal M^{\nabla_2}) \\
&= \mathcal M^{\nabla_0} \join \mathcal M^{\nabla_2}.
\end{align*}
By \corref{distributive varieties}, 
\begin{align*}
\mathcal M^{\nabla_{01}} \meet (\mathcal M^{\nabla_0} \join \mathcal M^{\nabla_2}) 
&= (\mathcal M^{\nabla_{01}} \meet \mathcal M^{\nabla_0}) \join (\mathcal M^{\nabla_{01}} \meet \mathcal M^{\nabla_2}) \\
&= \mathcal M^{\nabla_0} \join \mathcal K^{\nabla_1} \\
&= \mathcal M^{\nabla_0} \join \mathcal K^\nabla.
\end{align*}
Also,
\(\mathcal M^{\nabla_{01}} \join (\mathcal M^{\nabla_0} \join \mathcal M^{\nabla_2}) = (\mathcal M^{\nabla_{01}} \join \mathcal M^{\nabla_0}) \join \mathcal M^{\nabla_2} = \mathcal M^{\nabla_{01}} \join \mathcal M^{\nabla_2}\).

Finally to show \(\mathcal M^{\nabla_{01}} \join \mathcal M^{\nabla_2} \subset \mathcal M^\nabla\) let $a$ and $b$ be the fixed points of $\mathbf{\underline M}$, and let $\A$ be the following subalgebra of \(\mathbf{\underline M} \times \mathbf{\underline M}\):
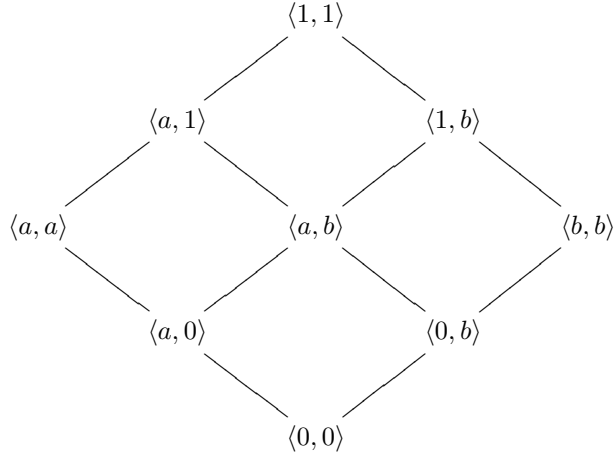
\[
\xymatrix{
							&											& \tuple{1,1} \ar@{-}[dl] \ar@{-}[dr] \\
							& \tuple{a,1} \ar@{-}[dl] \ar@{-}[dr]		&											& \tuple{1,b} \ar@{-}[dl] \ar@{-}[dr] \\
\tuple{a,a}	 \ar@{-}[dr]	&											&\tuple{a,b}  \ar@{-}[dl] \ar@{-}[dr]		&											&\tuple{b,b} \ar@{-}[dl] \\
							& \tuple{a,0}	 \ar@{-}[dr]				&											& \tuple{0,b} 	 \ar@{-}[dl] \\
							&											&\tuple{0,0}
}
\]
Observe that $\tuple{a,a}$ and $\tuple{b,b}$ are fixed points in $\A$ such that \(\tuple{a,a} \meet \tuple{b,b} = \tuple{0,0}\) and \(\tuple{a,a} \join \tuple{b,b} = \tuple{1,1}\). Let $\nabla$ be the type 2 quantifier whose range is 
\[
\set{\tuple{0,0}, \tuple{a,a}, \tuple{b,b}, \tuple{1,1}}.
\]
Then \(\tuple{\A, \nabla}\) is a monadic De Morgan algebra that does not belong to \(\mathcal M^{\nabla_{01}} \cup \mathcal M^{\nabla_2}\), on the one hand because $\nabla$ is a quantifier of type 2, and on the other hand because $\tuple{a,b}$ is a fixed point that is not fixed by $\nabla$. In order to apply \lemref{V1 join V2 -> V1 cup V2} we need to show that $\tuple{\A, \nabla}$ is subdirectly irreducible. In fact, $\tuple{\A, \nabla}$ is simple. We perform three representative calculations. Let $\cong$ be a congruence on $\tuple{\A, \nabla}$. If \(\tuple{0,0} \cong \tuple{a,0}\), then \(\tuple{0,0} = \nabla\tuple{0,0} \cong \nabla\tuple{a,0} = \nabla\tuple{a,a} = \tuple{a,a}\), in which case \(\tuple{0,0} \cong \tuple{a,a} = {\hneg\tuple{a,a}} \cong \tuple{1,1}\). Hence $\cong$ is the total congruence. Similarly, if \(\tuple{0,0} \cong \tuple{0,b}\), then $\cong$ is the total congruence. If \(\tuple{a,0} \cong \tuple{a,a}\), then \(\tuple{a,b} \cong \tuple{a,1}\), and \(\tuple{1,b} \cong \tuple{1,1}\). Hence \(\tuple{0,0} = {\hneg\tuple{1,1}} \cong {\hneg\tuple{1,b}} = \tuple{0,b}\), so $\cong$ is the total congruence. If \(\tuple{a,a} \cong \tuple{a,b}\), then \(\tuple{a,a} = \nabla\tuple{a,a} \cong \nabla\tuple{a,b} = \tuple{1,1}\), so $\cong$ is the total congruence. Thus $\tuple{\A, \nabla}$ is simple. Therefore $\tuple{\A, \nabla}$ does not belong to $\mathcal M^{\nabla_{01}} \join \mathcal M^{\nabla_2}$.
\end{proof}

%% The monadic De Morgan reduct
\subsection{The monadic De Morgan reduct}

Let $X$ be a one-dimensional pair of suits. It follows from \corref{C_0...C_N-1 pair of suits} that 
\begin{align*}
C_{0,\set{0}}(X) &= 
\begin{cases}
1      		& \text{if \(X \not\leq \Omega\) and \(X \not\leq \mho\)}, \\
\Omega		& \text{if \(0 < X \leq \Omega\)}, \\
\mho		& \text{if \(0 < X \leq \mho\)}, \\
0	    	& \text{if \(X = 0\)}.
\end{cases}
\intertext{If $X$ is a double suit, then }
C_{0,\set{0}}(X) &= 
\begin{cases}
1	    	& \text{if \(X \not\leq \Omega\)}, \\
\Omega		& \text{if \(0 < X \leq \Omega\)}, \\
0      		& \text{if \(X = 0\)}.
\end{cases}
\end{align*}
Thus the reduct of a suited IFG$_1$-algebra that includes $\Omega$ and $\mho$ to the signature \(\tuple{0,1,\n{},+_{\set{0}},\cdot_{\set{0}}, C_{0,\set{0}}}\) is a monadic De Morgan algebra with a quantifier of type 2, while the same reduct of a double-suited IFG$_1$-algebra that includes $\Omega$ is a monadic Kleene algebra with a quantifier of type 1. We will refer to the reduct of a suited IFG$_1$-algebra $\C$ to the signature \[\tuple{0,1,\n{},+_{\set{0}},\cdot_{\set{0}}, C_{0,\set{0}}}\] as the \emph{monadic De Morgan reduct}\index{monadic De Morgan reduct} of $\C$. It follows that the variety generated by the monadic De Morgan reducts of all suited IFG$_1$-algebras is contained in $\mathcal M^\nabla$.

\begin{proposition} % C(X * -X) < -C(X * -X)
\label{C(X * -X) < -C(X * -X)}
If $X$ is a double suit, then 
\[
C_{0,J_0} \ldots C_{N-1,J_{N-1}}(X \cdot_N \n{X}) \leq \n{C_{0,J_0} \ldots C_{N-1,J_{N-1}}(X \cdot_N \n{X})}.
\]
\end{proposition}

\begin{proof}
If $X$ is a double suit, then \(0 \leq X \cdot_N \n{X} \leq \Omega\). If \(X \cdot_N \n{X} = 0\), then 
\[
C_{0,J_0} \ldots C_{N-1,J_{N-1}}(X \cdot_N \n{X}) = 0 
\quad
\text{ and }
\quad
\n{C_{0,\emptyset} \ldots C_{N-1,\emptyset}(X \cdot_N \n{X})} = 1.
\]
If \(0 < X \cdot_N \n{X} \leq \Omega\), then 
\[
C_{0,J_0} \ldots C_{N-1,J_{N-1}}(X \cdot_N \n{X}) = \Omega = \n{C_{0,J_0} \ldots C_{N-1,J_{N-1}}(X \cdot_N \n{X})}. \qedhere
\]
\end{proof}

In particular, in the one-dimensional case \(C_{0,\set{0}}(X \cdot_{\set{0}} \n{X}) \leq \n{C_{0,\set{0}}(X \cdot_{\set{0}} \n{X})}\). Thus the variety generated by the monadic De Morgan reducts of all double-suited IFG$_1$-algebras is contained in $\mathcal K^{\nabla_1}\).

\begin{conjecture} % V(K) = K^\nabla_1
The variety generated by the monadic De Morgan reducts of all double-suited IFG$_1$-algebras is $\mathcal K^{\nabla_1}\).
\end{conjecture}

To prove the conjecture it would be sufficient to show that every monadic Kleene algebra in $\mathcal K^{\nabla_1}$ is a homomorphic image of a subalgebra of a product of De Morgan reducts of double-suited IFG$_1$-algebras. By \corref{subdirectly irreducible M-nabla-2} it suffices to consider monadic Kleene algebras with type 1 quantifiers and monadic Boolean algebras with type 0 quantifiers.

\begin{definition} % meet/join irreducible
Let $L$ be a lattice. An element \(a \in L\) is \emph{join irreducible}\index{join irreducible|mainidx} if \(a = x \join y\) implies \(a = x\) or \(a = y\). Dually, \(b \in L\) is \emph{meet irreducible}\index{meet irreducible|mainidx} if \(b = x \meet y\) implies \(b = x\) or \(b = y\). 
\end{definition}

\begin{proposition} % suited 0 1 irreducible
\label{suited 0 1 irreducible}
In any De Morgan reduct of a suited IFG$_N$-algebra 0 is meet irreducible and 1 is join irreducible.
\end{proposition}

\begin{proof}
Let $\C$ be a the De Morgan reduct of a double-suited IFG$_N$-algebra, and let \(X,Y \in \C\). If \(X +_N Y = 1\), then \(\^NA \in (X +_N Y)^+ = X^+ \cup Y^+\). If \(\^NA \in X^+\) then \(X = 1\), and if \(\^NA \in Y^+\) then \(Y = 1\).
\end{proof}

It follows that not every Kleene algebra is isomorphic to a double-suited IFG$_1$-algebra. For example, let $K$ be the following Kleene algebra with fixed point $c$:
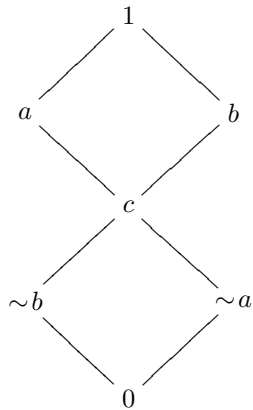
\[
\xymatrix{
					& 1 \ar@{-}[dl] \ar@{-}[dr]	&					\\
	a \ar@{-}[dr] 		& 					& b \ar@{-}[dl]			\\
					& c \ar@{-}[dl] \ar@{-}[dr] &				\\
	{\hneg b} \ar@{-}[dr]	&					& {\hneg a} \ar@{-}[dl]		\\
					& 0 		&					\\
}
\]
Observe that 0 is not meet irreducible, and 1 is not join irreducible, so $K$ cannot be embedded into any suited IFG$_N$-algebra. Perhaps every monadic Kleene algebra with a type 1 quantifier in which 0 is meet irreducible and 1 is join irreducible is isomorphic to the De Morgan reduct of a double-suited IFG$_1$-algebra. However, my advisor J.~Donald Monk and I have only been able to prove the following partial result.

\begin{theorem} % embedding monadic Kleene algebras with irreducible bounds
\label{embedding monadic Kleene algebras with irreducible bounds}
Every monadic Kleene algebra with a quantifier of type 1 in which  0 is meet irreducible and 1 is join irreducible is isomorphic to the De Morgan reduct of a rooted IFG$_1$-algebra.
\end{theorem}

\begin{proof}
Let $\tuple{K, \nabla}$ be a monadic Kleene algebra with a quantifier of type 1 in which 0 is meet irreducible and 1 is join irreducible. Let $c$ be the fixed point of $K$, and consider the interval \(L = [c, 1]\) as a bounded distributive lattice with minimum $c$ and maximum 1. Let $A$ be the set of prime filters on $L$, and let $\sigma$ be the Priestley isomorphism \(x \mapsto \setof{F \in A}{x \in F}\). Since 1 is join irreducible we have that $\set{1}$ is a prime filter, and for all \(x \in L\), \(\set{1} \in \sigma(x)\) if and only if \(x = 1\).

There is a partition $P$ of $\powerset(A)$ such that \(\abs{P} = \abs{A}\), every cell of the partition includes a singleton, but no two singletons belong to the same cell. Let \(f\colon A \to P\) be a bijection such that for all \(a \in A\), \(\set{a} \in f(a)\), and \(A \in f(\set{1})\). Define \(g\colon \powerset(A) \to \powerset(\powerset(A)) \setminus \set{\emptyset}\) by 
\[
g(U) = \set{\emptyset} \cup \bigcup_{a\in U} f(a).
\]
If \(U \not= V\), then without loss of generality there is a \(b \in U \setminus V\). Consequently \(f(b) \subseteq g(U) \setminus g(V)\). Since $f(b)$ is nonempty, \(g(U) \not= g(V)\). Hence $g$ is injective. Note that there is a singleton \(\set{a} \in g(U)\) if and only if $U$ is nonempty, and \(A \in g(U)\) if and only if \(\set{1} \in U\). 

Let \(G = g \circ \sigma\). Note that 
\[
G(c) = g(\sigma(c)) = g(\emptyset) = \set{\emptyset},
\]
\[
G(1) = g(\sigma(1)) = g(A) = \powerset(A),
\]
\[
G(x \join y) = g(\sigma(x \join y)) = g(\sigma(x) \cup \sigma(y)) = g(\sigma(x)) \cup g(\sigma(y)) = G(x) \cup G(y).
\]
To prove \(G(x \meet y) = G(x) \cap G(y)\) we need to show that 
\[
\set{\emptyset}\cup \hspace{-16 pt} \bigcup_{a \in \sigma(x) \cap \sigma(y)} \hspace{-16 pt} f(a)
= \left(\set{\emptyset} \cup \hspace{-5 pt}\bigcup_{a \in \sigma(x)} \hspace{-5 pt} f(a)\right) \cap 
\left(\set{\emptyset} \cup \hspace{-5 pt} \bigcup_{b \in \sigma(y)} \hspace{-5 pt} f(b)\right).
\]
Suppose $u$ is a nonempty member of the left-hand side. Then for some \(a \in \sigma(x) \cap \sigma(y)\), \(u \in f(a)\). That is, there is a prime filter $a$ such that \(x,y \in a\) and \(u \in f(a)\), which is enough to show that $u$ belongs to the right-hand side. Conversely, suppose $v$ is a nonempty member of the right-hand side. Then for some \(a \in \sigma(x)\) and some \(b \in \sigma(y)\), \(v \in f(a)\) and \(v \in f(b)\). However, if $a$ and $b$ were distinct, $f(a)$ and $f(b)$ would be disjoint. Thus \(a = b\), and $v$ belongs to the left-hand side.

The function $G$ is injective because $g$ and $\sigma$ are. Also note that there is a singleton \(\set{a} \in G(x)\) whenever \(x > c\), and \(A \in G(x)\) if and only if \(\set{1} \in \sigma(x)\) if and only if \(x = 1\).

Define a function $h$ from $\tuple{K,\nabla}$ to the monadic De Morgan reduct of $\Root_1(A)$ by \(h(x) = \tuple{G(x \join c),\, G({\hneg x} \join c)}\). If \(x \not= y\), then either \(x \join c \not= y \join c\) or \({\hneg x} \join c \not= {\hneg y} \join c\) because in a distributive lattice \(x \join c = y \join c\) and \(x \meet c = y \meet c\) imply \(x = y\). Thus $h$ is injective. 

To show that $h$ is a homomorphism, observe that 
\[
h(0) 	= \tuple{G(0 \join c),\, G(1 \join c)} \\
		= \tuple{G(c),\, G(1)} \\
		= \tuple{\set{\emptyset},\, \powerset(A)} = 0,
\]
\[
h(1) 	= \tuple{G(1 \join c),\, G(0 \join c)} \\
		= \tuple{G(1),\, G(c)} \\
		= \tuple{\powerset(A),\, \set{\emptyset}} = 1,
\]
\[
h({\hneg x}) = \tuple{G({\hneg x} \join c),\, G(x \join c)} = \n{h(x)},
\]
\begin{align*}
h(x \join y)	
	&= \tuple{G((x \join y) \join c),\, G({\hneg(x \join y)} \join c)} \\
	&= \tuple{G((x \join c) \join (y \join c)),\, G(({\hneg x} \join c) \meet ({\hneg y} \join c))} \\
	&= \tuple{G(x \join c) \cup G(y \join c),\, G({\hneg x} \join c) \cap G({\hneg y} \join c)} \\
	&= h(x) +_{\set{0}} h(y).
\end{align*}
Now we check the quantifier:
\[
h(\nabla 0) = h(0) = 0 = C_{0, \set{0}}(0) = C_{0,\set{0}}(h(0)).
\]
If \(0 < x \leq c\), then \(x \join c = c\) and \({\hneg x} \join c < 1\). Hence \(G(x \join c) = \set{\emptyset}\) and \(A \notin G({\hneg x} \join c)\). Thus
\[
h(\nabla x) = h(c) = \tuple{G(c),\, G(c)} = \Omega = C_{0,\set{0}}(\tuple{G(x \join c),\, G({\hneg x} \join c)}) = C_{0,\set{0}}(h(x)).
\]
If \(x \not\leq c\), then \(x \join c > c\) and \({\hneg x} \join c < 1\). Hence there is a singleton \(\set{a} \in G(x \join c)\) and \(A \notin G({\hneg x} \join c)\). Thus
\[
h(\nabla x) = h(1) = 1 = C_{0,\set{0}}(\tuple{G(x \join c),\, G({\hneg x} \join c)}) = C_{0,\set{0}}(h(x)).
\]
by \lemref{C_0...C_N-1 rooted}. Therefore $h$ is an embedding.
\end{proof}

Note that \thmref{embedding monadic Kleene algebras with irreducible bounds} does not resolve the conjecture because the elements in the range of $h$ are not double suits (or even pairs of suits).

%%%%%%%%%%%%%%%%%%%%%%%%%%%%%%%%%%%%%%%%%%%%%
	
%% Bibliography
\bibliographystyle{plain}

\end{document}